\documentclass[a4paper,10pt]{article}
\usepackage{amsmath}
\usepackage{amsfonts}
\usepackage{amssymb}
\usepackage{amsthm}
\usepackage{enumerate}
\usepackage{bbm}
\usepackage{hyperref}
\usepackage{color}
\usepackage{graphicx}
\usepackage[section]{placeins}
\graphicspath{ {images/} }
\usepackage[abs]{overpic}
\usepackage{pgfplots}
\usepackage{booktabs}
\usepackage[utf8]{inputenc}
\usepackage{lmodern}
\usepackage[T1]{fontenc}

\usepackage{float}

\usepackage{pgfplots}
\usepackage{caption,subcaption}

\usepackage{mathtools}
\usepackage[hang,flushmargin]{footmisc}

\usepackage{a4wide}

\newtheorem{theorem}{Theorem}[section]

\newtheorem{lemma}[theorem]{Lemma}

\theoremstyle{definition}

\theoremstyle{remark}
\newtheorem{remark}[theorem]{Remark}

\newcommand{\R}{\mathbb{R}}
\newcommand{\C}{\mathbb{C}}

\newcommand{\Z}{\mathbb{Z}}
\newcommand{\ii}{\mathrm{i}}
\newcommand{\ee}{\mathrm{e}}
\newcommand{\dd}{\mathrm{d}}

\newcommand{\LL}{\mathrm{L}}

\newcommand{\cw}{c_{n,\omega}}

\newcommand{\HH}{\mathrm{H}}

\newenvironment{Eqnarray}
{\arraycolsep=1.4pt
  \begin{eqnarray}}
  {\end{eqnarray}
    \hspace*{-4pt}}

\newenvironment{Eqnarray*}
{\arraycolsep=1.4pt
  \begin{eqnarray*}}
  {\end{eqnarray*}
    \hspace*{-4pt}}

\providecommand{\keywords}[1]
{
  \small	
  \textbf{Keywords and phrases} #1
}

\providecommand{\msc}[1]
{
  \small	
  \textbf{Mathematics subject classifications } #1
}

\numberwithin{equation}{section}

\def\email#1{email: \url{#1}}

\allowdisplaybreaks

\title{Approximation of Wave Packets on the Real Line}
\author{Arieh Iserles \footnote{Department of Applied Mathematics and Theoretical Physics, Centre for Mathematical Sciences, University of Cambridge, Wilberforce Rd, Cambridge CB4 1LE, United Kingdom, \email{ a.iserles@damtp.cam.ac.uk}\vspace*{0.25em}}
\and Karen Luong \footnote{Department of Applied Mathematics and Theoretical Physics, Centre for Mathematical Sciences, University of Cambridge, Wilberforce Rd, Cambridge CB4 1LE, United Kingdom, \email{ k.luong@damtp.cam.ac.uk}\vspace*{0.25em}}
\and Marcus Webb \footnote{Department of Mathematics, University of Manchester, Alan Turing Building, Manchester M13
9PL, United Kingdom, \email{marcus.webb@manchester.ac.uk}}}

\begin{document}
\maketitle
\begin{abstract}
In this paper we compare three different orthogonal systems in $\LL_2(\R)$ which can be used in the construction of a spectral method for solving the semi-classically scaled time dependent Schr\"odinger equation on the real line, specifically, stretched Fourier functions, Hermite functions and Malmquist--Takenaka functions. All three have banded skew-Hermitian differentiation matrices, which greatly simplifies their implementation in a spectral method, while ensuring that the numerical solution is unitary -- this is essential in order to respect the Born interpretation in quantum mechanics and, as a byproduct, ensures numerical stability with respect to the $\LL_2(\R)$ norm. We derive asymptotic approximations of the coefficients for a wave packet in each of these bases, which are extremely accurate in the high frequency regime. We show that the Malmquist--Takenaka basis is superior, in a practical sense, to the more commonly used Hermite functions and stretched Fourier expansions for approximating wave packets. 
\end{abstract}


\msc{30E10, 41A60, 42C05, 41A20, 42A16}

\section{Introduction \label{sec:intro}}

Highly oscillatory wave packets play an important role in quantum mechanics. A case in point is the semi-classically scaled time dependent Schr\"odinger equation (TDSE), which describes the quantum dynamics of the nuclei in a molecule \cite{lasser2020computing}. It has the form
\begin{Eqnarray}
        \ii \varepsilon \frac{\partial\psi(x,t)}{\partial t} = \left[-\frac{\varepsilon^{2}}{2} \Delta + V(x)\right]\psi(x,t), \label{eq:SE}
\end{Eqnarray}
where $\psi(x,t) \in \C$, $x \in \R^{d}$ and $t \in \R$. The Laplacian is $\Delta = \sum_{k=1}^d \partial^2/\partial x_k^2$, while $V$ is a smooth potential, $V(x):\R^{d} \rightarrow \R$. 
    
    This equation is important in a wide range of fields such as physical chemistry. However, solving TDSE comes with a number of computational challenges. A major reason is the small parameter $\varepsilon \ll 1$, where $\varepsilon^{2}$ describes the mass ratio of an electron and the nuclei. Consider equation \eqref{eq:SE} without $V(x)$ (the so-called {\em free Schr\"odinger equation\/}), then it can be shown that the  solution to this problem generates oscillations of frequency $\mathcal{O}\left(\varepsilon^{-1} \right)$ in space and time \cite{jin_markowich_sparber_2011}. Moreover, the current techniques used to solve equation \eqref{eq:SE} are mainly defined on a finite interval. In that instance Fourier spectral and pseudo-spectral methods work well in practice \cite{feit1982solution}. However, these approaches require the imposition of periodic boundary conditions, which may lead to non-physical behaviour over longer time frames, because the solution inevitably reaches the non-physical periodic boundary of the domain. Modern applications like quantum control also require much longer time frames. One solution to this problem, is to solve TDSE on the whole real line using spectral methods and this underlies the setting of this paper.
    
In designing a spectral method, we first need to pick basis functions that approximate the solution to equation \eqref{eq:SE}. In this paper, we restrict our focus to systems that are orthonormal in $\LL_{2}(\R)$ and have banded skew-Hermitian differentiation matrices. When the basis has a skew-Hermitian differentiation matrix, the solution operator is naturally unitary, so the Born interpretation of quantum mechanics is respected. As a result, these methods are also \emph{stable by design} \cite{hairer2016numerical,iserles2014skew}.  

Heller in 1976 proposed that the solutions to the semi-classical TDSE can be approximated by Gaussian type functions i.e.~wave packets \cite{heller1976time}.  Hagedorn proved error bounds for time dependent wave packets in 1980 \cite{hagedorn1980semiclassical}.  Lubich in 2008, proved that for travelling wave packets with fixed width $\sim \varepsilon^{1/2}$ that a representation of the solution by wave packets has an $\LL_{2}$ error of $\mathcal{O}\left(\varepsilon^{1/2}\right)$ for time intervals of size $t \sim 1$ \cite{lubich2008quantum}. Thus, wave packets are a good first approximation to the solution to semi-classical TDSE. In fact, for a quadratic potential with Gaussian initial data, they yield an exact solution to equation \eqref{eq:SE} \cite{lubich2008quantum}. Thus, a major yardstick for the suitability of basis functions for equation \eqref{eq:SE} is how well they approximate wave packets. 

A univariate wave packet function has the form $f(z)=\ee^{-\alpha(z-x_{0})^{2}}\cos(\omega z)$, where $\alpha>0$, $x_{0} \in \R$ and $\omega\in\R$.  The width of the wave packet is controlled by $\alpha$, $x_{0}$ is the position where the wave packet is centred and $\omega$ is the frequency of oscillations. For simplicity much of our analysis is restricted to $x_{0}, \omega >0$ and this should be assumed unless stated otherwise.
The expansion coefficients of a wave packet function, $f(x)$, in an orthonormal set of $\LL_{2}(\R)$ basis functions, $\{\varphi_n\}_{n \in I} $ where $I$ is a countable index set such as $\Z$, is given by the following inner product,
\begin{Eqnarray}
  a_{n}(\omega) &=&\int^{\infty}_{-\infty} \ee^{-\alpha(x-x_{0})^{2}}\cos(\omega x)\overline{ \varphi_{n}(x)} \dd x, \label{eq:generala} \\
  &=&\frac{1}{2} \int^{\infty}_{-\infty} \left( \ee^{\ii \omega x} + \ee^{-\ii \omega x}\right)  \ee^{-\alpha(x-x_{0})^{2}}\overline{ \varphi_{n}(x)} \dd x=b_{n}(\omega) + b_{n}(-\omega),\nonumber 
\end{Eqnarray}
where 
\begin{Eqnarray*}
	b_{n}(\omega)&:=&\frac{1}{2} \int_{-\infty}^{\infty} \ee^{-\alpha(x-x_{0})^{2}+\ii \omega x } \overline{\varphi_{n}(x)} \dd x.
\end{Eqnarray*}
Much of our analysis concerns the asymptotics of $b_n(\omega)$ in the large-$\omega$ regime, from which results about $a_n(\omega)$ follow. Note we use the notation $f(x) \sim g(x)$ for $x \rightarrow \infty$ to denote that functions $f(x)$ and $g(x)$ are asymptotically equivalent as $x \rightarrow \infty$. This can be rewritten as $f(x) = g(x)[ 1 + o(1)]$ \cite{de1981asymptotic}.

The three systems under consideration in this paper are 
\begin{itemize}
	\item stretched Fourier functions, which are simply the Fourier functions scaled to an interval $[-\lambda, \lambda]$,
	\begin{equation*}
\varphi_{n} (x) := \begin{cases}  (2\lambda)^{-1/2}\ee^{\ii \pi n  x/\lambda}&\mbox{if } x \in [-\lambda,\lambda], \\
0 & \mbox{for }x \in \mathbb{R}\setminus[-\lambda,\lambda]. \end{cases}
\end{equation*}
	\item Hermite functions  \cite[Table 18.3.1]{dlmf}
	\begin{equation*}
	\varphi_{n}(x) := \frac{1}{(2^{n}n!\pi^{1/2})^{1/2}} \HH_n(x)\ee^{-x^2/2}, \qquad n \in \mathbb{Z}_+.
\end{equation*}
	\item Malmquist--Takenaka (MT) functions
	\begin{equation*}
  \varphi_{n}(x) := \ii^{n}\sqrt{\frac{2}{\pi}}\frac{(1 + 2\ii x)^{n}}{(1 - 2\ii x )^{n+1}}, \qquad \qquad \qquad n \in \Z. 
\end{equation*}
\end{itemize}
In this paper, we derive an asymptotic estimate of the coefficients in each basis, which will allow us to compare and evaluate how well they approximate wave packets on the whole real line. For simplicity and without loss of generality  we restrict our analysis to $\alpha, x_{0}, \omega, n >0$.

For stretched Fourier functions, the asymptotic estimate of coefficients is described in Theorem \ref{thm:sf}. In its more streamlined version, for any given $\varepsilon > 0$ it is possible to choose the stretching constant $\lambda$ such that
\begin{equation*}
|a_{n,\lambda}(\omega)| \leq \left(\frac{\pi}{2\alpha\lambda}\right)^{1/2} \left[ \exp\!\left(\!-\frac{(\pi |n|-\lambda\omega)^2}{4\alpha\lambda^2}\right) + \varepsilon \right]\!.
\end{equation*}
A specific choice that satisfies this bound is
\begin{equation}
   \lambda = |x_{0}| + \sqrt{\frac{\log(\varepsilon^{-1})}{\alpha}}. \label{eq:sflambda}
\end{equation} 

\begin{figure}[tbh]
  \begin{center}
  	\includegraphics[width=.6\textwidth]{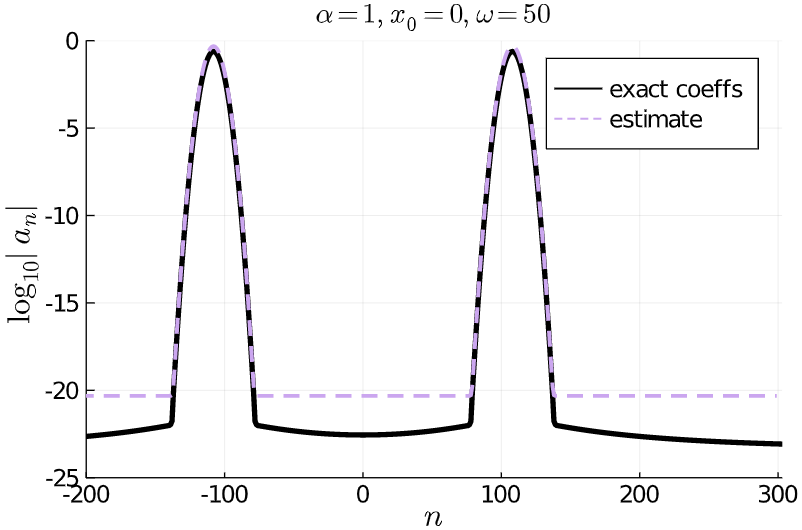}
   \end{center}
  \caption{\textit{Stretched Fourier functions:\/} Wave packet coefficients in a stretched Fourier basis. The solid line represents the exact coefficients and the (light) dashed line is a plot of the spectrally decaying part of our estimate in \eqref{eq:aSF} for $n>0$. We can also see the sombrero-like hump that is characteristic to stretched Fourier functions, which shows that our coefficients decay spectrally. The slow decay rate of coefficients (the brim of the sombrero) was purposefully chosen to start after $\varepsilon = 10^{-20}$ using \eqref{eq:sflambda}. \label{fig:sfcompare}}
\end{figure}

From this estimate we expect the logarithm of the magnitude of the expansion coefficients to resemble a hump, very similar to a sombrero, and we can observe these characteristics in Fig.~\ref{fig:sfcompare}. The sombrero-like hump tells us that the coefficients decay spectrally until the small constant $\varepsilon$ eventually dominates the estimate. This drop off in decay rate is real, and for $\alpha = 1, x_{0} =0$ and $\omega=50$ in Fig.~\ref{fig:sfcompare} occurs at about $\varepsilon=10^{-20}$, in line with equation \eqref{eq:sflambda}. We have a pessimistic certainty that for all $n$ satisfying,
\begin{displaymath}
 (\pi |n| - \lambda \omega)^2 \leq 4\alpha \lambda^{2} \log\!\left( \varepsilon^{-1}\left(2\left(2\alpha\lambda/\pi\right)^{1/2} - 1\right)^{-1} \right),
\end{displaymath}
we have  $|a_{n,\lambda}(\omega)| \geq 2\varepsilon$. The value of the stretching constant $\lambda$ is discussed further in section \ref{sec:sf}. One clear advantage of using the stretched Fourier basis is that we can take advantage of the structure seen in Fig.~\ref{fig:sfcompare}, where we can confine our calculation to approximating $N$ coefficients, the complement of which is guaranteed to satisfy $|a_{n}(\omega)| \leq \varepsilon$, by a single Fast Fourier Transform (FFT) in $\mathcal{O}(N\log_2 N)$ operations.

In summary, stretched Fourier functions hold a great deal of promise; we have proved that the coefficients exhibit spectral decay and can be approximated using FFT. However, in a more practical setting, our method of determining the optimal cut-off of coefficients  requires information about our wave packet which is typically unavailable. Moreover, a typical solution resembles an unknown linear combination of wave packets, whereas the nice performance e.g.~in Fig.~\ref{fig:sfcompare} is optimized for a single wave packet. This is discussed in more detail in section \ref{sec:sf}.

\begin{figure}[tbh]
  	\begin{center}
  		\includegraphics[width=.45\textwidth]{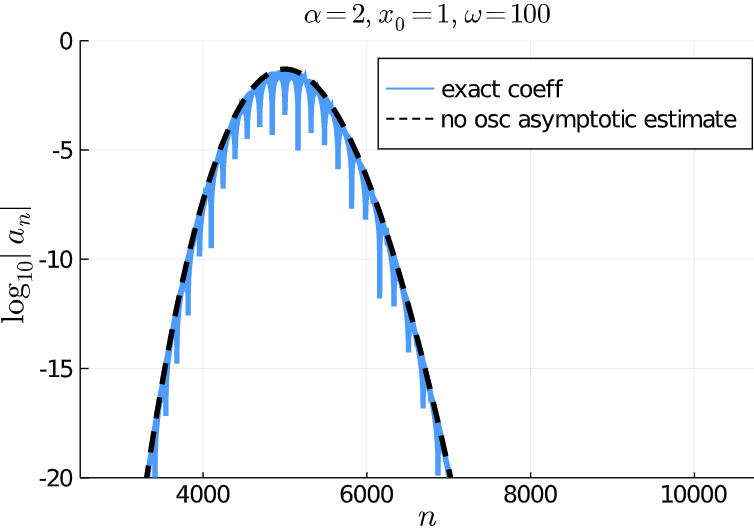}
  		\includegraphics[width=.45\textwidth]{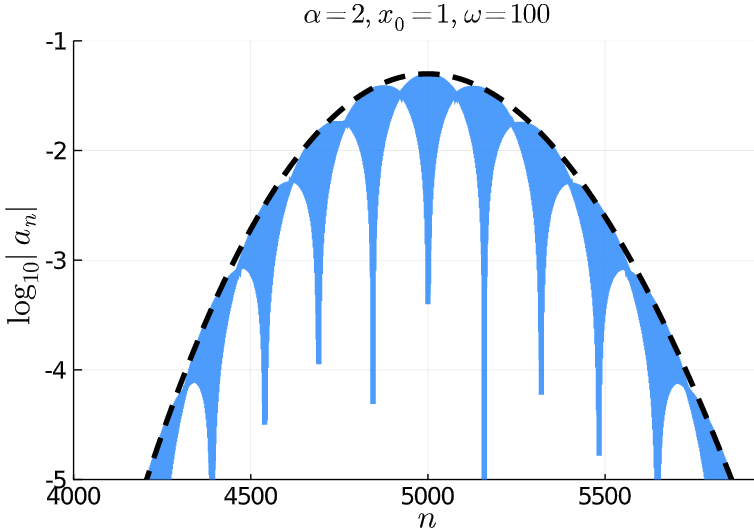}
  	\end{center}
  \caption{ \textit{Hermite functions:} The left figure shows the exact coefficients in the lighter colour and the asymptotic estimate with the bounded oscillatory term removed in the black dashed line. To further make the point that our asymptotic estimate captures the characteristics of the exact coefficients, the plot on the right is zoomed on  the coefficients near the peak.} 
  \label{fig:hermite1compare}
\end{figure}

 The next contender is the basis of Hermite functions \cite[Table 18.3.1]{dlmf}. Deferring the full asymptotic estimate to Theorem \ref{thm:hermitefunc}, a simplified version for $n > 0$ and $\alpha > 1/2$ is
			\begin{Eqnarray*}
			a_{n}(\omega) &\sim & \frac{\mathcal{C}_{\mathrm{osc}} \left[1 + \mathcal{O}\left(n^{-1}\right) \right] }{\sqrt{\omega(2\alpha +1)}} \left| \frac{1 + \sqrt{1 +2\frac{n}{\omega^{2}}\left(4\alpha^{2}-1\right)}}{\sqrt{2\frac{n}{\omega^{2}}}(2\alpha +1)}\right|^{n} \\
		& & \times \exp\!\left[-\frac{\omega^{2}}{2(2\alpha + 1)}+ \frac{n}{1+\sqrt{1+2\frac{n}{\omega^{2}}\left(4\alpha^{2}-1\right)}} +  \frac{\alpha x_{0}^{2}}{4\alpha^{2}-1} \left( 1 - \frac{2 \alpha}{\sqrt{1+2\frac{n}{\omega^{2}}\left(4\alpha^{2}-1\right)}}\right)\right]\!,
			\end{Eqnarray*}
			where $\mathcal{C}_{\mathrm{osc}}$ is a bounded oscillatory component. It is difficult to deduce the `shape' of these coefficients from this expression. However, to help convince the reader, let us set  $\cw = n/\omega^{2}$ to be constant -- we see in section \ref{sec:derivingasymphermitea12} that this is a reasonable assumption and note that this relationship keeps $n$ and $\omega$ large. Then we have
\begin{itemize}
			\item a factor of $1/\sqrt{\omega}$ which gently drives the size of $a_{n}(\omega)$ down.
			\item A factor, greater than one, to the power $n$ drives the size up geometrically. 
			\item The argument of the exponential is essentially dominated by the competing $\omega$ and $n$ term. If we replace $n = \cw \omega^{2}$ and note that for $\alpha > 1/2$ this term is negative, we obtain rapid decay to zero.
\end{itemize}
Thus, overall	we get geometric decay. The above expression matches the envelope of the $a_{n}$'s exceedingly well. In Fig.~\ref{fig:hermite1compare}, we can see that for large $\omega$ and small $\alpha$ and $x_{0}$  we can clearly observe spectral decay.
			
\begin{figure}[tbh]
   \begin{center}
  	\includegraphics[width=.6\textwidth]{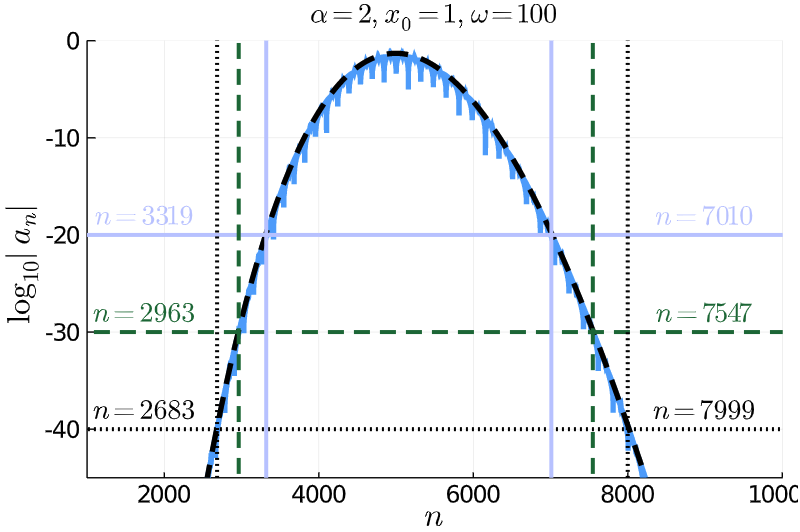}
  \end{center}
  \caption{ \textit{Hermite functions:\/} This is the same plot of the decay rate of coefficients as Fig.~\ref{fig:hermite1compare} but with a larger range. The plot shows that for varying values of $ \varepsilon = 10^{-20}, 10^{-30}, 10^{-40} $ the number of coefficients which are greater than $\varepsilon$ are 3691, 4584, and 5316 respectively. \label{fig:hermite4compare}}
\end{figure}

To understand how the number of coefficients grows with respect to $\omega$, consider the following crude inequality, the number of coefficients which are greater than $\varepsilon$ in absolute value are all the $n$'s that such that
\begin{Eqnarray*}
|n| &\gtrsim& \left[ \log \left| \varepsilon \right| +  \log \left| \frac{\sqrt{2 \alpha + 1}}{\mathcal{C}_{osc} \cw^{1/4}} \right|- \left| \frac{\alpha x_{0}^{2}}{4\alpha^{2}-1} \left( 1 - \frac{2 \alpha}{\sqrt{1+2\cw\left(4\alpha^{2}-1\right)}}\right) \right| \right]\!,
\end{Eqnarray*}
where $\cw$ is a constant, and $\alpha$ and $x_{0}$ are fixed. Note that number of coefficients increases linearly with $n$ but quadratically with $\omega$ (as $n=\cw \omega^{2}$). Fig.~\ref{fig:hermite4compare} displays the number of coefficients required for varying $\varepsilon$ and fixed $\omega$.
				
Although it is natural to think of Hermite functions in a quantum mechanical setting, since they are the eigenfunctions of the quantum harmonic oscillator, we show in Section \ref{sec:hermite}, that the convergence properties leave much to be desired. The decay of Hermite function coefficients for highly oscillatory wave packets with small $\alpha $ and $x_{0}$ is fast indeed, but as $\alpha$ and $x_0$ increase, this soon ceases to be true in practical terms.
			
\begin{figure}[tb]
  \begin{center}
  	\includegraphics[width=.45\textwidth]{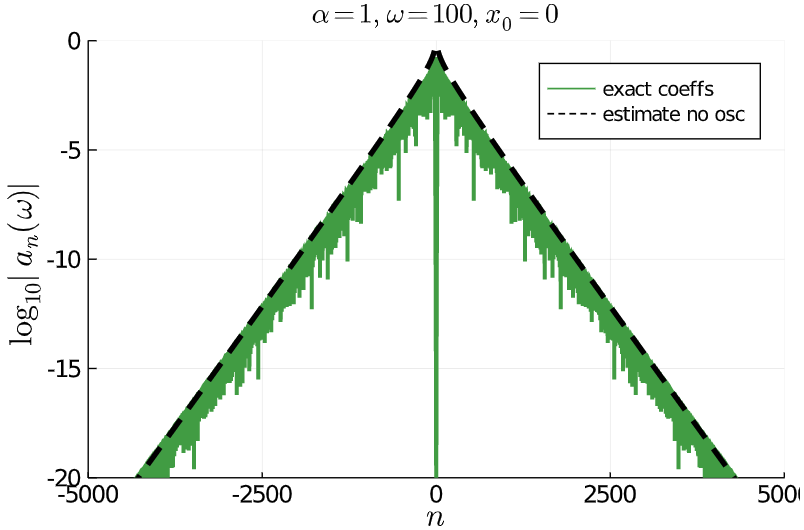}	
  	\includegraphics[width=.45\textwidth]{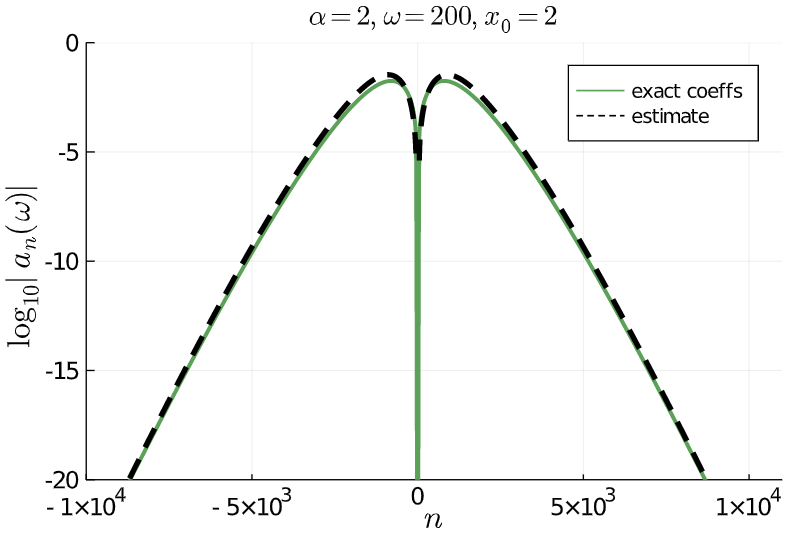}
  \end{center}
  \caption{\textit{MT functions:} Equation \eqref{eq:slimlinemt} suggests that for $x_{0}=0$ (left plot) we expect the rate of decay of coefficients to resemble a straight line and this agrees with the exact coefficients. Note that to aid the reader, we have removed the oscillations in our estimate in the left plot. As we increase $x_{0}$ (right plot), we observe that the coefficients form a hump. Our estimate is represented as a dark dashed line and the exact coefficients are denoted by a light solid line. \label{fig:mtcompare1} }
 \end{figure}

Out of our three bases, Malmquist--Takenaka functions are the strongest contender when approximating wave packets on the real line for two reasons: firstly, the asymptotic estimate in Theorem~\ref{thm:mtfunc} demonstrates spectral decay and secondly, the coefficients can be approximated using FFT. To give  a brief insight into our findings, the asymptotic estimate of coefficients is
\begin{equation}
	\left|a_{n}(\omega) \right| \leq \frac{\left[ 1 + \mathcal{O} \left(\omega^{-1}\right)\right]}{\sqrt{2 \omega} \left( \left| n / \omega \right| -1/4\right)^{1/4}}\exp\!\left(-\alpha \left(\sqrt{\left|\frac{n}{\omega} \right| -\frac14} - x_{0}\right)^{2} \right)\!, \label{eq:slimlinemt}
\end{equation}
for $\omega \geq 0 $. This holds for large $|\omega|$ regime and uniformly for $n$ satisfying
$1/4 \leq \left| n/\omega	\right| \leq c$ for any given constant $c$ (this is made more precise in Section 4). We can immediately deduce exponential decay with respect to $n$ within this regime relative to $\omega$. This is be confirmed experimentally in Fig.~\ref{fig:mtcompare1}. 

This bound highlights how the presence of high oscillations can be an advantage.  In the case of no oscillations, $\omega =0$, Weideman found in \cite{weideman1994computation} that the coefficients decay sub-exponentially. However, counterintuitively, as  $\omega >0$ grows, the estimate in \eqref{eq:slimlinemt} exhibits exponential convergence for $n \leq c \cdot \omega$ for some $c$. In other words,  the presence of oscillations means that we see exponential convergence within a domain and the range of the domain increases for larger values of $\omega$.  In the context of practical computation, for very highly oscillatory wave packets the coefficients will essentially exhibit exponential decay within the practical range of accuracy.

 Our asymptotic estimate shows that we should expect our coefficients to decay spectrally. From \eqref{eq:slimlinemt}, we can readily deduce the location of the peak which is given by 
 \begin{displaymath}
  n = \omega\! \left(x_{0}^{2} + \frac{1}{4}\right),  \quad \textnormal{for $\omega \geq 0$}.
\end{displaymath}

Interestingly, the location of the peak does not depend on how much we stretch the wave packet. The number of coefficients $n>0$ that are greater than $\varepsilon$ in absolute value is asymptotically equal to
\begin{displaymath}
 \omega \!\left[ \frac14 + \left( x_{0} + \sqrt{- \frac{1}{\alpha} \log \frac{\varepsilon}{c} }\right)^{2} \right].
\end{displaymath}

That is the number of coefficients required to resolve a wave packet in MT basis grows linearly with $\omega$, which beats Hermite functions where the number of coefficients grow quadratically in $\omega$!

Overall, MT functions seem to provide the most practical approach in approximating wave packets. The coefficients can be approximated rapidly using the FFT and we have proved that they exhibit asymptotic  exponential decay.  A more in-depth discussion about these bases can be found in section \ref{sec:mt}. 

Detailed derivation of the above asymptotic estimates is discussed in the sequel.

\setcounter{equation}{0}
\section{Stretched Fourier basis}\label{sec:sf}

Stretched Fourier functions are the familiar Fourier functions scaled to the interval $[-\lambda, \lambda]$, where $\lambda>0$ can be chosen based on the function $f$ to be approximated and a desired error tolerance $\varepsilon$. Fourier functions are the natural basis in the presence of periodic boundary conditions in a compact interval, but in our setting we are truncating a function which is defined on the whole real line to the interval $[-\lambda,\lambda]$, and must contend with an error introduced by the Gibbs effect at the endpoints. The basis functions, for a subspace of $\mathrm{L}_2(\mathbb{R})$, have the form 
\begin{equation}\label{eq:SFfunc}
\varphi_{n} (x) = \begin{cases}  (2\lambda)^{-1/2}\ee^{\ii \pi n  x/\lambda}&\mbox{if } x \in [-\lambda,\lambda], \\
0 & \mbox{for }x \in \mathbb{R}\setminus[-\lambda,\lambda]. \end{cases}
\end{equation}

\setcounter{section}{1}
 \setcounter{section}{2}
 \setcounter{figure}{0}
		
For a wave packet, $f(x) = \ee^{-\alpha(x-x_{0})^{2}}\cos(\omega x)$, where $\alpha > 0$, $x_{0},\omega \in \mathbb{R}$, an integral expression for the expansion coefficient is found by substituting \eqref{eq:SFfunc} into \eqref{eq:generala} to get
\begin{equation}
a_{n,\lambda}(\omega)=\frac{1}{(2\lambda)^{1/2}} \int_{-\lambda}^\lambda \ee^{-\alpha(x-x_0)^2}\cos(\omega x) \exp\!\left(-\frac{\pi\ii nx}{\lambda}\right)\!\dd x.  \label{eq:aSF}
\end{equation}
It is clear that a set of $N$ of these coefficients can be approximated by a single FFT in $\mathcal{O}(N\log_2 N)$ operations.

\begin{figure}[tbh]
  \begin{center}
    \includegraphics[width=.60\textwidth]{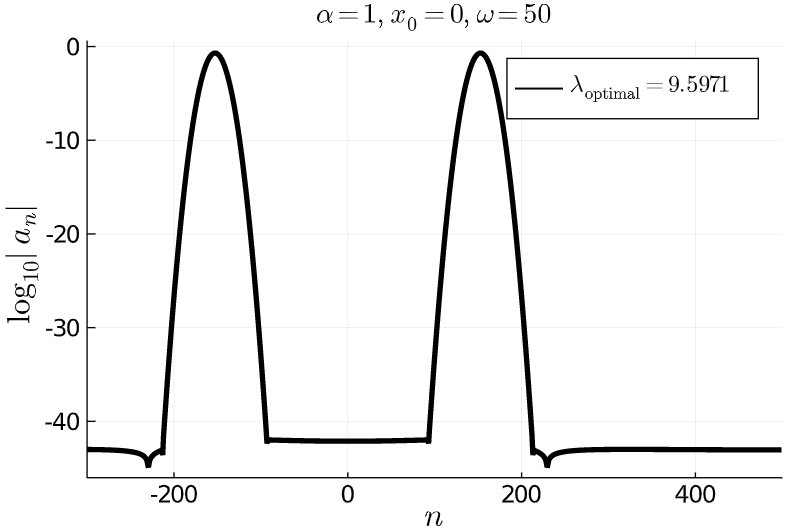}
  \end{center}
  \caption{Plot of  stretched Fourier coefficients $a_{n,\lambda}(\omega)$ in \eqref{eq:aSF} of a wave packet. Here $\lambda= 9.5971$ originates in equation \eqref{eq:sflambda} with $\varepsilon = 10^{-40}$.  \label{fig:sombrero}}
\end{figure}

Dietert and Iserles in their unpublished report \cite{dietert17far} showed that stretched Fourier expansions for functions of the form $f(x) = g(x) \exp(-\alpha x^2)$ (where $g(x)$ is an entire function which grows at most exponentially) converge at spectral speed down to a specified error tolerance, provided that care is taken when choosing $\lambda$. An appealing result shown therein is the `sombrero phenomenon', whereby the logarithm of the absolute value of the expansion coefficients is shaped like a sombrero. For the special case of the wave packets considered in this paper, we apply their method to arrive at the result in Theorem \ref{thm:sf}, and give an example to demonstrate the sombrero effect in Figure \ref{fig:sombrero}.

\begin{theorem}[Coefficients for Stretched Fourier] \label{thm:sf}
  Suppose that $f(x) = \ee^{-\alpha( x-x_{0})^{2}}\cos(\omega x)$, where $\alpha >0$, $x_{0}, \omega \in \mathbb{R}$. The coefficients of $f$ in the stretched Fourier basis (as in equation \eqref{eq:SFfunc}) with stretching parameter $\lambda \geq |x_0|$ satisfy
  \begin{Eqnarray}
	|a_{n,\lambda}(\omega)|&\leq& \left(\frac{\pi}{2\alpha\lambda}\right)^{1/2} \left[ \exp\!\left(\!-\frac{(\pi |n|-\lambda\omega)^2}{4\alpha\lambda^2}\right) + \exp\!\left(\!-\alpha(\lambda - |x_0|)^2\right) \right]\!. \label{eq:sfbound}
  \end{Eqnarray}
for all $n \in \mathbb{Z}$.
\end{theorem}

Before proving the theorem, let us discuss this sombrero phenomenon in more detail. We do not expect the reader to be familiar with millinery terminology, so allow us to explain that the main part of a hat in which the head sits is called the \emph{crown\/} and the projecting edge is called the \emph{brim\/}. Drawing the analogy with a sombrero seen in Figure 2.1, the crown of the sombrero exhibits super-exponential decay from the central peak and the brim of the sombrero is where the decay rate is very slow -- specifically, $\mathcal{O}(n^{-1})$, as predicted by standard Fourier theory for a function of bounded variation. Note however, that if $\lambda$ is sufficiently large, as it turns out to be the case in Fig.~\ref{fig:sombrero}, it is also a region where the magnitude of the coefficients is smaller than the accuracy we need for practical computation. Thus we can truncate the series for an approximation scheme and attain effective spectral convergence down to our desired accuracy. 

\begin{proof}[Proof of Theorem \ref{thm:sf}]
  
  Essentially, we evaluate the integrals in Lemma~3 and~4 from \cite{dietert17far} exactly.  Consider the modified coefficient given by an integral over the whole real line, which is much more amenable to direct calculation:
  \begin{Eqnarray*}
    \check{a}_{n,\lambda}(\omega)&=&\frac{1}{(2\lambda)^{1/2}} \int_{-\infty}^\infty \ee^{-\alpha(x-x_0)^2} \cos(\omega x) \exp\!\left(-\frac{\pi\ii nx}{\lambda}\right)\! \dd x\\
    &=& \left(\frac{\pi}{2\alpha\lambda}\right)^{1/2}\frac12\left[\ee^{\ii x_0(\pi n-\lambda\omega)/\lambda} \exp\!\left(\!-\frac{(\pi n-\lambda\omega)^2}{4\alpha\lambda^2}\right) +\ee^{\ii x_0(\pi n+\lambda\omega)/\lambda} \exp\!\left(\!-\frac{(\pi n+\lambda\omega)^2}{4\alpha\lambda^2}\right)\!\right]\!.
  \end{Eqnarray*}
  From this we have the bound,
  \begin{Eqnarray*}
  | \check{a}_{n,\lambda}(\omega)| &\leq&\frac12 \left(\frac{\pi}{2\alpha\lambda}\right)^{1/2} \left[ \exp\!\left(\!-\frac{(\pi n-\lambda\omega)^2}{4\alpha\lambda^2}\right) + \exp\!\left(\!-\frac{(\pi n+\lambda\omega)^2}{4\alpha\lambda^2}\right)\! \right] \\
   &\leq& \frac12 \left(\frac{\pi}{2\alpha\lambda}\right)^{1/2} \exp\!\left(\!-\frac{(\pi |n|-\lambda\omega)^2}{4\alpha\lambda^2}\right) \left[1  + \exp\!\left(\!-\frac{\pi |n|\omega}{\alpha\lambda}\right)\! \right] \\
   &\leq& \left(\frac{\pi}{2\alpha\lambda}\right)^{1/2} \exp\!\left(\!-\frac{(\pi |n|-\lambda\omega)^2}{4\alpha\lambda^2}\right).
  \end{Eqnarray*}
  A bound for the original coefficient can then be found as follows.
  \begin{Eqnarray*}
    |a_{n,\lambda}(\omega)|&\leq& | \check{a}_{n,\lambda}(\omega)|+\frac{1}{(2\lambda)^{1/2}} \left[\int_\lambda^\infty \ee^{-\alpha(x-x_0)^2}\dd z+\int_{-\infty}^{-\lambda} \ee^{-\alpha(x-x_0)^2}\dd x\right]\\
    &=&|\check{a}_{n,\lambda}(\omega)|+\frac{1}{(2\alpha\lambda)^{1/2}} \left[\int_{\alpha^{1/2}(\lambda-x_0)}^\infty \ee^{-x^2}\dd x+\int_{\alpha^{1/2}(\lambda+x_0)}^\infty \ee^{-x^2}\dd x\right]\! \\
    &\leq& |\check{a}_{n,\lambda}(\omega)|+\left(\frac{2}{\alpha\lambda}\right)^{1/2} \int_{\alpha^{1/2}(\lambda - |x_0|)}^\infty \ee^{-x^2} \dd x \\
    &\leq& |\check{a}_{n,\lambda}(\omega)|+\left(\frac{\pi}{2\alpha\lambda}\right)^{1/2} \ee^{-\alpha(\lambda - |x_0|)^2}.
  \end{Eqnarray*}
  In the last step we have used the known inequality,
  \begin{displaymath}
  \int_t^\infty \ee^{-x^2}\dd x\leq \frac{\sqrt{\pi}}{2}\ee^{-t^2},
  \end{displaymath}
  which holds for $t \geq 0$. Substituting the bound for $|\check{a}_{n,\lambda}(\omega)|$ obtained above yields the result of the theorem.
  \end{proof}

\subsection{The choice of stretching parameter $\lambda$}

\begin{figure}[tbh]
  \begin{center}
     \includegraphics[width=.60\textwidth]{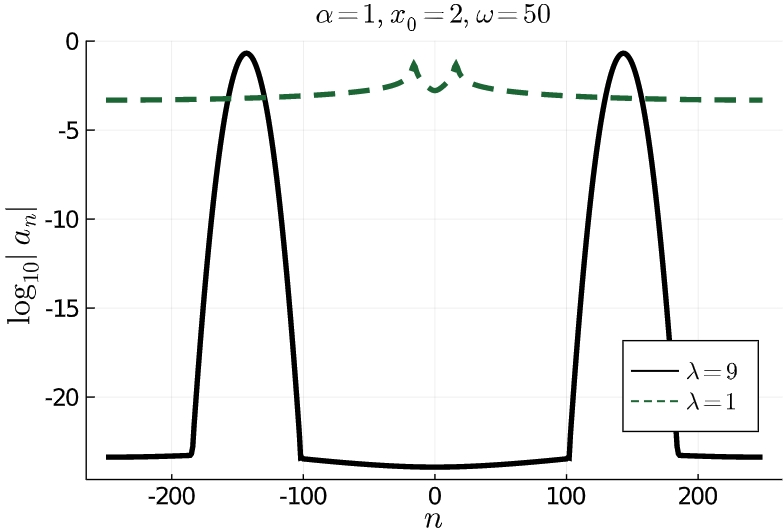}
  \end{center}
  \caption{\textit{Stretched Fourier functions:} To understand the characteristics of our bound in \eqref{eq:sfbound}, which requires us to choose $\lambda \geq |x_{0}|$, this plot shows the effect of choosing $\lambda$ less than and greater than $|x_{0}|$. 
%The left figure shows the plot of a wave packet for $\alpha =1, x_{0}=2$ and $\omega = 50$. The green dashed line represents the interval from -1 to 1, and the black solid line shows the interval from -9 to 9. The green dashed interval only covers the left tail of the wave packet, whereas the interval from -9 to 9, the black solid line, contains most of the mass of the highly oscillatory wave packet. %
The figure shows the plot of the logarithm of the absolute value of coefficients. For $\lambda = 1$, green dashed line, we see slow rate of decay as the majority of the `mass' of the wave packet is outside the interval [-1,1]. For $\lambda = 9$, we contain the majority of the wave packet in the interval [-9,9] and we see that the coefficients decay spectrally down beyond $10^{-20}$.  \label{fig:sflambdacompare}}
\end{figure}

It was found in \cite{dietert17far} that the position of the sombrero brim depends on our choice of $\lambda$. In order to choose an appropriate $\lambda$, we need to understand the characteristics of our bound in \eqref{eq:sfbound}. It requires us to choose $\lambda \geq |x_{0}|$, if we instead choose $\lambda$ less than $|x_0|$ then we obtain a poor approximation to the original wave packet. This can be seen in Fig.~\ref{fig:sflambdacompare}, where the wave packet has much of its `mass' outside the interval of approximation. The bound also suggests that we are limited by the magnitude of $\exp\left(-\alpha (\lambda -|x_0|)^2\right)$. If this quantity is not smaller than some desired threshold $\varepsilon$ then the coefficients are not guaranteed to decay spectrally below $\mathcal{O}(\varepsilon)$, but instead we expect them to decay like $\mathcal{O}\left(n^{-1}\right)$, until they reach such a threshold. This can be seen in Figure \ref{fig:sflambdacompare}, for $\lambda = 1$, the $\exp\left(-\alpha (\lambda -|x_0|)^2\right)$ term for $\alpha = 1$ and $x_{0}=2$ becomes $\exp(-1) \approx 0.368$, resulting in slow decay. Whereas, for $\lambda=9$, $\exp(-49) \approx 5.24\times10^{-22}$ so we see spectral decay past $10^{-20}$.

There are many effective ways to define $\lambda$; for a given $\epsilon > 0$, a simple choice is,
\begin{equation}\label{eqn:lambdachoice}
\lambda_{\textnormal{optimal} }= |x_0| + \sqrt{\frac{\log(\varepsilon^{-1})}{\alpha}},
\end{equation}
in order to guarantee
\begin{equation}\label{eqn:sfcoefflambdachoice}
|a_{n,\lambda}(\omega)| \leq \left(\frac{\pi}{2\alpha\lambda}\right)^{1/2} \left[ \exp\!\left(\!-\frac{(\pi |n|-\lambda\omega)^2}{4\alpha\lambda^2}\right) + \varepsilon \right]\!. 
\end{equation}

\begin{figure}[th]
  \begin{center}
    \includegraphics[width=.60\textwidth]{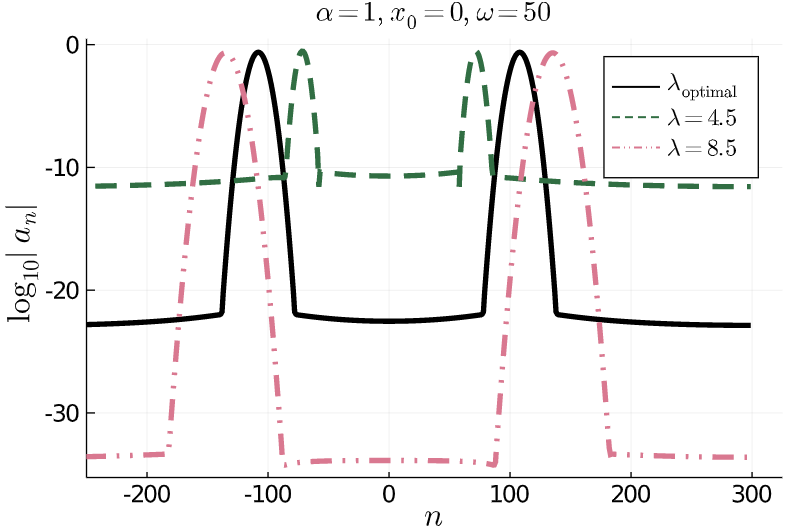}
  \end{center}
  \caption{\textit{Stretched Fourier functions:} For $|a_{n}(\omega)| \geq 10^{-20}$, $\lambda_{\mathrm{optimal}}$ is roughly 6.7861. The plot shows the effect of under- and overshooting. If we use a smaller $\lambda$, no matter how many coefficients we use, we would not arrive at our desired accuracy in reasonable time. Once  $\lambda>\lambda_{\mathrm{optimal}}$, though, we reach ultimately  excessive accuracy -- but require more coefficients to reach the $10^{-20}$ tolerance. \label{fig:sombrero_compare}}
\end{figure}

Fig.~\ref{fig:sombrero_compare} demonstrates the dependence of the wave packet coefficients on the choice of $\lambda$. If we choose $\lambda$ to be smaller than $\lambda_{\mathrm{optimal}}$ we would need to choose $n={\mathcal O}(\varepsilon^{-1})=\mathcal{O}(10^{20})$ in order to reach our desired accuracy, which is unrealistically large. For $\lambda$ larger than our optimal value the brim of the sombrero lies lower but the thickness of the crown is greater and we need to choose large (at least, only linearly) $n$ to reach $10^{-20}$ accuracy. Thus even though we get spectral decay with stretched Fourier bases, this approach is extremely sensitive to how well we can define the optimal $\lambda$. This effect can also be seen theoretically in Theorem \ref{thm:sf}.

In a more practical setting, a linear combination of wave packets is a more realistic approximate solution to the semi-classical TDSE and due to our current method of determining $\lambda_{\mathrm{optimal}}$ we need more knowledge about the exact solution before we can approximate it. Our analysis indicates that it is safer to overestimate $\lambda$ than to underestimate it, providing a clue to a good strategy for a linear combination of wave packets. The current practical limitations of stretched Fourier basis notwithstanding, the functions still hold significant potential due to their simplicity, fast approximation of coefficients using FFT and spectral rate of decay.

\subsection{Number of coefficients required for a given tolerance}

 Given a tolerance $\varepsilon > 0$, once we choose $\lambda$ as in equation \eqref{eqn:lambdachoice} we can guarantee that the coefficients satisfy \eqref{eqn:sfcoefflambdachoice}.  To gain a rough idea on the number of coefficients required for a given $\varepsilon$, we do the following overestimation by considering $|a_{n,\lambda}(\omega)| \geq 2\varepsilon$ and finding all the $n$'s that satisfy this bound.  That is all $n$ satisfying
\begin{Eqnarray*}
	(\pi |n| - \lambda \omega)^2 &\leq& 4\alpha \lambda^{2} \left[  \log\left( \varepsilon^{-1} \right) - \log \left( (2(2\alpha\lambda/\pi)^{1/2} - 1) \right)\right]\!.
\end{Eqnarray*}
 
Neglecting the complicated term inside the logarithm (i.e. bounding this inequality from above), we have 
\begin{displaymath}
  |n|\leq \frac{\omega+2\sqrt{\alpha}\sqrt{\log\varepsilon^{-1}}}{\pi}\lambda
\end{displaymath}
and, substituting $\lambda=\lambda_{\mathrm{optimal}}$ from \eqref{eqn:lambdachoice}, we conclude with 
\begin{displaymath}
  |n|\leq \frac{\omega+2\sqrt{\alpha}\sqrt{\log\varepsilon^{-1}}}{\pi} \left(|x_0|+\sqrt{\frac{\log\varepsilon^{-1}}{\alpha}}\right)\!,
\end{displaymath}
which emphasises the role different parameters play in the bound. 

The important factor to note here is the linear dependence on both $\omega$ and $|x_0|$, which is extremely favourable compared to the other two bases considered in this paper. The downside, though, is that we require prior knowledge of $x_0$, $\alpha$ and $\omega$ to choose $\lambda$ and that the approximation of a linear combination of wave packets is likely to result in less favourable behaviour.

\setcounter{equation}{0}
\section{Hermite functions}\label{sec:hermite}

It is natural to think of Hermite functions in the context of quantum mechanics since they are  solutions of the quantum harmonic oscillator. Hermite functions have a skew-symmetric and tridiagonal differentiation matrix, which implies that their associated spectral methods are stable and preserve unitarity \cite{hairer2016numerical}. However, computation of coefficients in an expansion in Hermite functions is associated with a number of practical challenges. An obvious method of computing expansion coefficients in Hermite functions with quadrature is by explicit integration of each coefficient at the overall cost of $\mathcal{O}(N^{2})$ operations. Nonetheless they are still a contender when approximating wave packets and, less intuitively, it is possible to compute the expansion coefficients in $\mathcal{O}(N(\log N)^2)$ operations using the fast multipole algorithm \cite{boyd01cfs}.

\setcounter{section}{2}
\begin{figure}[tb]
  \begin{center}
  \includegraphics[width=.45\textwidth]{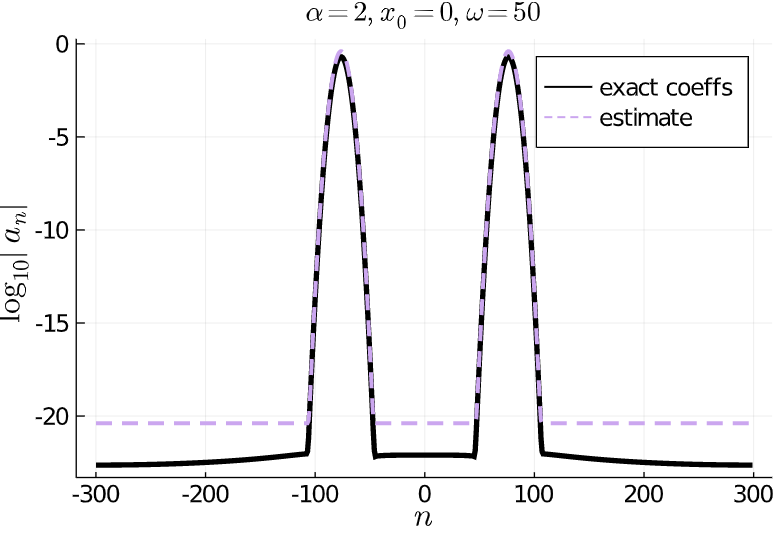}\hspace*{20pt}
  \includegraphics[width=.45\textwidth]{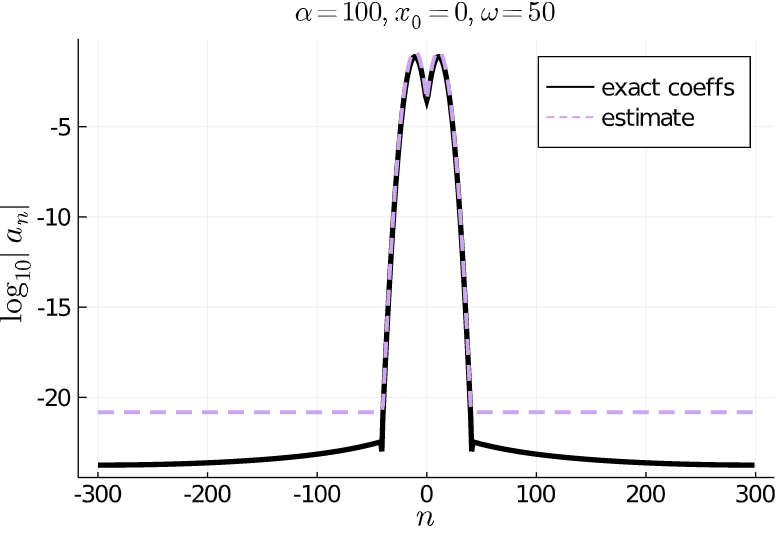}
  \includegraphics[width=.45\textwidth]{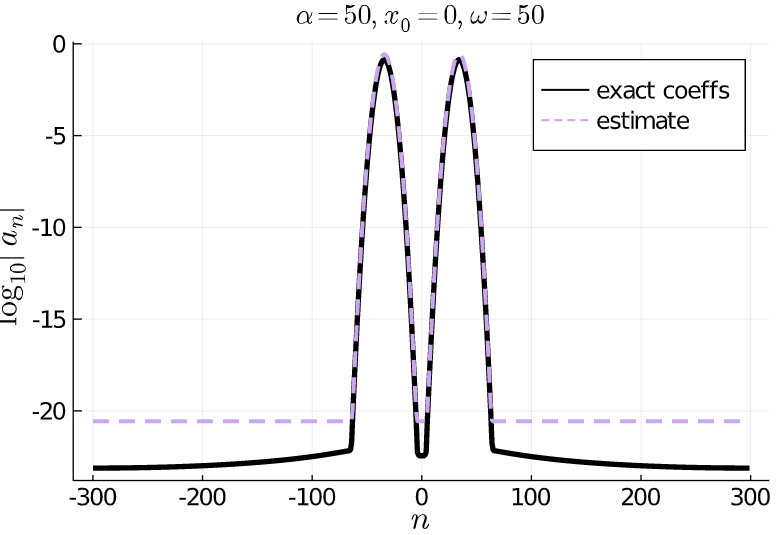}\hspace*{20pt}
  \includegraphics[width=.45\textwidth]{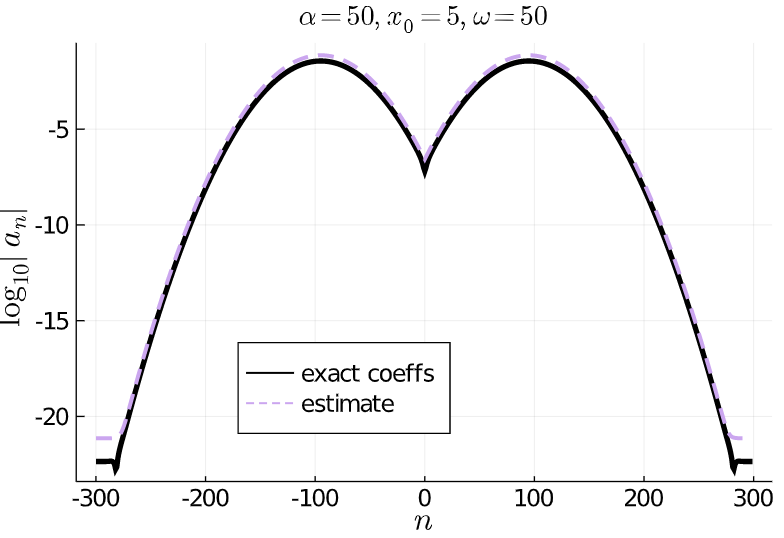}
    \end{center}
  \caption{\textit{Stretched Fourier functions:\/} Plot of \eqref{eqn:sfcoefflambdachoice} where the optimal $\lambda$ is used, for fixed accuracy of $\varepsilon=10^{-20}$, and varying $\alpha$ and and $x_{0}$. This demonstrates experimentally that the bound in Theorem \ref{thm:sf} is representative of the behaviour of the coefficients in the crown of the sombrero, although the estimate of the beginning of the brim is universally pessimistic. Changes in the $\alpha$ parameter have minimal effect on the coefficients, whereas changes in $x_0$ have a significant effect.}
\end{figure}

\setcounter{section}{3}
\setcounter{figure}{0}
We consider normalised Hermite functions
\begin{Eqnarray}\label{eqn:Hermitebasis}
	\varphi_{n}(x) := \frac{1}{(2^{n}n!\pi^{1/2})^{1/2}} \HH_n(x)\ee^{-x^2/2}, \qquad n \in \mathbb{Z}_+.
\end{Eqnarray}
The corresponding expansion integral is
\begin{Eqnarray}
  \nonumber
  a_{n}(\omega)&=&\frac{1}{(2^{n}n!\pi^{1/2})^{1/2}}\int_{-\infty}^\infty \ee^{-\alpha(x-x_0)^2} \cos(\omega x) \HH_n(x)\ee^{-x^2/2}\dd x,\qquad n\in\Z_+,  \\
  &=& \frac{1}{(2^{n}n!\pi^{1/2})^{1/2}} \textnormal{Re}\left[ \int_{-\infty}^\infty \ee^{-\alpha(x-x_0)^2-x^2/2 + \ii \omega x}  \HH_n(x)\dd x\right]. \label{eq:hermitean}
\end{Eqnarray}
where the $\HH_n$s are Hermite polynomials,
\begin{displaymath}
	\int_{-\infty}^\infty \HH_m(x)\HH_n(x)\ee^{-x^2}\dd x = 0 \quad \textnormal{ for } m \neq n, \quad \textnormal{ and } \quad \int_{-\infty}^\infty \HH_n^2(x)\ee^{-x^2}\dd x = \pi^{1/2}2^{n}n!.
\end{displaymath}

\begin{theorem}[Coefficients for Hermite functions]\label{thm:hermitefunc}
	Suppose that $f(z) = \ee^{-\alpha( z-x_{0})^{2}}\cos(\omega z)$, where $\alpha > 0$, $x_{0}$ and $\omega$ are real. The asymptotic behaviour of the coefficients, $a_{n}$, with Hermite basis functions in equation \eqref{eqn:Hermitebasis} is given by 

\begin{Eqnarray*}
  \textnormal{For } \alpha >\frac12:\quad a_{n} &\sim &\frac{1}{2 \sqrt{\omega(2\alpha +1)}}\! \left[ \frac{4\alpha^{2} -1}{1 + 2\frac{n}{\omega^{2}}(4\alpha^{2}-1)} \right]^{\!1/4}\!\left| \frac{1 +\sqrt{1+2\frac{n}{\omega^{2}}\left(4\alpha^{2}-1\right)}}{\sqrt{2\frac{n}{\omega^{2}} \left(4\alpha^{2}-1 \right)}} \right|^{1/2} \\
  & & \mbox{}\times \left| \frac{1 + \sqrt{1 +2\frac{n}{\omega^{2}}\left(4\alpha^{2}-1\right)}}{\sqrt{2\frac{n}{\omega^{2}}}(2\alpha +1)}\right|^{n} \exp\! \left( - \frac{\omega^{2}}{2(2\alpha + 1)} +\frac{n}{1+\sqrt{1 +2\frac{n}{\omega^{2}}\left(4\alpha^{2}-1\right)}} \right.\\
  &&\hspace*{25pt}\left.\mbox{}+ \frac{\alpha x_{0}^{2}}{4\alpha^{2}-1}\left[1 -  \frac{2 \alpha}{\sqrt{1+2\frac{n}{\omega^{2}}\left(4\alpha^{2}-1\right)}}\right]+ \mathcal{O}\left(\omega^{-2}\right)\right) \\
  & &\mbox{} \times\textnormal{Re}\!\left\{  \left[ 1 + \frac{\ii \alpha	x_{0}}{\omega \left[ 1 + 2\frac{n}{\omega^{2}}\left(4\alpha	^{2} -1\right) \right]} + \mathcal{O}(\omega^{-2}) \right] \right. \\
  & &\hspace*{12pt}\mbox{}  \times \exp\! \left( \frac{\ii n \pi}{2} - \frac{2\ii \alpha x_{0} \omega \left(-2\alpha + \sqrt{1 + 2\frac{n}{\omega^{2}}\left(4\alpha^{2}-1\right)}\right)}{4\alpha^{2}-1} \right.\\
  &&\hspace*{25pt}\left.\left. \mbox{}+ \frac{\ii x_{0} \alpha \left[-3 + 6\frac{n}{\omega^{2}} + 8 \frac{n}{\omega^{2}}\alpha^{2}\left(x_{0}^{2} - 3\right)\right]}{3\omega\left[1 + 2\frac{n}{\omega^{2}}\left(4\alpha^{2}-1\right)\right]^{3/2}}\right) \right\}\!\left[1 + \mathcal{O}\left(n^{-1}\right) \right], \\[4pt]
  \textnormal{For }  \alpha=\frac12: \quad \quad \ & & \\
   x_{0} \neq 0: \quad a_{n}&=& \sqrt{\frac{\pi^{1/2}}{2^{n}n!}}(x_0^2 + \omega^2)^{n/2} \ee^{-(x_0^2+\omega^2)/4} \cos\left(\frac{\omega x_{0}}{2} + n \tan^{-1}\left(\frac{\omega}{x_0}\right) \right),\\
  x_{0} = 0: \quad a_{n}  &=& \sqrt{\frac{\pi^{1/2}}{2^{n}n!}}\omega^{n} \ee^{-\omega^2/4} \cos\left(\frac{n\pi}{2}\right),\\
   \nonumber
    \textnormal{For } \alpha <\frac12: \quad \quad \ & & \\
  n \neq \frac{\omega^{2}}{2 \left(1-4\alpha^{2}\right)}: \quad a_{n} &\sim & \textnormal{Re} \left\{-\frac{1}{\sqrt{\omega(1+2\alpha)}}\left| \frac{2\alpha -1}{2 \alpha + 1}\right|^{n/2} \left(\frac{|1 - 4\alpha^{2}|}{2 \frac{n}{\omega^{2}} |1-4\alpha^{2}| -1}\right)^{1/4}  \right. \\
	& &\times  \left(1 - \frac{\ii \alpha x_{0}}{\omega \left(2\frac{n}{\omega^{2}} |1 - 4\alpha^{2}| -1\right)} + \mathcal{O} \left(\omega^{-2}\right)\right)  \left( \frac{1+ \sqrt{1-2\frac{n}{\omega^{2}} |1-4\alpha^{2} |}}{\sqrt{2 \frac{n}{\omega^{2}} |1-4\alpha^{2}|}}\right)^{n+1/2}\\
		& & \times \exp\left\{\frac{\alpha}{|1-4\alpha^{2}|}\left( \omega^{2}- x_{0}^{2} +2 x_{0} \omega \sqrt{2 \frac{n}{\omega^{2}}|1-4\alpha^{2}| -1} \right) \right. \\
		& & + \frac{\alpha x_{0} \left(8 \alpha^{2} x_{0}^{2}\frac{n}{\omega^{2}} + 6 \frac{n}{\omega^{2}} |1-4\alpha^{2}| - 3\right)}{3\omega \left(2\frac{n}{\omega^{2}}|1-4\alpha^{2}|-1\right)^{3/2}} \\
		& & +  \ii \left[\left(n + \frac{1}{2} \right)\frac{\pi}{2} - \frac{ 4 \alpha^{2} x_{0} \omega}{|1-4\alpha^{2}|} + \frac{\omega^{2}}{2 |1-4\alpha^{2}|}\sqrt{2 \frac{n}{\omega^{2}} |1-4\alpha^{2}| -1}  \right. \\
		& & \left.\left. \left. + \frac{2 \alpha^{2} x_{0}^{2}}{|1-4\alpha^{2}|\sqrt{2\frac{n}{\omega^{2}} |1-4\alpha^{2}|-1}}\right]\right\} \right\}\!\left[1 + \mathcal{O}\left(n^{-1}\right) \right].%,\\
%n =  \frac{\omega^{2}}{2 \left(1-4\alpha^{2}\right)}: \quad a_{n} &\sim &  \textnormal{Re} \left\{ - \frac{(2n)^{1/4}}{\sqrt{1+2\alpha}} \left| \frac{1 - 2 \alpha}{1+ 2 \alpha} \right|^{n/2} \left(\frac{\ii^{7/4} |1-4\alpha^{2}|}{\sqrt{2} \left(\alpha x_{0}\right)^{1/4} \omega^{3/4}}\right) \left(1 - \frac{4\alpha^{2}x_{0}^{2} + |1-4\alpha^{2}|}{\ii 16 \alpha x_{0} \omega }\right) \right. \\
%	& & \left. \times \exp \left[ \frac{\ii \pi }{4} - \frac{\alpha x_{0}^{2}}{|1-4\alpha^{2}|}+ \frac{(\ii \pi + 4 \alpha)\omega^{2}}{4 |1-4\alpha^{2}|} - \frac{4 \ii \alpha^{2} x_{0} \omega}{|1-4\alpha^{2}|} + \frac{8 \ii^{1/2} \left(\alpha x_{0}\right)^{3/2} \omega^{1/2}}{3 |1-4\alpha^{2}|}\right. \right.\nonumber \\
%	& & \left. \left. - \frac{\ii^{3/2}}{5 \sqrt{\omega}}  \left(5 \sqrt{\alpha x_{0}} + \frac{4\left(\alpha x_{0}\right)^{5/2}}{|1-4\alpha^{2}|} \right) \right]\right\}\!\left[1 + \mathcal{O}\left(\omega^{-2}\right) \right]. \nonumber
\end{Eqnarray*}
The asymptotic estimate for $\alpha >1/2$ holds uniformly for $n > \cw \omega^{2}$ for any given constant $\cw$. For $\alpha <1/2$, this estimate holds uniformly for $n < \cw \omega^{2}$, $\cw$ is some constant, and $n$ outside of any given neighbourhood of $\omega^{2}/(2(1-4\alpha^{2}))$. These bounds hold uniformly as $\omega \rightarrow \infty$.
\end{theorem}

The behaviour of the coefficients is hard to deduce by a direct examination of this expansion. It is helpful to break down the result into three regimes. Note that in a more practical setting we would expect $\alpha>1/2$, so the following analysis considers the asymptotics for this case: 
\begin{enumerate}
\item When $\omega$ is large, 
\begin{Eqnarray*}
  a_{n}(\omega) \sim \frac{\mathcal{C}_{\mathrm{osc}}}{\sqrt{\omega}} \exp \left(-\frac{\omega^{2}}{2(2\alpha + 1)}\right),
\end{Eqnarray*}
where $\mathcal{C}_{\mathrm{osc}}$ is a bounded oscillatory term. The asymptotic expression only depends on $\omega$ and we expect exponential decay.
\item When $n$ is large, 
\begin{Eqnarray*}
  a_{n}(\omega) \sim \mathcal{C}\left| \frac{1 + \sqrt{1 +2\frac{n}{\omega^{2}}\left(4\alpha^{2}-1\right)}}{\sqrt{2\frac{n}{\omega^{2}}}(2\alpha +1)}\right|^{n}\exp\left(\frac{n}{1+\sqrt{1+2\frac{n}{\omega^{2}}\left(4\alpha^{2}-1\right)}} \right).
\end{Eqnarray*}
The expression is still fairly opaque but in Fig.~\ref{fig:hermite2compare} we can see that the envelope of the logarithm of the absolute value of coefficients is approaching a straight line. 

\begin{figure}[tbh]
  \begin{center}
  	\includegraphics[width=.60\textwidth]{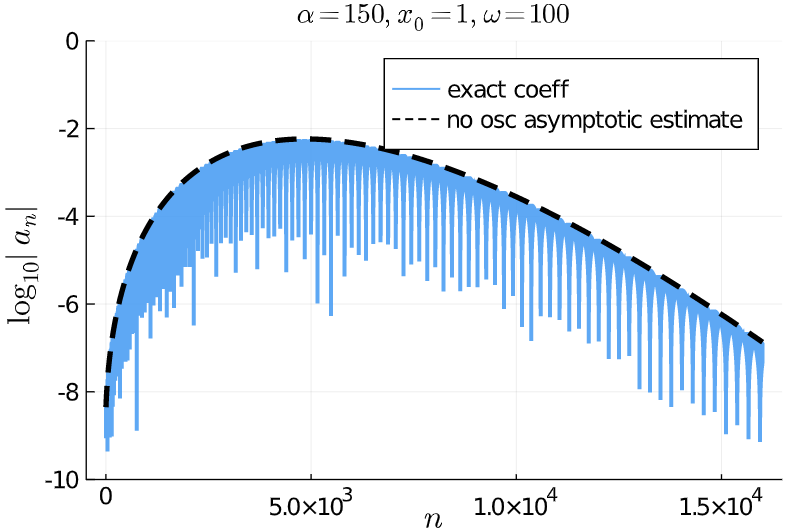}
  \end{center}
  \caption{\textit{Hermite functions:} Increasing $\alpha$ makes our coefficients decay much slower. Just as in Fig.~\ref{fig:hermite1compare} we can see that the our asymptotic estimate approximates the exact coefficients very well. The figure shows how well the asymptotic estimate envelopes the exact coefficients.  \label{fig:hermite2compare}}
\end{figure}
 
\item When both $n$ and $\omega$ are large subject to a power law,  $n=\omega^\lambda$ for some $\lambda>0$, then
\begin{Eqnarray*}
  a_{n}(\omega) &\sim & \frac{\mathcal{C}}{\sqrt{n^{1/\lambda}(2\alpha +1)}} \left| \frac{1 + \sqrt{1 +2n^{1-2/\lambda}\left(4\alpha^{2}-1\right)}}{\sqrt{2n^{1-2/\lambda}}(2\alpha +1)}\right|^{n} \\
  & & \times \exp\left(-\frac{n^{2/\lambda}}{2(2\alpha + 1)}+ \frac{n}{1+\sqrt{1+2n^{1-2/\lambda}\left(4\alpha^{2}-1\right)}}\right),
\end{Eqnarray*}

\begin{figure}[tbh]
  \begin{center}
  	\includegraphics[width=.60\textwidth]{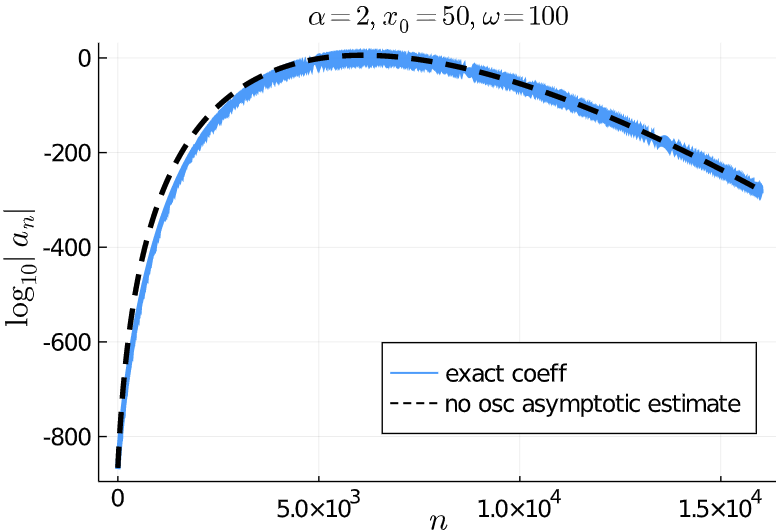}
  \end{center}
  \caption{\textit{Hermite functions:}  For a large shift in $x_{0}$, the estimate does not do too well for small $n$. As $n$ increases the estimate gets better. \label{fig:hermite3compare}}
\end{figure} 

\end{enumerate}

Comparing Fig.~\ref{fig:hermite2compare} to Fig.~\ref{fig:hermite1compare} we observe that increasing $\alpha$ renders the rate of decay of coefficients worse. The same happens when we increase $x_{0}$ in Fig.~\ref{fig:hermite3compare}. 

\subsection{Fourier-type integral of Hermite functions}

The main idea behind the proof is to consider a Fourier type integral of Hermite functions,

\begin{lemma}\label{lemma:fourierhermite}
  It is true that
  \begin{equation}
    \label{eq:8.18}
    \mathcal{H}_{n,c,p}=\frac{1}{\sqrt{2\pi c}}\int_{-\infty}^\infty \ee^{-x^2/(2c)}\HH_n(x)\ee^{\ii px}\dd x  = \left(1-2c\right)^{n/2} \ee^{-cp^2/2} \HH_n\!\left(\frac{cp}{\sqrt{2c-1}}\right) 
  \end{equation}
  for $\textnormal{Re}\, c>0$, $c\neq\frac12$ and
  \begin{equation}
    \label{eq:8.20}
      \mathcal{H}_{n,\frac12,p}=\ee^{-p^2/4} (\ii p)^n.
   \end{equation}
   The branch choice in eqn \eqref{eq:8.18} for $\sqrt{2 c - 1}$,  $0 < c < \frac{1}{2}$, is unimportant because only the even powers of $\sqrt{2 c - 1 }$ occur on the right-hand side due to symmetry of the Hermite polynomials.
\end{lemma}

\begin{proof}

Consider  Hermite functions in a Cauchy integral form,
\begin{displaymath}
	\HH_{n}(z) = \frac{n!}{2 \pi \ii} \int_{\gamma} \ee^{2 w z - w^{2}} \frac{\dd w}{w^{n + 1}},
\end{displaymath}
where the contour $\gamma$ encircles the origin. Substituting this into  \eqref{eq:8.18} and changing the order of integration,
\begin{Eqnarray*}
  \mathcal{H}_{n,c,p}&=& \frac{1}{\sqrt{2 \pi c}}\frac{n!}{2 \pi \ii} \int_{\gamma} \int_{-\infty}^{\infty} \ee^{-x^{2}/(2c) + \ii p x + 2 w x} \dd x \ee^{-w^{2}} \frac{\dd w}{w^{n + 1}} \\
  &=& \frac{n!}{2 \pi \ii} \int_{\gamma} \ee^{-\left(cp^{2}\right)/2 + 2 \ii c p w - (1 - 2 c)w^{2}} \frac{\dd w}{w^{n + 1}} \\
  &=& \ee^{-\left(cp^{2}\right)/2} \frac{n!}{2 \pi \ii }\int_{\gamma} \ee^{2\left(cp/\sqrt{2 c - 1} \right)\xi - \xi^{2}} \frac{\dd \xi}{\xi^{n + 1}}\left(1-2c\right)^{n/2} \qquad \left(	\xi = (1-2 c)^{1/2} w\right) \\
  &=&\left(1-2c\right)^{n/2}\ee^{-\left(cp^{2}\right)/2} \HH_{n}\left(\frac{cp}{\sqrt{2c-1}}\right)\!.
\end{Eqnarray*}
For $c = 1/2$, the proof of \eqref{eq:8.20} is identical (and easier!) to the steps above: we deduce that 
\begin{equation*}
      \mathcal{H}_{n,\frac12,p}=\ee^{-p^2/4} (\ii p)^n.
   \end{equation*}
\end{proof} 

\begin{figure}[tbh]
  \begin{center}
  	\includegraphics[width=.45\textwidth]{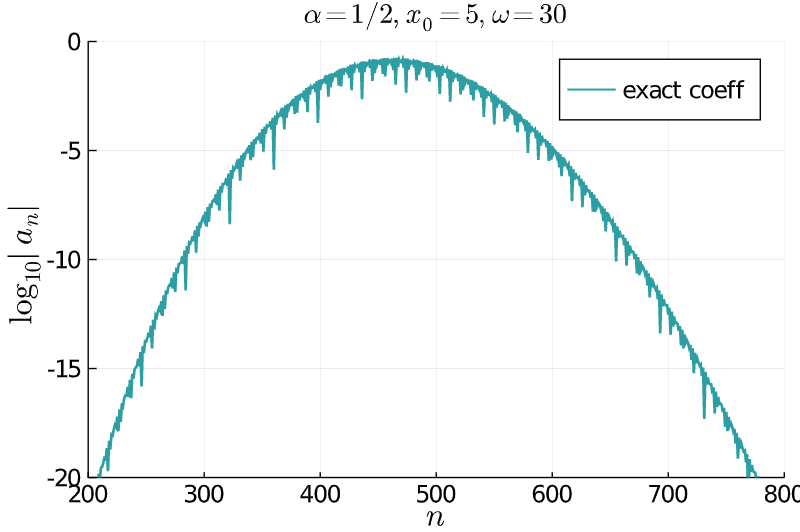}
  	\includegraphics[width=.45\textwidth]{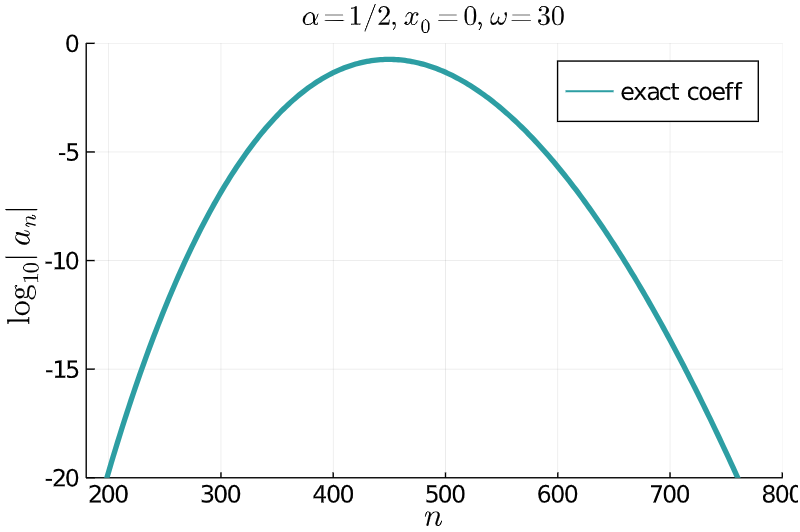}
  \end{center}
  \caption{\textit{Hermite functions:} Plot of \eqref{eq:8.23} for different $x_{0}$. We see that the decay rate of the exact hermite coefficients for $\alpha = 1/2$ resembles a hump. \label{fig:hermitea12plots}}
\end{figure}

We proceed by comparing \eqref{eq:8.18} with \eqref{eq:hermitean} to deduce that
\begin{displaymath}
  c=\frac{1}{1+2\alpha},\qquad  p=\omega-2\ii\alpha x_0,
\end{displaymath}
and
\begin{equation}\label{eqn:Hermiteantilde}
  a_{n}(\omega)=\frac{1}{\sqrt{2^{n}n! \pi^{1/2}}}\ee^{-\alpha x_0^2} \left(\frac{2\pi}{1+2\alpha}\right)^{\!1/2} \textnormal{Re}\,\mathcal{H}_{n,(1+2\alpha)^{-1},\omega-2\ii\alpha x_0}.
  \end{equation}
After some algebraic manipulation (see \ref{sec:alphacases}) and setting 
\begin{Eqnarray*}
  X=\frac{2\alpha x_0}{|4\alpha^2-1|^{1/2}},\qquad Y=\frac{\omega}{|4\alpha^2-1|^{1/2}},
\end{Eqnarray*}
where for simplicity we restrict our analysis to $\alpha,  x_{0}, \omega >0$ so that $X$ and $Y$ are positive real numbers, we note that our expressions are different for different values of $\alpha$,
\begin{Eqnarray}
  \nonumber
  \alpha >\frac12:\quad a_{n}&=&\frac{\pi^{1/4}}{(2^nn!)^{1/2}} \left(\frac{2}{1+2\alpha}\right)^{\!1/2} \left|\frac{1-2\alpha}{1+2\alpha}\right|^{\!n/2}\!\! \exp\!\left(\frac{|1-2 \alpha
 |}{4\alpha}(X^2+2\alpha Y^2) \!\right)\\
  \label{eq:8.22}
  &&\hspace*{30pt}\mbox{}\times \textnormal{Re}\left[\ee^{-\ii|2\alpha-1|XY} \HH_n(X+\ii Y)\right]\!,\\[4pt]
  \label{eq:8.23}
  \alpha=\frac12:\quad a_{n} &=& \frac{\pi^{1/4}}{(2^{n}n!)^{1/2}}(x_0^2 + \omega^2)^{n/2} \ee^{-(x_0^2+\omega^2)/4} \cos\left(\frac{\omega x_{0}}{2} + n \tan^{-1}\left(\frac{\omega}{x_0}\right) \right)\!, \\[4pt]
   \nonumber
  \alpha < \frac12:\quad a_{n}&=& \frac{\pi^{1/4}}{(2^nn!)^{1/2}} \left(\frac{2}{1+2\alpha}\right)^{\!1/2} \left|\frac{1-2\alpha}{1+2\alpha}\right|^{n/2} \exp\!\left(-\frac{|1-2\alpha|}{4\alpha}(X^2+2\alpha Y^2)\right)\\
   \label{eq:8.24}
  &&\hspace*{30pt}\mbox{}\times \textnormal{Re}\left[\ee^{\ii|2\alpha-1|XY + \ii n \pi } \HH_n(Y-\ii X)\right]\!.
\end{Eqnarray}
Thereafter we focus on the asymptotics for $\alpha\neq\frac12$, the remaining case can be also obtained by a limiting process. Note that \eqref{eq:8.23} are the exact coefficients for $\alpha = 1/2$ and the plot of \eqref{eq:8.23} can be seen in Fig.~\ref{fig:hermitea12plots}.

\subsection{An asymptotic estimate}
\subsubsection{The method of steepest descent \label{sec:msdexplain}}
The method of steepest descent (MSD), also known as the saddle point approach, is a key tool used in our Hermite and Malmquist--Takenaka estimates. We give a brief overview here; for a more in-depth review  we refer to \cite{de1981asymptotic, flajolet2009analytic, olver1997asymptotics, wong2001asymptotic}. It is an approximation method that can be applied to integrals of the form 
\begin{Eqnarray*}
	I(t) = \int_{a}^{b} G(z)\ee^{-t F(z)} \dd z,
\end{Eqnarray*}
where $t$ is a large parameter. Let us assume that $F(z)$ is a complex-valued meromorphic function and $G(z)$ is an entire function.

The main idea is to deform the path of integration over a complex contour $\Gamma$ to pass through (or near) certain saddle points of $F(z)$, that is, points $z^{*}$ such that $F'(z^{*})=0$. As well as passing through a subset of saddle points, the path of integration is also designed so that $\textnormal{Im}\,F(z)$ is constant, i.e.~$\textnormal{Im}\, F(z) = \textnormal{Im}\, F(z^{*})$, which ensures that $\exp(-t F(z))$ is non-oscillatory along the contour.
If $F''(z^{*}) > 0$ then $\exp(-t F(z))$ decays rapidly away from $z^{*}$ along this contour -- the path of steepest descent. We assume at  the outset that $G(z)$ can be expanded as convergent power series about the neighbourhood of $z^{*}$.

Skipping many technical details, suppose that there exists an ordered set of contours $\Gamma_1, \Gamma_2,\ldots,\Gamma_s$ such that each contour contains a single saddle point ($z^{*}_{k} \in \Gamma_k$), within each contour $F$ satisfies $\textnormal{Im}\, F(z) = \textnormal{Im}\, F(z^{*}_{k})$, each saddle point in question satisfies $F''(z^{*}_k) > 0$, $G(z)$ can be expanded as a convergent power series about each $z^{*}_{k}$, and such that $\int_a^b \ee^{-t F(z)} \dd z = \sum_{k=1}^s \int_{\Gamma_k} G(z) \ee^{-t F(z)} \dd z$. Then as $t \to \infty$ and $G(z_{k}^{*}) \neq 0$, to first order approximation our integral can be expressed as
\begin{equation}\label{eq:MSD}
  I(t) \sim \sum_{ k =1}^s G(z^{*}_{k})\ee^{-t F(z^{*}_{k})} \sqrt{\frac{2\pi}{ t F''(z^{*}_{k})}} \left[ 1+ \mathcal{O}\left(t^{-1}\right) \right]
  \end{equation}
(\cite[eq.~5.7.2]{de1981asymptotic}). 

\begin{remark}\label{msd_remark}
The power of this approach is that the first approximation depends only on the values of $F$, $F''$ and $G$ at a certain subset of the saddle points, and nowhere else.  Note that the terms hidden in the $\mathcal{O}$ error term in \eqref{eq:MSD} is a polynomial in $G^{'}/G, G^{''}/G, F^{'''}, F^{(IV)}$ and $1/F^{''},$ see \cite[eq.~2.4.16]{olver2010nist},  evaluated at a subset of saddle points.  These terms depend on $z$ and $\kappa$, where $\kappa$ is a carefully chosen parameter with some dependence on $t$ -- more importantly they must remain controlled as $t \rightarrow \infty$ thus imposing constraints on $\kappa$. This ensures that \eqref{eq:MSD} holds uniformly as $t \rightarrow \infty$.
\end{remark}

\subsubsection{Steepest descent: Hermite functions}

To find the asymptotic behaviour of  \eqref{eqn:Hermiteantilde}, we consider $\HH_{n}(z)$ for $z$ given in \eqref{eq:8.22} and \eqref{eq:8.24}, where $n$ is a large parameter and $Y$ is large in the complex argument, $z$. We start by considering Hermite functions in a Cauchy integral form,
\begin{displaymath}
  \HH_n(z)=\frac{n!}{2\pi\ii} \int_\gamma \frac{\ee^{2zw-w^2}}{w^{n+1}}\dd w,
\end{displaymath}
where $\gamma$ is a contour encircling the origin. Let $\zeta=z/\nu$, so that
\begin{displaymath}
  \HH_n(\nu\zeta)=\frac{n!}{2\pi\ii \nu^n} \int_\gamma \ee^{-\nu^2\phi(w)}\frac{\dd w}{w^{1/2}},
\end{displaymath}
where $\phi(w)=w^2-2\zeta w+\frac12\log w$, $\nu=\sqrt{2n+1}$. The saddle points are the zeros of $\phi'$, namely $w_\pm=\frac12(\zeta\pm\ii\sqrt{1-\zeta^2})$. 

Deforming our path of integration to pass through the dominant saddle point, $w_{*}$, then applying \eqref{eq:MSD}, gives us
\begin{equation}
  \label{eq:8.21}
  \HH_n(z)\sim \frac{1}{2\pi i}  \frac{n!}{\nu^{n+1}}\exp \left( -\nu^2\phi(w_*) \right) \!\left[\!\frac{2\pi}{w_*\phi''(w_*)}\!\right]^{\!1/2}\!\left[1+\mathcal{O}\left(n^{-1}\right)\right]
\end{equation}
-- note that $\phi''(w_*)=2-1/(2 w_{*}^{2})$. This asymptotic approximation is valid for large values of $n$ and $\omega$. For small values of $\omega$ a similar contribution from the remaining saddle point should also be included.

By remark \ref{msd_remark}, in order for the estimate in \eqref{eq:8.21} to hold uniformly with respect to $\zeta$ for $\omega \rightarrow \infty$, we require the existence of $c > 0$ such that $c < |w_{*}|<c^{-1}$, which is guaranteed as long as $\zeta$ is bounded. We also require $|\zeta - 1| > \delta$, $\delta$ being a small positive number, to avoid the coalescing of saddle points, in which case we would need to use uniform asymptotic methods, see \cite[\S 23.4]{temme2014asymptotic}.  This leads to the combined conditions  $\delta <|\zeta - 1| < \delta^{-1}$.

We have the key ingredient to find the asymptotic expansion of the Hermite coefficients. For $\alpha > 1/2$,  we want to apply \eqref{eq:8.21} to $\HH_{n}(z)$ in \eqref{eq:8.22}, the dominant root is $w_{-}$ and $w_-\phi''(w_-)=-2\ii(1-\zeta^2)^{1/2}$.  Recall that for this case, $\zeta = (X + \ii Y)/\nu$, where 
\begin{Eqnarray}\label{eq:XandY}
  X=\frac{2\alpha x_0}{|4\alpha^2-1|^{1/2}},\qquad Y=\frac{\omega}{|4\alpha^2-1|^{1/2}},
\end{Eqnarray}
$\omega$ is a large parameter, while $\alpha$ and $x_{0}$ are fixed. 
 For $\alpha < 1/2$,  we want to do the same, apply \eqref{eq:8.21} to $\HH_{n}(z)$ in \eqref{eq:8.24}, $w_+$ root accounts for the main contribution to the asymptotic estimate and $w_+\phi''(w_+)=2\ii(1-\zeta^2)^{1/2}$. In this case, $\zeta = (Y - \ii X)/\nu$.

As highlighted in remark \ref{msd_remark}, care must be taken to constrain $\zeta$ so that our approximation in \eqref{eq:8.21} holds uniformly as $\omega \rightarrow \infty.$ For $\alpha \neq 1/2$ and large $\omega$,
\begin{Eqnarray*}
	\zeta \sim \frac{ \omega}{\sqrt{(1-4\alpha^{2})(2n + 1)}}.
\end{Eqnarray*}
We have two competing large parameters $n$ and $\omega$, and we want to find a substitution which preserves the relationship between the two terms, 
\begin{Eqnarray}\label{eq:hermiteconstrain1}
	\sqrt{n} \propto \omega,\qquad c_{1} = \frac{n}{\omega^{2}}, \textnormal{ for some } c_{1} >0.
\end{Eqnarray}
This satisfies the condition in remark \ref{msd_remark} for $\zeta$ to remain bounded for large $n$ and $\omega$.  The additional condition in the remark for non-coalescing saddle points, $\delta < |\zeta - 1|$ is equivalent to  
\begin{Eqnarray}\label{eq:hermiteconstrain2}
	\left|\omega^{2}/n - 2\left(1- 4\alpha^{2} \right) \right|> c_{2}, \textnormal{ for some } c_{2}>0.
\end{Eqnarray}

\subsubsection{$\alpha > 1/2$: Deriving an asymptotic estimate \label{sec:derivingasymphermitea12}}
In a more practical setting, we would expect $\alpha > 1/2$ and we commence by finding the asymptotics for this case. Equation \eqref{eq:8.21} is instrumental in deriving first-order asymptotics when large degree coexists with large argument: higher order estimates can be obtained from \cite{olver1959uniform}, but we do not pursue this theme further.

\begin{figure}[htb]
  \begin{center}
      \includegraphics[width=.47\textwidth]{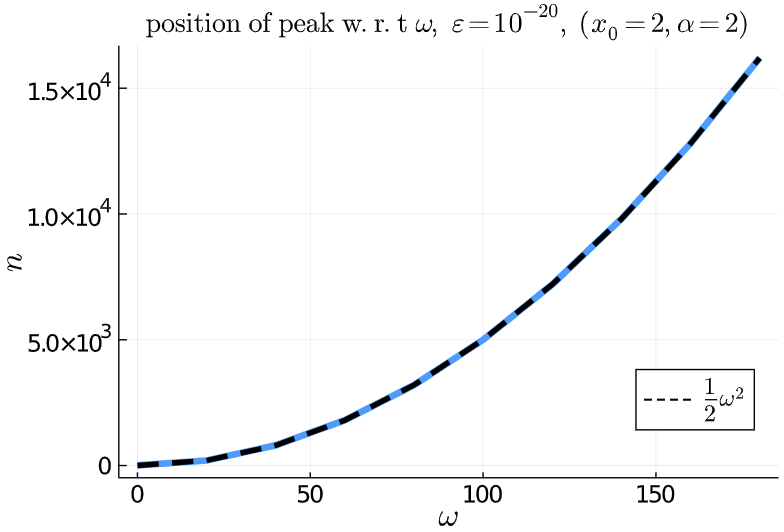}
      \includegraphics[width=.47\textwidth]{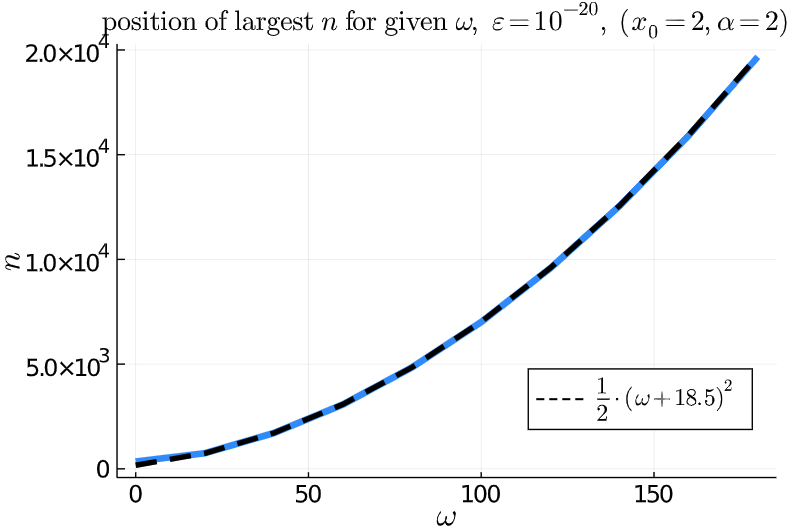}
      \includegraphics[width=.6\textwidth]{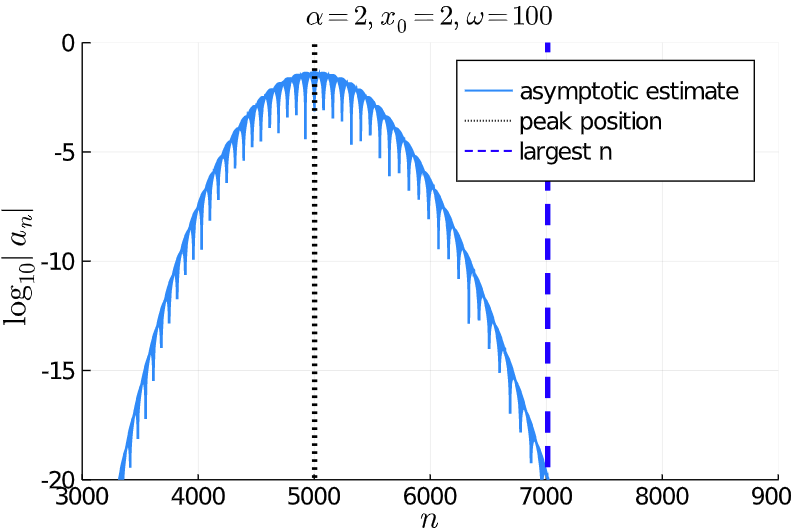}
    \end{center}
    \caption{\textit{Hermite functions:\/} Numerical experiments to determine the relationship between $n$ and $\omega$. The first figure is a plot of the position of the peak, $n$, as $\omega$ varies. The relationship follows $n = \omega^{2}/2$. The second figure is the experiment where we found the largest of coefficients such that $|a_{n}| > 10^{-20}$. This relationship follows $n = \frac{1}{2} (\omega + 18.5)^{2}$. The two experiments are consistent with $\sqrt{n} \propto \omega$ relationship. The last plot shows the peak position in black dotted line and largest $n$ in dashed line for $\alpha=2, x_{0}=2$ and $\omega=100$.\label{fig:hermitea12asymp}}
  \end{figure}

We observe that for reasonable values of $x_{0}$ we can simplify things by requiring $Y\gg X$ so $\zeta\sim \ii Y/\nu$. However, for more accurate asymptotics we want to retain $x_{0}$ explicitly in our expansions. As discussed in the last section, to satisfy remark \ref{msd_remark}, \eqref{eq:hermiteconstrain1} must hold for $\zeta$ to remain bounded for large $n$ and $\omega$.  Moreover, as $\alpha > 1/2$ we have $2(1-4\alpha^{2}) <0$, so condition \eqref{eq:hermiteconstrain2} is automatically satisfied. Thus the MSD approximation in \eqref{eq:8.21} holds uniformly for $n > \cw \omega^{2}$ for any given constant $\cw$. The relationship is demonstrated in our numerical experiments, shown in Fig.~\ref{fig:hermitea12asymp}.

Applying \eqref{eq:8.21} to \eqref{eq:8.22} (with identical notation), we have 
\begin{Eqnarray*}
  a_{n}&\sim& \textnormal{Re} \left[ \frac{\pi^{1/4}}{(2^nn!)^{1/2}} \left(\frac{2}{1+2\alpha}\right)^{\!1/2} \left|\frac{2\alpha-1}{2\alpha+1}\right|^{n/2} \exp\!\left(\frac{|1-2\alpha|}{4\alpha}(X^2+2\alpha Y^2)\right)\right.\\
  &&\left.\mbox{}\times \exp\!\left(\frac{2\ii\alpha\omega x_0}{1+2\alpha}\right) \times \frac{1}{2\pi i}  \frac{n!}{\nu^{n+1}}\exp \left( -\nu^2\phi(w_-) \right) \left(\!\frac{2\pi}{w_-\phi''(w_-)}\!\right)^{\!1/2}\right]\!\!\left[1+\mathcal{O}\left(n^{-1}\right)\right] =\textnormal{Re}\,(AB), 
\end{Eqnarray*}
where
\begin{subequations}
\begin{Eqnarray}
  A&=&\frac{1}{\ii 2^{1/2}\pi^{1/4}\nu(1+2\alpha)^{1/2}}\exp\!\left(-\frac{\alpha x_0^2}{1+2\alpha}+\frac{2\ii\alpha\omega x_0}{1+2\alpha}\right)\frac{1}{[w_-\phi''(w_-)]^{1/2}}, \label{eq:A1}\\
  B&=&\left(\frac{n!}{2^n}\right)^{\!1/2} \frac{1}{\nu^n} \left|\frac{2\alpha-1}{2\alpha+1}\right|^{n/2} \exp\!\left(\frac{|1-2\alpha|}{2}Y^2-\nu^2\phi(w_-)\right)\!. \label{eq:B1}
\end{Eqnarray}
\end{subequations}
$A$ is slowly varying for large $n$, $Y$ and $\nu = \sqrt{2 n + 1}$, while $B$ is the dominant contribution, the fast component.

At this point, the reader can refer to section \ref{sec:hermiteworking1} for the full derivation of the asymptotics. To summarise, we substitute terms $w_{-}\phi''(w_{-})$, $\phi(w_{-})$ and \eqref{eq:XandY} into \eqref{eq:A1} and \eqref{eq:B1}. 

\begin{figure}[tbh]
  \begin{center}
    \includegraphics[width=.40\textwidth]{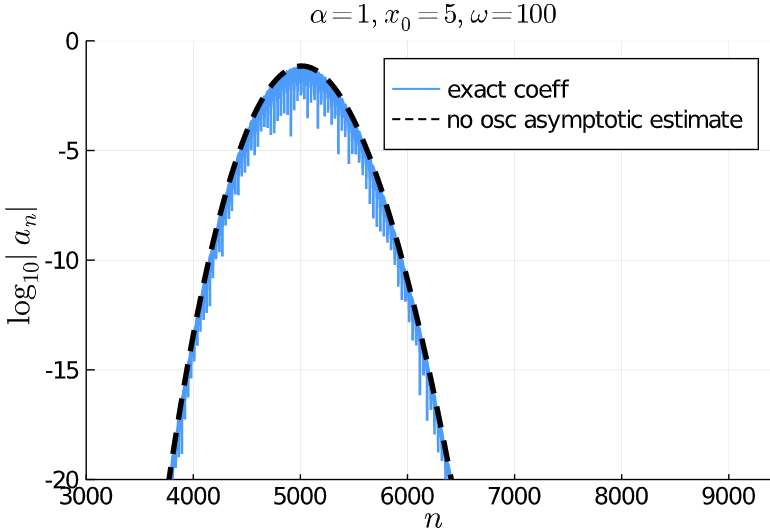}\hspace*{30pt}
    \includegraphics[width=.40\textwidth]{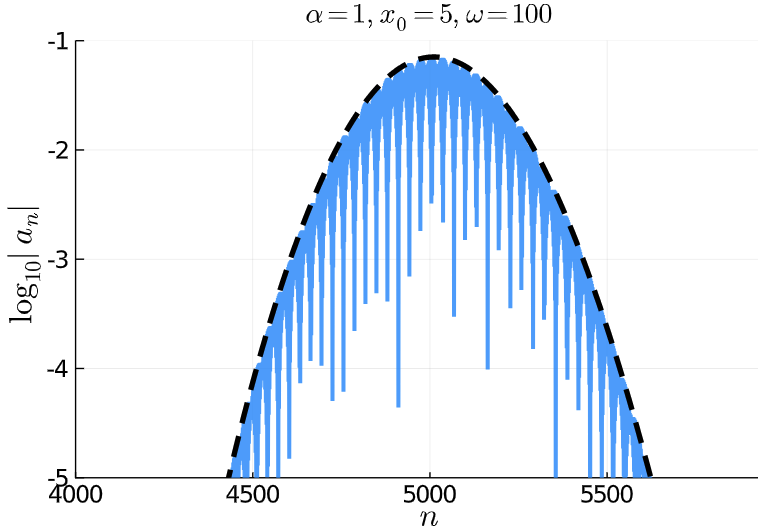}
    \includegraphics[width=.40\textwidth]{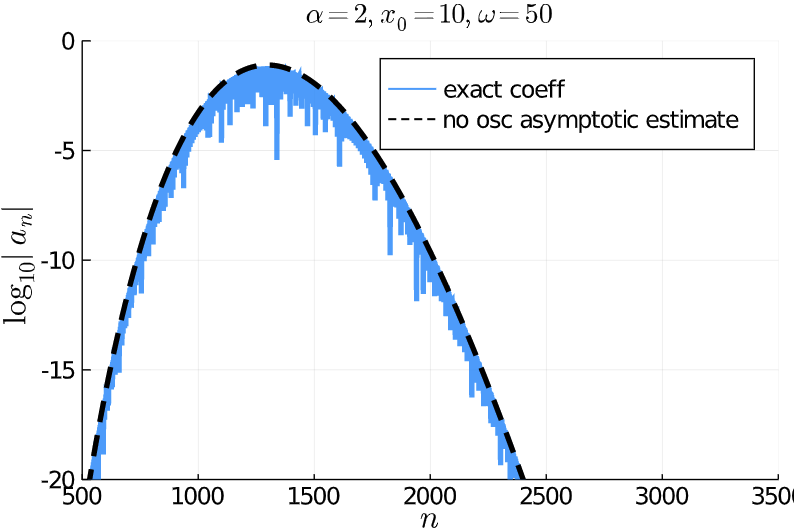}\hspace*{30pt}
    \includegraphics[width=.40\textwidth]{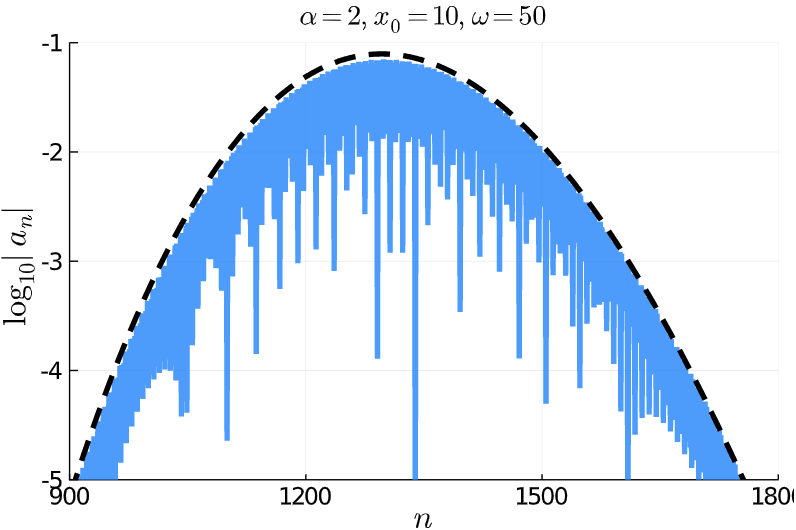}
    \includegraphics[width=.40\textwidth]{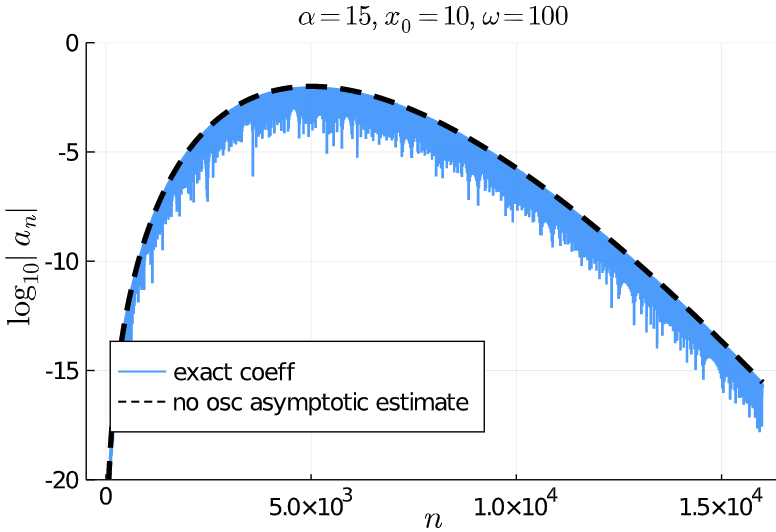}\hspace*{30pt}
   \includegraphics[width=.40\textwidth]{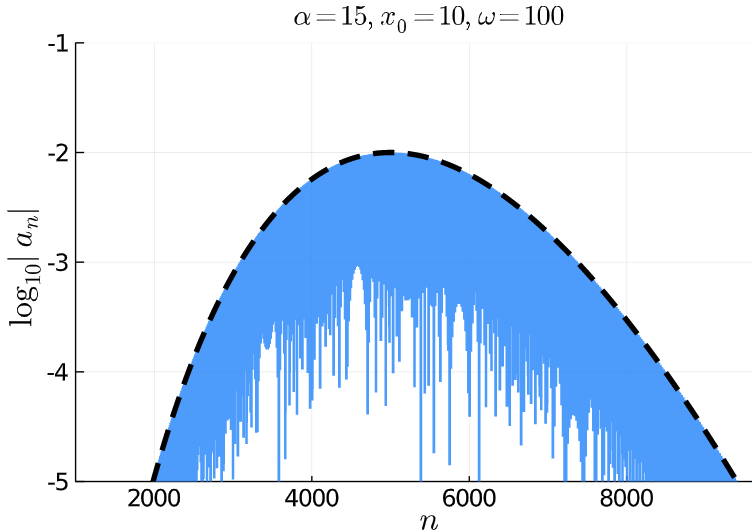}
   \end{center}
  \caption{\textit{Hermite functions:} The dashed line represents our estimate without oscillations. This line envelopes the oscillations in the lighter colour, which further confirm that our estimate provides an accurate representation of the exact coefficients. The left column of figures are plots of $\log_{10}|a_{n}|$ against $n$, where the vertical scale goes from $-20$ to 0. The figures in the right column are close-up plots to show how well we are able to bound from above the exact coefficients. \label{fig:hermiteplots}}
  \end{figure}
  
  Setting $n = \cw \omega^{2}$ where $\cw$ is some constant, then expanding about $\omega$ we arrive at 
\begin{Eqnarray*}
  A &\sim & \frac{1}{\ii \pi^{1/4} \sqrt{-2\ii (2\alpha + 1)}}\exp\left(\frac{{-\alpha x_{0}^{2} + 2\ii \omega\alpha x_{0}}}{2\alpha + 1}\right) \left[ \frac{4\alpha^{2}-1}{1+2c_{n,\omega}\left(4\alpha^{2}-1\right)}\right]^{1/4} \frac{1}{(2c_{n,\omega})^{1/4}\omega}  \\
  & &\mbox{} \times \left\{ 1 + \frac{\ii \alpha x_{0}}{\omega \left[1 + 2c_{n,\omega}\left(4\alpha	^{2}-1\right)\right]} + \mathcal{O}\left(\omega^{-2}\right) \right\} \quad \textnormal{and }\\
  B&\sim&(2\pi n)^{1/4} 2^{1/2} \ee^{\ii \pi\left(n +1/2\right)/2} \\
  & & \mbox{}\times \exp\! \left( \omega^{2}\! \left[-\frac{1}{2(2\alpha + 1)}+ c_{n,\omega} \!\left(\frac{1}{1+\sqrt{1+2c_{n,\omega}\left(4\alpha^{2}-1\right)}} + \log \left|\frac{1\!+\!\sqrt{1\!+\!2c_{n,\omega}\left(4\alpha^{2}\!-\!1\right)}}{\sqrt{2 c_{n,\omega}}\left(2\alpha\! + \!1\right)} \right| \right)\right] \right. \\
  & & \left. \quad \mbox{}- \frac{2\ii\alpha x_{0} \omega \left( -1+\sqrt{1+2c_{n,\omega}\left(4\alpha^{2}-1\right)}\right)}{4\alpha^{2}-1} + \frac{2 x_{0}^{2} \alpha^{2}}{ 4\alpha^{2} -1} \left(1 - \frac{1}{\sqrt{1 + 2c_{n,\omega} \left(4\alpha^{2} -1\right)}}\right) \right. \\
  & & \quad \left.\mbox{}+ \frac{1}{2} \log \left| \frac{1 + \sqrt{1 + 2 \cw \left(4\alpha^{2} -1\right)}}{\sqrt{2\cw \left(4\alpha^{2} -1\right)}} \right|+ \frac{\ii x_{0}\alpha \left[-3 + 6c_{n,\omega} + 8 c_{n,\omega}\left(-3 + x_{0}^{2}\right)\alpha^{2}\right]}{3\omega\left[1+2c_{n,\omega}\left(4\alpha^{2}-1\right)\right]^{3/2}} + \mathcal{O}\left(\omega^{-2}\right)\right)\!.
\end{Eqnarray*}
Substituting $A$ and $B$ back into $a_{n}$ and restoring  $c_{n,\omega} = n/\omega^{2}$, we obtain the asymptotic estimate of coefficients, 
\begin{Eqnarray}
  a_{n}(\omega) &\sim & \frac{1}{\sqrt{\omega(2\alpha +1)}} \left[ \frac{4\alpha^{2} -1}{1 + 2\frac{n}{\omega^{2}}(4\alpha^{2}-1)} \right]^{1/4} \nonumber \\
  & &\mbox{} \times \left| \frac{1 + \sqrt{1 +2\frac{n}{\omega^{2}}\left(4\alpha^{2}-1\right)}}{\sqrt{2\frac{n}{\omega^{2}}}(2\alpha +1)}\right|^{n} \left| \frac{1 +\sqrt{1+2\frac{n}{\omega^{2}}\left(4\alpha^{2}-1\right)}}{\sqrt{2\frac{n}{\omega^{2}} \left( 4\alpha^{2}-1 \right) }} \right|^{1/2} \nonumber \\
  & & \mbox{}\times \exp \!\left(- \frac{\omega^{2}}{2(2\alpha + 1)} + \frac{n}{1+\sqrt{1 +2\frac{n}{\omega^{2}}\left(4\alpha^{2}-1\right)}}   \right. \nonumber \\
  & & \left. \quad \mbox{}+  \frac{\alpha x_{0}^{2}}{4\alpha^{2}-1} \left[ 1 - \frac{2 \alpha}{\sqrt{1+2\frac{n}{\omega^{2}}\left(4\alpha^{2}-1\right)}}\right]  -\frac{\left[-1 + 2\frac{n}{\omega^{2}} + 8\frac{n}{\omega^{2}}\alpha^{2}\left(x_{0}^{2}-1\right)\right]^{2}}{16 n \left[1+ 2\frac{n}{\omega^{2}}\left(4\alpha^{2}-1\right)\right]^{5/2}} \right) \nonumber\\
  & & \mbox{}\times \left\{ \cos \left[ \frac{n \pi}{2} - \frac{2\alpha x_{0} \omega \left(-2\alpha + \sqrt{1 + 2\frac{n}{\omega^{2}}\left(4\alpha^{2}-1\right)}\right)}{4\alpha^{2}-1} \right. \right. \nonumber\\
  & & \left. \left. \quad \quad \mbox{}+ \frac{x_{0} \alpha \left[-3 + 6\frac{n}{\omega^{2}} + 8 \frac{n}{\omega^{2}}\alpha^{2}\left(x_{0}^{2} - 3\right)\right]}{3\omega\left[1 + 2 \frac{n}{\omega^{2}}\left(4\alpha^{2}-1\right)\right]^{3/2}}\right] \right.\nonumber\\
  & & \left. \quad \mbox{}- \frac{\alpha	x_{0}}{\omega \left[ 1 + 2\frac{n}{\omega^{2}}\left(4\alpha	^{2} -1\right) \right]} \sin\! \left( \frac{n \pi}{2} - \frac{2\alpha x_{0} \omega \left(-2\alpha + \sqrt{1 + 2\frac{n}{\omega^{2}}\left(4\alpha^{2}-1\right)}\right)}{4\alpha^{2}-1} \right. \right. \nonumber\\
  & & \left. \left. \quad \quad \mbox{}+ \frac{x_{0} \alpha \left[-3 + 6\frac{n}{\omega^{2}} + 8 \frac{n}{\omega^{2}}\alpha^{2}\left(x_{0}^{2} - 3\right)\right]}{3\omega\left[1 + 2 \frac{n}{\omega^{2}}\left(4\alpha^{2}-1\right)\right]^{3/2}} \right)\!\right\}\!\left[1 + \mathcal{O}\left(n^{-1}\right) \right]. \label{eq:hermiteest}
\end{Eqnarray}
Fig.~\ref{fig:hermiteplots} demonstrates how remarkably well our estimate models the size of the coefficients.

The dominant term in the Hermite estimate with the bounded oscillations (the terms in the curly brackets in \eqref{eq:hermiteest}) removed has the  form
\begin{Eqnarray}\label{eq:hermiteestnoosc}
  a_{n}(\omega) &\sim & \frac{\mathcal{C}_{osc}}{\sqrt{\omega(2\alpha+1)}} \left| \frac{1 + \sqrt{1 +2\frac{n}{\omega^{2}}\left(4\alpha^{2}-1\right)}}{\sqrt{2\frac{n}{\omega^{2}}}(2\alpha +1)}\right|^{n} \\
  & & \mbox{}\times \exp\!\left(\!-\frac{\omega^{2}}{2(2\alpha\! +\! 1)}+ \frac{n}{1\!+\!\sqrt{1\!+\!2\frac{n}{\omega^{2}}\!\left(4\alpha^{2}-1\right)}} +  \frac{\alpha x_{0}^{2}}{4\alpha^{2}\!-\!1} \!\!\left[ 1\! -\! \frac{2 \alpha}{\sqrt{1+2\frac{n}{\omega^{2}}\!\left(4\alpha^{2}\!-\!1\right)}}\right]\right)\!\!\left[1\! +\! \mathcal{O}\left(n^{-1}\right) \right]\!.  \nonumber
\end{Eqnarray}

\begin{figure}[tbh]
  \begin{center}
    \includegraphics[width=.60\textwidth]{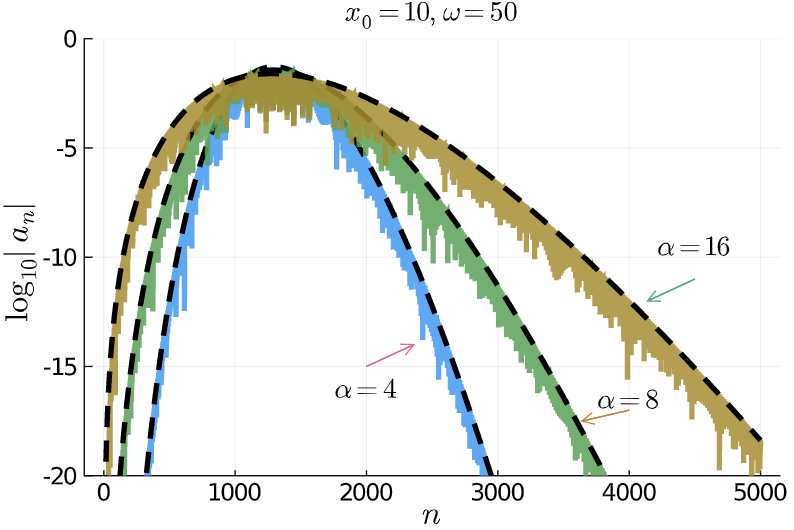}
   \end{center}
  \caption{\textit{Hermite functions:} For increasing $\alpha$, we see that the rate of decay becomes slower.  The dashed line indicate our asymptotic estimate as described by \eqref{eq:hermiteestnoosc}. Note how well our result envelopes the exact coefficients in the lighter plots.  \label{fig:hermitealphaplots}}
  \end{figure}

\subsubsection{$\alpha < 1/2$ asymptotics }\label{sec:alessthan12}
We have used the same method as for $\alpha > 1/2$ so we have omitted some technical details and the bulk of the algebra has been moved to the appendix. 

To apply the MSD, we first need to satisfy the conditions in remark \ref{msd_remark} which, in our case, translate to \eqref{eq:hermiteconstrain1} and \eqref{eq:hermiteconstrain2}. For $\alpha < 1/2$,  \eqref{eq:hermiteconstrain2} requires $\omega^{2}/n$ to be kept away from $2(1-4\alpha^{2})$.  Combining the two constraints, the MSD estimates will hold uniformly for $n > \cw \omega^{2}$ and values of $n$ outside any given neighbourhood of $\omega^{2}/(2(1-4\alpha^{2}))$.

We proceed by applying, \eqref{eq:8.21}. For $\alpha < 1/2$
\begin{Eqnarray*}
\label{eq:AN}
  a_{n}&\sim& \textnormal{Re} \left[ \frac{\pi^{1/4}}{(2^nn!)^{1/2}} \left(\frac{2}{1+2\alpha}\right)^{\!1/2} \left|\frac{2\alpha-1}{2\alpha+1}\right|^{n/2} \exp\!\left(-\frac{|1-2\alpha|}{4\alpha}(X^2+2\alpha Y^2)\right)\right.\\
  &&\left.\mbox{}\times \exp\!\left(\frac{2\ii\alpha\omega x_0}{1+2\alpha}\right) \times \frac{n!}{\nu^{n+1}} \frac{1}{2\pi\ii} \ee^{-\nu^2\phi(w_+)} \left(\frac{2\pi}{w_+\phi''(w_+)}\right)^{\!1/2}\right]\!\!\left[1+\mathcal{O}\left(n^{-1}\right)\right] =\textnormal{Re}\,(AB),
\end{Eqnarray*}
Note that for this case we apply the MSD on $\HH_{n}(\zeta)$ for $\zeta = (Y-\ii X)/\nu$ where our contour of integration passes from the saddle point $w_{+}$.
Like before, we can rewrite $a_{n}$ in terms of $A$ and $B$, the slow and fast components respectively,
\begin{Eqnarray}
	A&=&\frac{1}{\ii \pi^{1/4}\sqrt{2\ii (1+2\alpha)}}\exp\!\left(-\frac{\alpha x_0^2}{1+2\alpha}+\frac{2\ii\alpha\omega x_0}{1+2\alpha}\right)\frac{1}{\sqrt{\nu}(X^{2}+2\ii XY - Y^{2} + \nu^{2})^{1/4}}, \label{eq:AA1a}\\
	B &\sim&(2\pi n )^{1/4}2^{1/2}\sqrt{-\ii} \times \exp\!\left(-\frac{\omega^{2}}{2(2\alpha+1)} + \frac{n}{2}\log\left|\frac{2\alpha-1}{2\alpha+1}\right| \right. \label{eq:stirling}\\
	& & \left.\mbox{} +\frac{(2n + 1)(X + \ii Y)}{2 (X + \ii Y +\sqrt{X^{2}-Y^{2} + 2\ii X Y + \nu^{2}})} \right. \label{eq:BB1a} \\
	& & \left. \mbox{} + \frac{(2n +1)}{2} \log \left| \frac{X + \ii Y +\sqrt{X^{2}-Y^{2} + 2\ii X Y + \nu^{2}}}{\nu} \right|\right). \label{eq:BB1b} 
\end{Eqnarray}
Here, we have substituted $w_{+}\phi''(w_{+})$ and $\phi(w_{+})$ in terms of $X$ and $Y$, into $A$ and $B$ (see section \ref{sec:hermiteworking2}). In \eqref{eq:stirling} we have used Stirling's formula.

From equations \eqref{eq:AA1a}, \eqref{eq:BB1a} and \eqref{eq:BB1b}, the key term we need to approximate is 
\begin{Eqnarray*}
	X^{2} - Y^{2} + 2 \ii X Y  + \nu^{2}&=& 1 + \omega^{2} \left(2\cw -\frac{1}{1 - 4\alpha^{2}} \right) + \frac{4 \ii \alpha x_{0} \omega}{1-4\alpha^{2}} + \frac{4\alpha^{2} x_{0}^{2}}{1-4\alpha^{2}} .
\end{Eqnarray*}
Here we have substituted $n=\cw \omega^{2}$, like before, and we want to expand around $\omega$. First thing to note is that the large $\omega$ term is given by
\begin{Eqnarray*}
	\omega^{2}\!\left(2\cw -\frac{1}{1 - 4\alpha^{2}} \right)\!.
\end{Eqnarray*}
Observe that there is a change in behaviour at 
\begin{Eqnarray*}
	\cw^{*} = \frac{1}{2\left(1-4\alpha^{2}\right)}.
\end{Eqnarray*}
This tells us that for a small neighbourhood around $\cw^{*}$, we would expect this region to have different asymptotics to the the region away from it. This agrees nicely with our constraints discussed at the start of this section,  where we have identified that the non-coalescing saddle point regime in remark \ref{msd_remark} is equivalent to considering $\cw$ away from $\cw^{*}$.
%\begin{Eqnarray*}
%\cw \neq \cw^{*}: & & \textnormal{the large parameter is } \omega^{2}\!\left(2 \cw - \frac{1}{1-4\alpha^{2}} \right)\\
%\cw = \cw^{*}: & & \textnormal{we substitute the critical value and expand around the large $\omega$ term.} 
%\end{Eqnarray*}

Going through the same steps as in the previous section and sparing the reader more algebra (see section \ref{sec:neqcstar} for  full derivation), for $\cw \neq \cw^{*}$ we obtain
\begin{Eqnarray*}
	a_{n}(\omega )= \textnormal{Re}(AB) &\sim& \textnormal{Re} \left\{-\frac{1}{\sqrt{\omega(1+2\alpha)}}\left| \frac{2\alpha -1}{2 \alpha + 1}\right|^{n/2} \left(\frac{|1 - 4\alpha^{2}|}{2 \frac{n}{\omega^{2}} |1-4\alpha^{2}| -1}\right)^{1/4}  \right. \\
	& &\times  \left(1 - \frac{\ii \alpha x_{0}}{\omega \left(2\frac{n}{\omega^{2}} |1 - 4\alpha^{2}| -1\right)} + \mathcal{O} \left(\omega^{-2}\right)\right)  \left( \frac{1+ \sqrt{1-2\frac{n}{\omega^{2}} |1-4\alpha^{2} |}}{\sqrt{2 \frac{n}{\omega^{2}} |1-4\alpha^{2}|}}\right)^{n+1/2}\\
		& & \times \exp\left(\frac{\alpha}{|1-4\alpha^{2}|}\left( \omega^{2}- x_{0}^{2} +2 x_{0} \omega \sqrt{2 \frac{n}{\omega^{2}}|1-4\alpha^{2}| -1} \right) \right. \\
		& & + \frac{\alpha x_{0} \left(8 \alpha^{2} x_{0}^{2}\frac{n}{\omega^{2}} + 6 \frac{n}{\omega^{2}} |1-4\alpha^{2}| - 3\right)}{3\omega \left(2\frac{n}{\omega^{2}}|1-4\alpha^{2}|-1\right)^{3/2}} \\
		& & +  \ii \left[\left(n + \frac{1}{2} \right)\frac{\pi}{2} - \frac{ 4 \alpha^{2} x_{0} \omega}{|1-4\alpha^{2}|} + \frac{\omega^{2}}{2 |1-4\alpha^{2}|}\sqrt{2 \frac{n}{\omega^{2}} |1-4\alpha^{2}| -1}  \right. \\
		& & \left.\left. \left. + \frac{2 \alpha^{2} x_{0}^{2}}{|1-4\alpha^{2}|\sqrt{2\frac{n}{\omega^{2}} |1-4\alpha^{2}|-1}}\right]\right)\!\right\}\!\left[1+\mathcal{O}\left(n^{-1}\right)\right] .
\end{Eqnarray*}

%For $\cw = \cw^{*}$ (see section \ref{sec:eqcstar} for full derivation), expanding around $\omega$,
%\begin{Eqnarray}
%	a_{n}(\omega) &\sim& \textnormal{Re} \left\{ - \frac{(2n)^{1/4}}{\sqrt{1+2\alpha}} \left| \frac{1 - 2 \alpha}{1+ 2 \alpha} \right|^{n/2} \left(\frac{\ii^{7/4} |1-4\alpha^{2}|}{\sqrt{2} \left(\alpha x_{0}\right)^{1/4} \omega^{3/4}}\right) \left(1 - \frac{4\alpha^{2}x_{0}^{2} + |1-4\alpha^{2}|}{\ii 16 \alpha x_{0} \omega }\right) \right. \label{eq:alessthancritc}\\
%	& & \left. \times \exp \left( \frac{\ii \pi }{4} - \frac{\alpha x_{0}^{2}}{|1-4\alpha^{2}|}+ \frac{(\ii \pi + 4 \alpha)\omega^{2}}{4 |1-4\alpha^{2}|} - \frac{4 \ii \alpha^{2} x_{0} \omega}{|1-4\alpha^{2}|} + \frac{8 \ii^{1/2} \left(\alpha x_{0}\right)^{3/2} \omega^{1/2}}{3 |1-4\alpha^{2}|}\right. \right.\nonumber \\
%	& & \left. \left. - \frac{\ii^{3/2}}{5 \sqrt{\omega}}  \left[5 \sqrt{\alpha x_{0}} + \frac{4\left(\alpha x_{0}\right)^{5/2}}{|1-4\alpha^{2}|} \right] \right)\!\right\}\!\left[1+\mathcal{O}\left(n^{-1}\right)\right]. \nonumber
%\end{Eqnarray}

\begin{figure}[tbh]
  \begin{center}
   	\includegraphics[width=.40\textwidth]{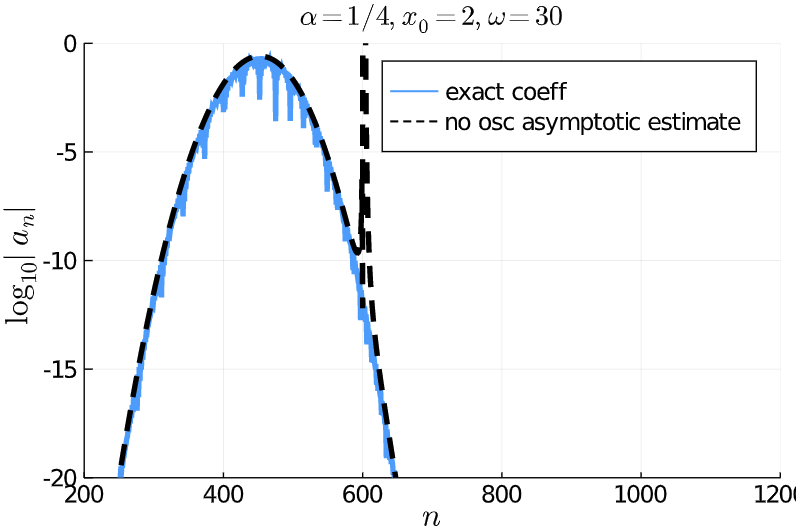}\hspace*{30pt}
      \includegraphics[width=.40\textwidth]{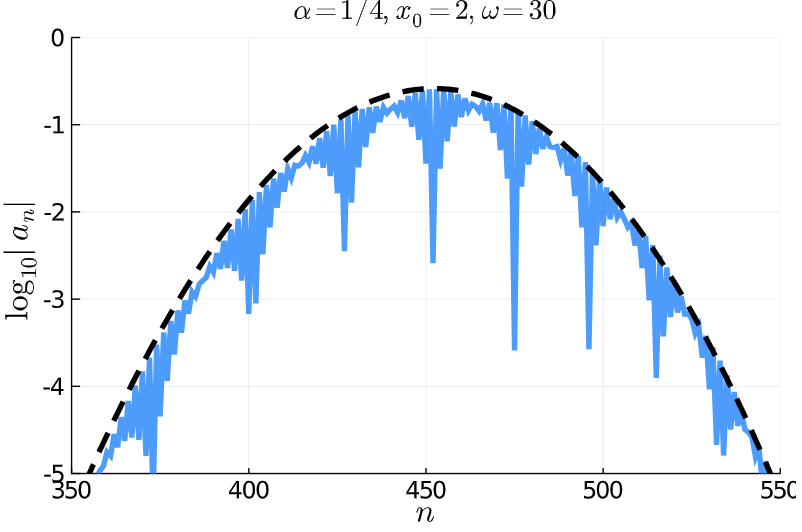}
 	 \includegraphics[width=.40\textwidth]{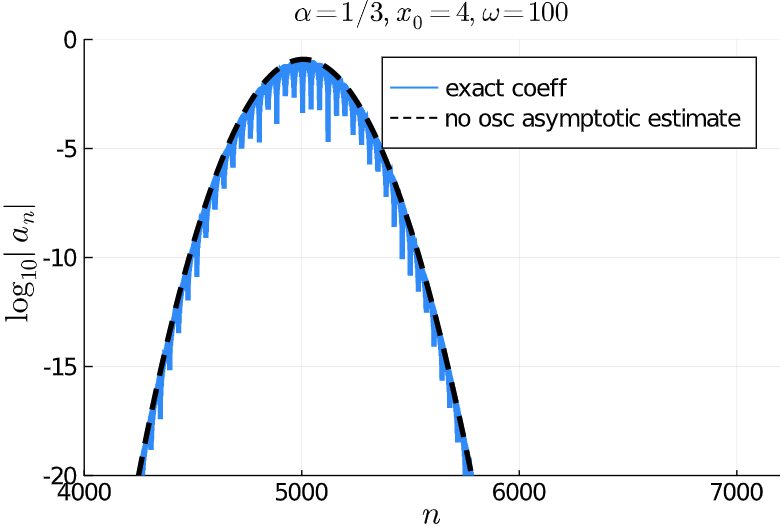}\hspace*{30pt}
      \includegraphics[width=.40\textwidth]{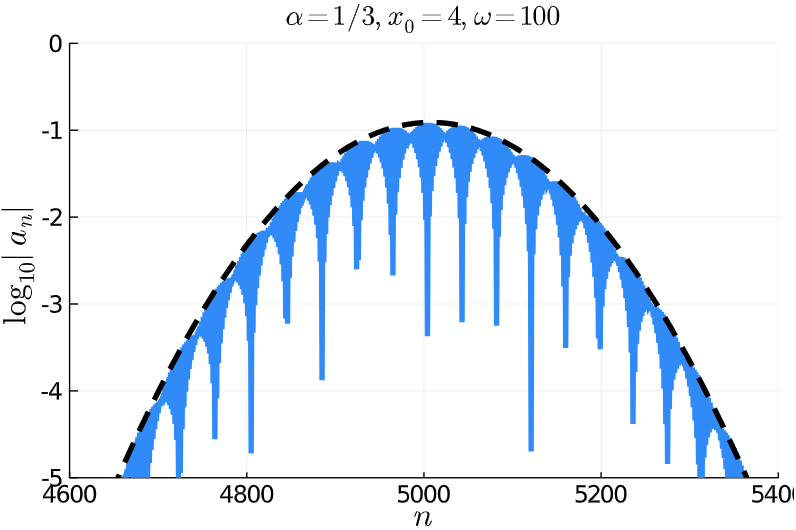}
       \includegraphics[width=.40\textwidth]{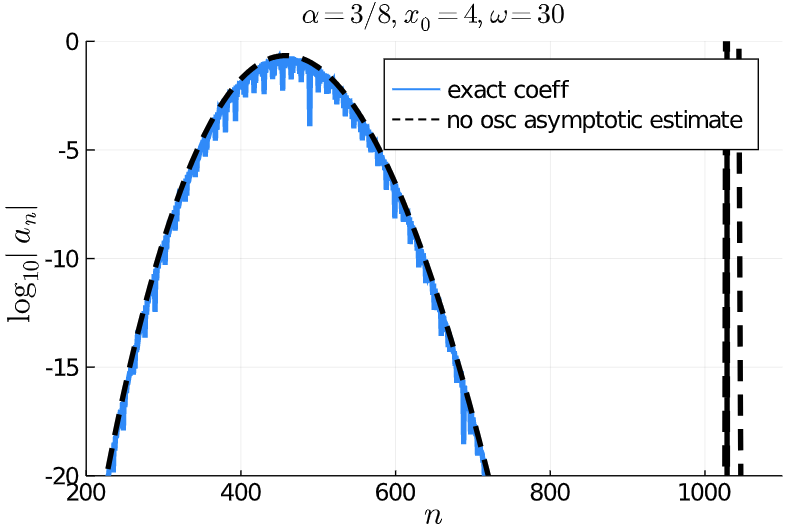}\hspace*{30pt}
       \includegraphics[width=.40\textwidth]{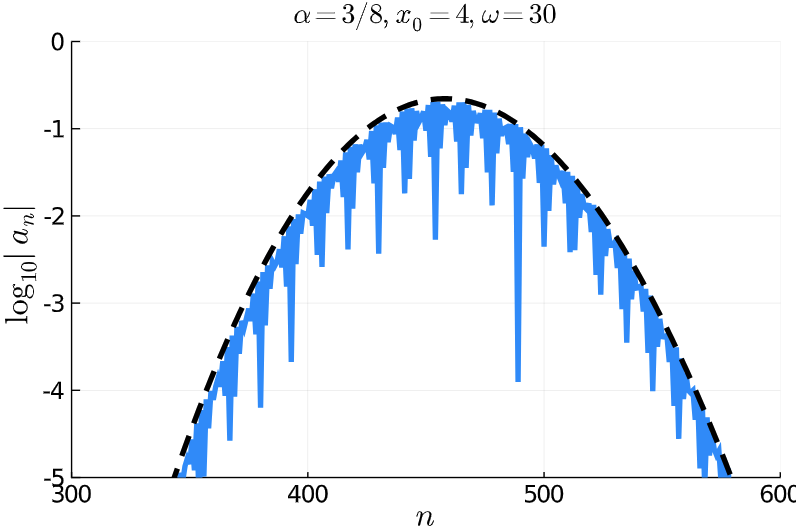}
    \end{center}
    \caption{\textit{Hermite functions:} The figures on the left show the decay rate of the coefficients where $\log_{10}|a_{n}| > 10^{-20}$. In the top row we have $\alpha = 1/4, x_{0}=2$ and $\omega = 30$. In the first figure (on the left) we can see that our estimate starts to blow up around $n = 600$, this is expected as $\cw^{*}$ for these parameters occurs at $n=600$. The plots in the right column are a closeup of the figures to the left. The plots in the second and third row, further show that our asymptotic estimates are accurate for varying $\alpha, x_{0}$ and $\omega$. Note, for these figures the $\cw^{*}$ occurs outside of the plotted region of $\log_{10}|a_{n}| > 10^{-20}$. \label{YetAnotherFigure}}
  \end{figure}

The asymptotics are displayed in Fig.~\ref{YetAnotherFigure}, which confirms the remarkable fit of our formul\ae{} with computed values of the coefficients. 

\setcounter{equation}{0}
\setcounter{figure}{0}
\section{Malmquist--Takenaka functions}\label{sec:mt}

\begin{figure}[htb]
  \begin{center}
    \includegraphics[width=.43\textwidth]{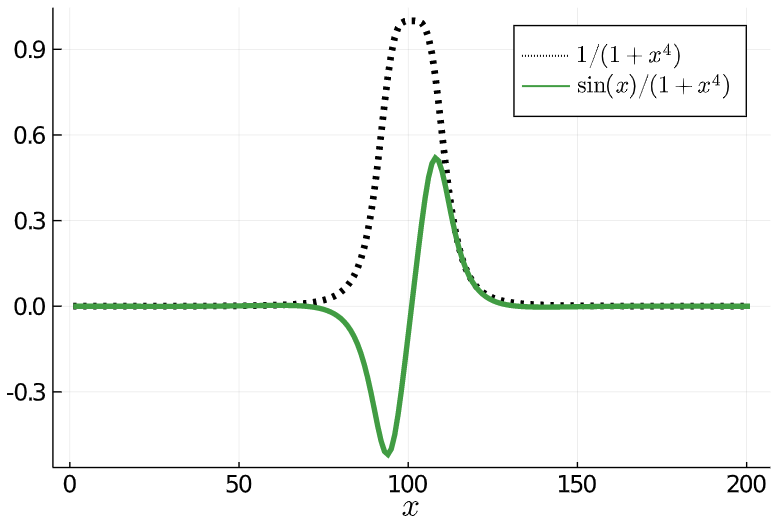}\hspace*{20pt}
     \includegraphics[width=.43\textwidth]{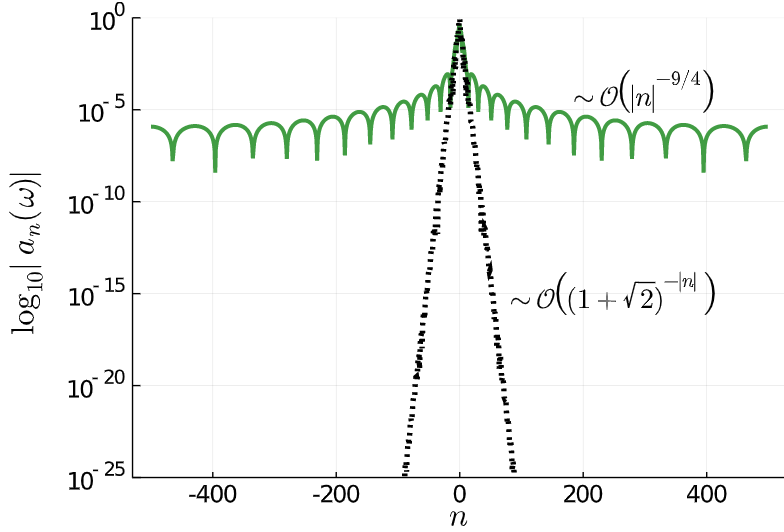}
   \end{center}
   \caption{\textit{MT functions:} The left figure shows $f(x)=\sin x/(1+x^4)$ (solid green line) and $f(x)=1/(1+x^4)$ (dotted black line). The right figure shows the corresponding logarithm of the absolute value of MT coefficients for $\sin x/(1+x^4)$ (solid green line) and $1/(1+x^4)$ (dotted black line).}
  \label{fig:algfuncompare}
\end{figure}

MT functions have a long history. Discovered simultaneously by Malmquist \cite{malmquist1925determination} and Takenaka \cite{takenaka1925orthogonal} in 1925, they have since been studied in fields ranging from signal processing (Wiener in 1949  \cite{WienerNorbert1949Eias}) to spectral methods (Christov in 1982 \cite{christov1982complete}). In \cite{iserles2019family}, Iserles and Webb proved that this is essentially the only orthonormal and complete system in $\LL_{2}(\R)$ with a banded tridiagonal skew-Hermitian differentiation matrix, and the coefficients can be computed using the Fourier transform.  This is part of Iserles and Webb's ongoing work classifying bases with banded skew-Hermitian differentiation matrices \cite{iserles2019orthogonal, iserles2020differential,iserles2019fast}. It was the work of Boyd in \cite{boyd1987spectral}  and Weideman in \cite{weideman1994computation,weideman1995computing} that highlighted the interesting approximation properties of MT functions. The approximating characteristics of MT functions do not obey simple rules. As an example, for $f(x)= 1/(1 + x^{4})$ the coefficients decay at a spectral rate of $a_{n} \sim \mathcal{O}\left((1 + \sqrt{2})^{-|n|}\right)$, however if we include an oscillatory factor so that our function is in the form $f(x) = \sin(x)/(1 + x^{4})$, our approximation collapses from spectral to algebraic convergence, $a_{n}\sim \mathcal{O}\left( |n|^{-9/4} \right)$, see Fig.~\ref{fig:algfuncompare}. Why does multiplying by $\sin(x)$, an entire function, have such a large effect on the rate of convergence? How do we get exponential convergence? 

In \cite{boyd1987spectral} and \cite{weideman1995computing} (respectively), Boyd and Weideman proved that MT expansion coefficients decay exponentially if and only if $f(z)$ decays at least linearly as $|z| \to \infty$ in the complex plane, and $f$ is analytic in a region which may exclude $\pm \ii/2$ but which must include the point at infinity. Analyticity at infinity is a strong condition to impose, since $\sin z$, for example, blows up as $z\rightarrow \infty$ along the imaginary axis. This is not unexpected because the proof of exponential convergence in a compact interval (e.g.~with orthogonal polynomials)
depends on functions being analytic within a Bernstein ellipse surrounding this interval. 

\begin{figure}[htb]
  \begin{center}
     \includegraphics[width=.60\textwidth]{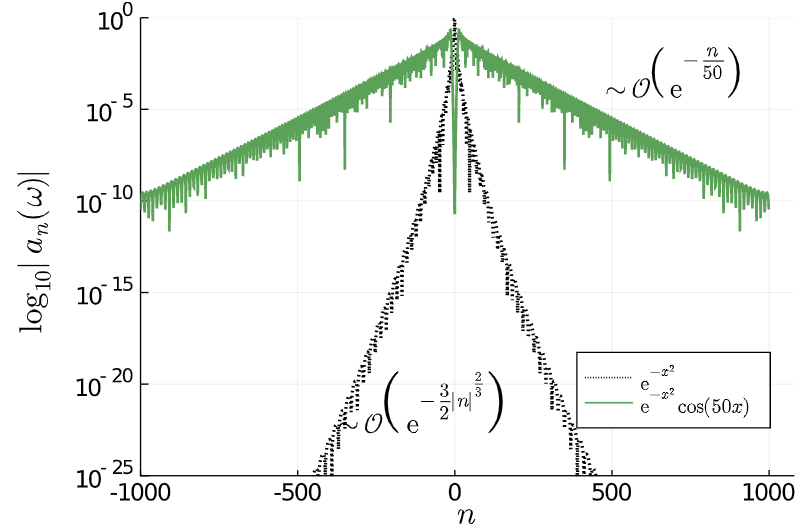}
   \end{center}
  \caption{\textit{MT functions:} The plot compares the decay rate of the MT coefficients of the Gaussian (dotted black line) and the wave packet (solid green line). Note the slower decay rate of the wave packets relative to the decay rate of the Gaussian coefficients.}
  \label{fig:expfuncompare}
\end{figure}

Interestingly, Weideman (in \cite{weideman1994computation}) found that the coefficients, $a_{n}$, of a Gaussian function asymptotically decay close to exponential convergence, $a_{n} \sim \mathcal{O}\!\left(\ee^{-\frac{3}{2}|n|^{2/3}}\right)$. Note, we do not have exact spectral decay as the Gaussian function has an essential singularity at $\infty$. What does this mean for wave packet coefficients? We proved that for wave packet functions the coefficients exhibit asymptotic exponential convergence. In other words as $\omega$ gets larger, we get closer to exponential convergence up to some threshold. For example, if $f(x) = \ee^{-x^{2}}\cos(\omega x)$ the coefficients decay like $a_{n} \sim \mathcal{O}\left(\ee^{-n/\omega}\right)$.This can be seen in Figure \ref{fig:expfuncompare}, where even though the $\cos(50 x)$ factor decreases the rate of convergence compared to just the Gaussian function,  we still see asymptotic spectral decay, at least in the regime $|a_{n}(\omega)| \geq 10^{-10}$.

To see how we have proved this, let us start by defining our system. From \cite{iserles2019family} the Malmquist--Takenaka system is defined as
\begin{equation}
  \varphi_{n}(z) := \ii^{n}\sqrt{\frac{2}{\pi}}\frac{(1 + 2\ii z)^{n}}{(1 - 2\ii z )^{n+1}}, \quad n \in \Z. \label{eq:mtfunc}
\end{equation}
Substituting  \eqref{eq:mtfunc} into  \eqref{eq:generala} yields the expansion coefficients 
\begin{Eqnarray}
  a_{n}(\omega) &=& (-\ii)^{n}\sqrt{\frac{2}{\pi}}\int_{-\infty}^{\infty} \ee^{-\alpha(x- x_{0})^{2}}\cos(\omega x) \left(\frac{1-2\ii x}{1 +2 \ii x}\right)^{\!n}\!\frac{\dd x}{1 +2 \ii x}, \label{eq:anmtfunc}\\
  &=& \frac{(-\ii)^{n}}{\sqrt{2\pi}}\int_{-\infty}^{\infty} \left(\ee^{\ii \omega x} + \ee^{-\ii \omega x} \right)  \ee^{-\alpha(x- x_{0})^{2} -2 \ii n \tan^{-1}(2x)} \frac{\dd x}{1 + 2 \ii x}, \nonumber \\
  & =& b_{n}(\omega) + b_{n}(-\omega), \nonumber
\end{Eqnarray}

where 
\begin{Eqnarray}
  b_{n}(\omega) &=&\frac{(-\ii)^{n}}{\sqrt{2\pi}}\int_{-\infty}^{\infty} \ee^{- \omega g(x)}\frac{\dd x}{1 + 2 \ii x}, \label{eq:breint} \\
g(x)&=& - \ii x + \beta (x- x_{0})^{2} + 2\ii \cw \tan^{-1}(2 x), \label{eq:gtanlog}
\end{Eqnarray}
where $\alpha = \beta \omega, n = \cw \omega$ and we take the principal value of $\tan^{-1}$. We initially restrict ourselves to the case $n \geq 0, \ \omega\geq 0$. The remaining cases will be considered in section \ref{sec:reduction}.

\subsection{The case $\omega \geq 0$, $n \geq0$: the method of steepest descent}

This section draws inspiration from Weideman in \cite{weideman1994computation}, who found the rate of decay of MT coefficients for the Gaussian function, $\ee^{-x^{2}}$ by applying MSD to the coefficient integral. Our analysis will be considerably more involved than that of Weideman’s for the Gaussian function, because we have two large parameters, $n$ and $\omega$. A brief explanation of the MSD can be found in section \ref{sec:msdexplain}.

Since $\beta=\alpha/\omega$, it is bounded for large $\omega$. We  stipulate that $n$ and $\omega$ vary so that $\cw$ is bounded: details will be discussed further in the next section.The main step of the MSD is to deform our integral in \eqref{eq:breint} over a complex contour that passes through (or near) a saddle point of \eqref{eq:gtanlog}. This represents the path of steepest descent. We can find the saddle points by solving
\begin{Eqnarray}
	g'(z) = - \ii +  2 \beta (z- x_{0}) + \frac{4 \ii \cw}{4 z^{2} +1} =0, \nonumber
\end{Eqnarray}
with some manipulation i.e. multiplying through by $1+ 4z^{2}$ then collecting $z$ terms,  this reduces to the cubic
\begin{Eqnarray}\label{eq:note03}
	8 \beta z^{3} - 4 \left( \ii + 2 \beta x_{0}\right)z^{2} + 2 \beta z + \ii \left(4 \cw -1\right) -2 x_{0} \beta = 0
\end{Eqnarray}
Our formulation in \eqref{eq:breint} has a nice symmetry aspect that comes from the $\tan^{-1}$. For example, when $\beta=0$ there are two saddle points
\begin{equation}\label{eq:note04}
z_-=-\sqrt{\cw-\tfrac14},\quad z_+=\sqrt{\cw-\tfrac14}.
\end{equation}
As $\alpha$ is not an asymptotic parameter, $\beta$ will be small for large $\omega$, and the two shown values $z_\pm$ in \eqref{eq:note04} will be close to two exact saddle points, while the third one  approaches infinity when $\beta\to 0$.

We  have for two saddle points $z_2$ and  $z_3$ expansions
\begin{equation}\label{eq:note05}
z_2=z_++c^{(2)}_1\beta+c^{(2)}_2\beta^2+\ldots,\quad z_3=z_-+c^{(3)}_1\beta+c^{(3)}_2\beta^2+\ldots,
\end{equation}
where $\cw$ should be bounded away from $\frac14$ and 
\begin{Eqnarray*}
c^{(*)}_1&=&\frac{\ii\cw(x_0-z_*)}{z_*}, \\
c^{(*)}_{2}&=& - \frac{\cw (x_{0} - z_*)}{2 z_*^{3}} \left[ -  \cw x_{0} + z_* \left( 4 x_0 z_* - 5\cw + 1 \right) \right].
\end{Eqnarray*}
We obtain this from substituting \eqref{eq:note05} into \eqref{eq:note03}, then collecting $\beta$ terms to find the $c_{k}^{(*)}$s, the coefficients $c_{k}$s at saddle point $z_{*} \in \{z_{1},z_{2},z_{3}\}$. For the third saddle point, $z_1$, we use a property of the cubic equation. Writing the equation as $a_3z^3+a_2z^2+a_1z+a_0=0$,
\begin{equation*}
z_1+z_2+z_3=-\frac{a_2}{a_3}= x_0+\frac{\ii}{2\beta},
\end{equation*}
and $z_3$ has the expansion
\begin{equation}\label{eq:note08}
z_1=-\frac{a_2}{a_3}-z_2-z_3=x_0+\frac{\ii}{2\beta}+2\ii\cw\beta- \left(c^{(2)}_2+c^{(3)}_2\right)\beta^2-\dots\,.
\end{equation}

The two saddle points $z_2$ and $z_3$ become nearly equal when $\cw\to\frac14$. That is,  there is a special case, when $\omega\sim4n$.  In this case, we cannot use the MSD and will have to use uniform asymptotic methods -- this method can be found in \cite{olver1997asymptotics,temme2014asymptotic, wong2001asymptotic}. We do not consider this case in the paper and assume that $\cw \geq \frac{1}{4} + \delta$ where $\delta > 0$ is small.  This is a necessary constraint to satisfy remark \ref{msd_remark} and we discuss this further in the next section.

\subsection{Large $\omega$ asymptotics}

Our main goal is to find out how many MT functions are required to approximate a wave packet with a frequency of $\omega$ to prescribed accuracy. We want to estimate the number of $\{a_{n}: |a_{n}| = |b_{n}(\omega) + b_{n}(-\omega)| \geq \varepsilon \}$ for any given $\varepsilon > 0$. Our asymptotic expansion requires a subtle approach as both $n$ and $\omega$ are large parameters. 

We have learnt from numerical experiments that it is a poor idea to fix one parameter and let the other become large i.e.~fix $\omega$ and let $n \gg 1$, as the asymptotics from this regime miss on the spectral decay behaviour that we wish to capture. A simple explanation for this alludes to the relationship between $n$ and $\omega$ in the square-root term in equation \eqref{eq:note04}. 

We found that there was an intermediate range where $n$ and $\omega$ followed a relationship of the form
\begin{Eqnarray*}
	\cw < n < c \qquad \cw = \frac{n}{\omega} > \frac{1}{4},
\end{Eqnarray*}
for any given constant $c$. This range satisfies remark \ref{msd_remark} and ensures that the MSD approximation holds uniformly for $\omega \rightarrow \infty$.

Writing \eqref{eq:note05} and \eqref{eq:note08} in orders of $\omega$,  i.e. substituting $\beta = \alpha/ \omega$, the roots have the form
\begin{Eqnarray}
	z_{1} &\sim & x_0+\frac{\ii \omega }{2\alpha}+\frac{2\ii\cw \alpha}{\omega} + \mathcal{O}\!\left(\omega^{-2}\right)\,\!,\label{eq:z1est} \\
	z_{2} &\sim & \sqrt{\cw-\tfrac14} +  \frac{\ii\alpha \cw(x_0-\sqrt{\cw-\tfrac14})}{ \omega\sqrt{\cw-\tfrac14}}  + \mathcal{O}\!\left(\omega^{-2}\right),\label{eq:z2est}\\
	z_{3} & \sim & -\sqrt{\cw-\tfrac14} -  \frac{\ii \alpha\cw(x_0+\sqrt{\cw-\tfrac14})}{\omega\sqrt{\cw-\tfrac14}} + \mathcal{O}\!\left(\omega^{-2}\right)\label{eq:z3est}.
\end{Eqnarray}

We require $\cw > \frac{1}{4}$ so $\sqrt{\cw - \frac{1}{4}}$ is always  positive. For $\alpha, x_{0}, n, \omega \geq 0$, $z_{3}$ stays in the third quadrant, $z_{2}$ can be located in the first or fourth quadrant and $z_{1}$ is located in first quadrant of the complex plane. Moreover,  if $x_{0} =0$, then $z_{2}$ would only be in the fourth quadrant. From $z_{1}$, we see that  the root is shifted from the pure imaginary axis by $x_{0}$.  By requiring $\textnormal{Re}\,z_{2} > x_{0}$,  we get $\cw>\frac{1}{4} + x_{0}^{2}$ and past this point we see that $\textnormal{Im}\,z_{2} < 0$, so in this regime $z_{2}$ is confined to the fourth quadrant.

The plot below suggests that $z_2$ and $z_3$ may be appropriate saddle points through which to deform the contour. 

\begin{figure}[tbh]
  \begin{center}
  	\includegraphics[width=.45\textwidth]{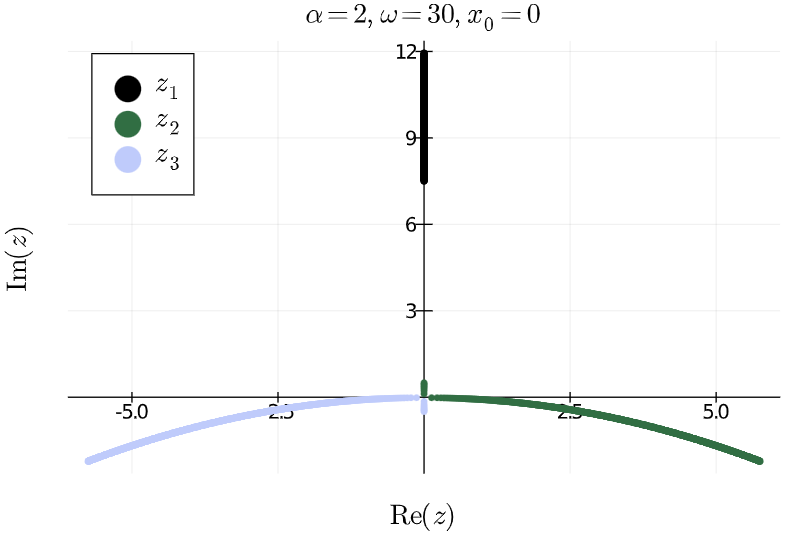}
    \includegraphics[width=.45\textwidth]{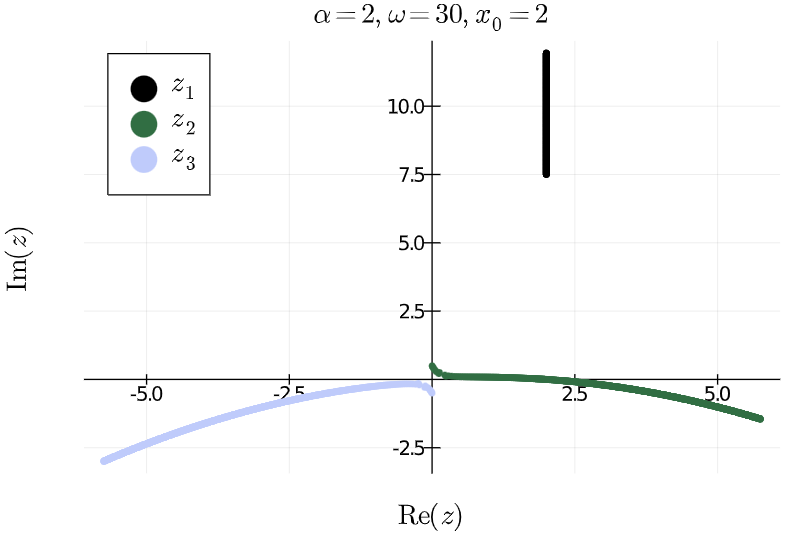}
   \end{center}
  \caption{\textit{MT functions:} Plot of the roots  \eqref{eq:z1est}, \eqref{eq:z2est} and \eqref{eq:z3est} for varying $n$. \textit{Left plot:} Plot of saddle points when $x_{0}=0$. All the saddle points are on the imaginary line when $\cw < \frac{1}{4}$,  when $\cw > \frac{1}{4}$, $z_{2}$ is the reflection of $z_{3}$ on the $y$ axis.  \textit{Right plot:} Plot of saddle points when $x_{0}=2$.  $z_{1}$ is shifted $x_{0} = 2$ away from the imaginary axis and $z_{2}$ travels from the first quadrant to the fourth quadrant as $n$ increases.}
\end{figure}

\subsubsection{Paths of steepest descent}
\begin{figure}[tbh]
  \begin{center}
    \includegraphics[width=0.48\textwidth]{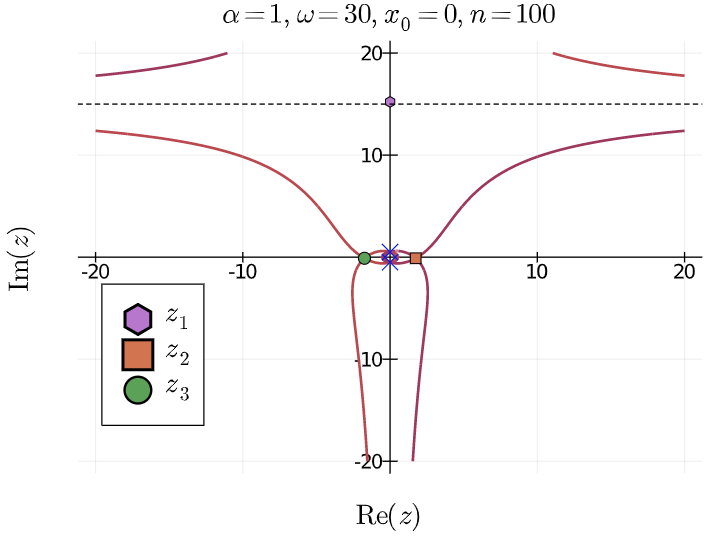}
    \includegraphics[width=0.48\textwidth]{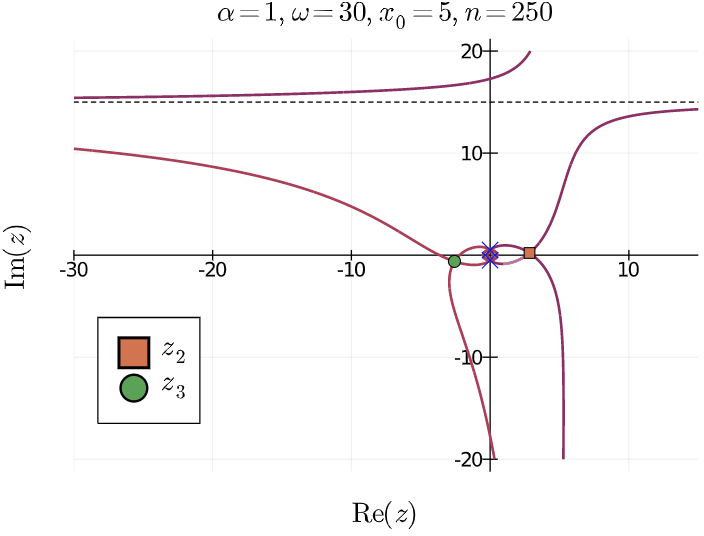}
    \includegraphics[width=0.60\textwidth]{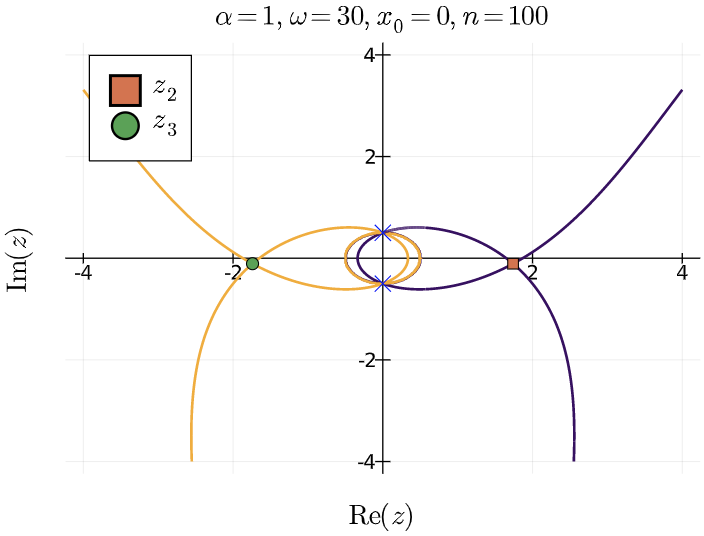}
    \end{center}
  \caption{\textit{MT functions:} Plot of  \eqref{eq:pathofsteeptestdescent} for given $\alpha, \omega$, $x_{0}$ and $n$. The $\times$s on the imaginary axis represent the pole at $\ii/2$ and zero at $-\ii/2$ for $n \geq 0$. It can be seen in the top plots that the contour approaches an asymptote that is independent of $x_{0}$. In the bottom plot, we see a zoomed in plot of the top left plot, the steepest descent contours pass through the pole at $\ii/2$ and zero at $-\ii/2$ for $n \geq 0$.}\label{fig:contourplot}
\end{figure}

The path of steepest descent is chosen so that $\textnormal{Im}\,g(z)$ is constant and passes through the saddle point $z^{*}$ of $g(z)$.  Recall from  equation~\eqref{eq:gtanlog} that  
\begin{Eqnarray*}
	g(z) = - \ii z + \beta (z- x_{0})^{2} + 2\ii \cw \tan^{-1}(2 z).
\end{Eqnarray*}
Substituting $z = x+ \ii y$, and using the identity 
\begin{Eqnarray*}
	\tan^{-1} \!\left(a + \ii b\right) &=& \frac{1}{2} \tan^{-1} \!\left(\frac{2 a }{1 - (a^{2} + b^{2})}\right) + \frac{\ii}{2} \tanh^{-1}\! \left(\frac{2 b}{1 + a^{2}+b^{2}}\right),
\end{Eqnarray*}
the imaginary part of $g(z)$ has the form
\begin{Eqnarray} \label{eq:pathofsteeptestdescent}
	\textnormal{Im}g(x + \ii y) &=& -x +  2\beta y (x - x_{0})  + \cw\tan^{-1}\! \left(\frac{4x}{1 - 4x^{2} - 4y^{2}}\right)\!.
\end{Eqnarray}
We choose the path of steepest descent to be the path where $\textnormal{Im}\,g(z)=\textnormal{constant}$,  the constant being chosen so that the path passes through the saddle point $z^{*}$ i.e.\ path $= \{z: \textnormal{Im}\,g(z) = \textnormal{Im}\,g(z^{*}) \}$.  Figure \ref{fig:contourplot},  shows a plot of steepest descent paths, equation~\eqref{eq:pathofsteeptestdescent}, for certain parameter values.

To construct our steepest descent contour $\Gamma$, we split $\Gamma = \Gamma^{(3)} \cup \Gamma^{(2)}$, where $\Gamma^{(2)}$ is the contour that passes through $z_{2}$, i.e.~the integral that goes from $0$ to $\infty$.  Note that for $n>0$ \eqref{eq:anmtfunc} has a pole at $\ii/2$ and a zero at $-\ii/2$ and it is sensible to avoid the pole in the upper half plane. From our discussion on the location of our roots in the last section,  $z_{2}$ can be in the first or fourth quadrant. For $\cw > \frac{1}{4} + x_{0}^{2}$,  $z_{2}$ is in the fourth quadrant. For simplicity we will consider the root to be located in the fourth quadrant, but the same principle can be applied when the root is in the first quadrant. We start by defining contour $\Gamma^{(2)}$ at $-\ii/2$:
\begin{itemize}
\item This was chosen simply because we want to avoid the pole at $\ii/2$.  We can show that $-\ii/2$ is always on our steepest descent path by substituting $-\ii/2$ into \eqref{eq:pathofsteeptestdescent} which gives $\textnormal{Im}\,g(-\ii/2) = \beta x_{0}= \alpha x_{0}/\omega$ where $\alpha, x_{0}, \omega > 0$ with no dependence on $n$.  For a visual representation see the bottom plot of Fig.~\ref{fig:contourplot},  where the steepest descent contour through $z_{2}$ passes through both the pole and the zero. We can simply choose to construct our contour from $-\ii/2$ as there is nothing to look for along the path from 0 to $-\ii/2$.  Moreover, the point $z = -\ii/2$ is the central point of a valley which is why the path of steepest descent may visit this point.
\item By starting our contour from $-\ii/2$ means our curve looks nice and smooth. The imaginary part converges to $\omega/2\alpha$ whilst the real part tends to infinity. As the real part of our contour tends to infinity, the integrand tends to zero very fast.  Although this step is not necessary,  we can bring the tail back to the real line by travelling straight down, only varying the imaginary part of our path, to reach the real line and then carrying on towards infinity. 
\end{itemize}

\begin{figure}[tbh]
  \begin{center}
    \includegraphics[width=0.75\textwidth]{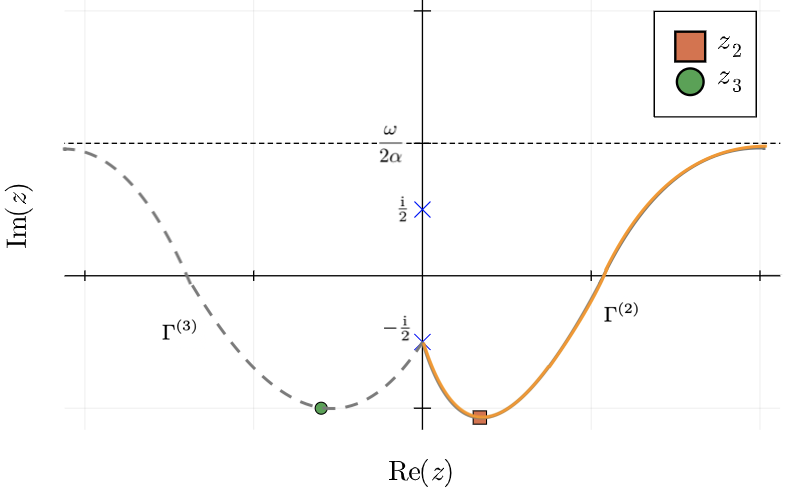}
    \end{center}
  \caption{\textit{MT functions:} Schematic representation of the contour of integration $\Gamma$.  The steepest descent path looks like a Gaussian centred around our root.}\label{fig:contourconstruct}
\end{figure}

The same idea can be applied to the remaining integral, $\Gamma^{(3)}$, from $-\infty$ to  $0$. For a schematic representation see Fig.~\ref{fig:contourconstruct}. 

The key argument used in constructing $\Gamma^{(2)}$ is that the imaginary part of the path of steepest descent converges to $\omega / 2 \alpha$ as $\textnormal{Re}\,z \rightarrow \infty.$ Note, this also eliminates the choice of $z_{1}$ as $ \textnormal{Im}(z_{1}) > \omega/2 \alpha$ which can be seen from eqn~\eqref{eq:z1est}. The same idea holds for $\textnormal{Re}\,z \rightarrow -\infty$. We can observe this in Fig.~\ref{fig:contourplot} for a given $\alpha$, $\omega$ $x_{0}$ and $n$.

\begin{lemma}
    For sufficiently large $x$ there exists a function $y = \gamma(x)$ such that $x + \ii y \in \Gamma$. Furthermore, $\lim_{x \rightarrow \pm \infty} \gamma(x) = \omega/2\alpha.$
\end{lemma}
\begin{proof}
	
  Our steepest descent contour satisfies
  \begin{Eqnarray*}
	-x +  2\beta y (x - x_{0})  + \cw\tan^{-1}\! \left(\frac{4x}{1 - 4x^{2} - 4y^{2}}\right) = \textnormal{Im}\,g(z_{*}),
  \end{Eqnarray*}
  where for some fixed $n$ and $\omega$ the RHS is constant.
  Let $z_{*} = c + \ii d$, where $c$ and $d$ are constants. We have
  \begin{Eqnarray*}
	0 &=& x\left( 2\beta y -1 \right) - c\left(2 \beta d -1 \right)+2 \beta x_{0} (d-y)\\ 
	& &  \mbox{}+\cw \left[\tan^{-1}\!\left(\frac{4 x}{1 - 4x^{2} - 4y^{2}} \right)- \tan^{-1}\!\left(\frac{4 c}{1 - 4 c^{2} - 4d^{2}}\right)\right] \!.
  \end{Eqnarray*}
  We retrieve the $\gamma(x)$ by the implicit function theorem. Dividing through by $x$ and letting $x \rightarrow \pm \infty$,
  \begin{Eqnarray*}
	0 &=& 2\beta y -1+ \lim_{x \rightarrow \pm \infty}\left[\frac{1}{x} \left( -c \left(2 \beta d -1 \right) +2 \beta x_{0} (d-y) \right) \right] \nonumber \\
	& & \mbox{}+ \left(\lim_{x \rightarrow \pm \infty} \frac{1}{x}\right)\cw\left[\lim_{x \rightarrow \pm \infty} \tan^{-1}\!\left(\frac{4}{1/x - 4x - 4y^{2}/x}\right)- \tan^{-1}\!\left(\frac{4 c}{1 - 4 c^{2} - 4d^{2}}\right)\right]\!. \nonumber 
   \end{Eqnarray*}
 Taking the limit we are left with $y = 1/(2\beta)= \omega/(2 \alpha)$, (where $\beta = \alpha/\omega$) as required.
\end{proof}

Once the desired path is chosen, we can apply the MSD as discussed in section \ref{sec:msdexplain}. The integral in \eqref{eq:breint} can be rewritten as 
\begin{Eqnarray*}
b_n(\omega)=\frac{(-\ii)^n}{\sqrt{2\pi}}\left(\int_{\Gamma^{(3)}}\ee^{-\omega g(z)}\,\frac{\dd z}{1+2\ii z}+
\int_{\Gamma^{(2)}}\ee^{-\omega g(z)}\,\frac{\dd z}{1+2\ii z}\right),
\end{Eqnarray*}
where the path of integration $\Gamma^{(3)}$ starts at $-\infty$, passes through the saddle point $z_3$ and terminates at $-\frac12\ii$, and $\Gamma^{(2)}$ starts at that point, runs through the saddle point $z_2$ and terminates at $+\infty$.

The transformations
\begin{equation}\label{eq:note10}
 g(z)-g(z_2)=\tfrac12 t^2,\qquad g(z)-g(z_3)=\tfrac12 t^2
\end{equation}
result in
\begin{equation*}
b_n(\omega)=
\frac{(-\ii)^n}{\sqrt{2\pi}}\left(\ee^{-\omega g(z_3)}\!\int_\R e^{-\frac12\omega t^2}f^{(3)}(t)\,\dd t+
e^{-\omega g(z_2)}\!\int_\R \ee^{-\frac12\omega t^2}f^{(2)}(t)\,\dd t\right)\!,
\end{equation*}
where\footnote{$f^{(*)}(t)$ denotes $f^{(2)}(t)$ or $f^{(3)}(t)$, $z^{*}$ denotes $z_2$ or $z_3$. Similar notation will be used throughout. }
\begin{equation}\label{eq:note12}
f^{(*)}(t)=\frac{\dd z}{\dd t}\frac{1}{1+2\ii z},\quad \frac{\dd z}{\dd t}=\frac{t}{g^\prime(z)}.
\end{equation}

In the transformations in \eqref{eq:note10} we assume that  $z$ on the infinite part of the steepest descent path from $-\infty$ to $z_3$ corresponds to a nonpositive  $t$-value, and that $z$ on the finite part from $z_3$ to $-\frac12\ii$ with a nonnegative $t$-value. Similarly for the integral through $z_2$: the points $z$ on the part from $-\frac12\ii$ to $z_2$ correspond to $t\le 0$, the points $z$ on the part from $z_2$ to $+\infty$ correspond to $t\ge0$. In this way there is  a one-to-one relation of points $z\in \Gamma^{(*)}$ and $t\in\R$.
Observe that on the paths of steepest decent the quantities $g(z)-g(z^{*})$ are real and non-negative and their graphs are convex, like the graph of $\frac12t^2$ on $\R$.

Substituting the Taylor expansion $\displaystyle{f^{(*)}(t)=\sum_{k=0}^\infty a^{(*)}_k t^k}$, the asymptotic expansion becomes
\begin{equation} \label{eq:asympticbn}
b_{n}(\omega) \sim
 \frac{(-\ii)^n}{\sqrt{2 \pi}}
 \left( 
\ee^{-\omega g(z_3)}\sum_{k=0}^\infty a^{(3)}_{2k}\!\int_\R \ee^{-\frac12\omega t^2}  t^{2k}\,\dd t+
\ee^{-\omega g(z_2)}\sum_{k=0}^\infty a^{(2)}_{2k}\!\int_\R \ee^{-\frac12\omega t^2}  t^{2k}\,\dd t
 \right)\!,
\end{equation}
where the odd $k$ terms disappear as we are integrating an odd integrand over the real line.  The integrals now look like Gamma integrals of the form
\begin{Eqnarray*}
	\Gamma(z) = \int_{0}^{\infty} x^{z-1} e^{-x}  \, \dd x, \qquad \R(z) >0.
\end{Eqnarray*}
Let $y = \omega t^{2}/2$,  $\left(y^{-1/2}/\sqrt{2 \omega} \right)\dd y = \dd t $ and $t^{2k} = \left(2y / \omega\right)^{k}$, then
\begin{Eqnarray}
	\int_\R \ee^{-\frac12\omega t^2}  t^{2k}\,\dd t &=& \sqrt{\frac{2}{\omega} }\frac{2^{k}}{\omega^{k}} \int_{0}^{\infty} \ee^{-y} y^{k-1/2} \dd y,  \nonumber \\
		&=& \sqrt{ \frac{2}{\omega}} \frac{2^{k}}{\omega^{k}} \Gamma(k + 1/2),  \nonumber \\
		&=&\sqrt{	\frac{2 \pi}{\omega} }\frac{2^{k}}{\omega^{k}} \left(\frac{1}{2}\right)_{k}, \label{eq:gammaintform}
\end{Eqnarray}
where $\left( \frac{1}{2} \right)_{k}$ is the Pochhammer symbol, $\left(x\right)_{n} = \Gamma(x+n)/ \Gamma(x)$.
By substituting \eqref{eq:gammaintform} into \eqref{eq:asympticbn} we obtain
\begin{equation}\label{eq:note13}
b_n(\omega)\sim
\frac{(-i)^n}{\sqrt{\omega}}
\left(
e^{-\omega g(z_3)}\sum_{k=0}^\infty a^{(3)}_{2k}\frac{2^k\left(\frac12\right)_k}{\omega^k}+
e^{-\omega g(z_2)}\sum_{k=0}^\infty a^{(2)}_{2k}\frac{2^k\left(\frac12\right)_k}{\omega^k}
\right)\!,
\end{equation}
as $\omega\to\infty$, and large $n=\cw\omega$ such that  $\cw\ge \frac14+\delta$, with $\delta$ a small positive number.
The next step is to determine the $a_{2k}$'s, the full derivation can be found in section \ref{sec:a2ks}.
The first coefficients have the form 
\begin{Eqnarray*}
	a_{0}^{(3)} = \frac{1}{(1 + 2 \ii z_{3})\sqrt{g''(z_{3})}}, \qquad a_{0}^{(2)} = \frac{1}{(1 + 2 \ii z_{2})\sqrt{g''(z_{2})}}.
\end{Eqnarray*}
Substituting into \eqref{eq:note13} we obtain the first-order approximation
\begin{Eqnarray}\label{eq:mt_bn_firstapprox}
b_n(\omega)\sim
(-\ii)^n
\left\{
 \frac{\ee^{-\omega g(z_3)}}{(1 + 2 \ii z_{3})\sqrt{\omega g''(z_{3})}}\left[1 + \mathcal{O}\left(\omega^{-1}\right) \right]  
+  \frac{\ee^{-\omega g(z_2)}}{(1 + 2 \ii z_{2})\sqrt{\omega g''(z_{2})}}\left[1 + \mathcal{O}\left(\omega^{-1}\right) \right]  
\right\}\!.
\end{Eqnarray}

\subsubsection{Asymptotic expansion} \label{sec:mtasyexpan}
To complete the asymptotic expansion,  we substitute \eqref{eq:z2est}, \eqref{eq:z3est} and $\beta = \alpha /\omega$ into the integrand in \eqref{eq:mt_bn_firstapprox} and expand for large $\omega$. This yields
\begin{Eqnarray}
	-\omega g(z_{2}) &\sim&- \frac{\alpha}{4 }\left(\sqrt{4\cw -1} - 2x_{0} \right)^{\!2} + \frac{\ii \omega}{2}\sqrt{4\cw -1} - 2 \ii \cw \omega \tan^{-1}\!\left(\sqrt{4\cw -1}\right)  + {\cal O}\left(\omega^{-1}\right)\!, \nonumber \\
	\omega g''(z_{2}) &\sim & 8 \alpha - \frac{2\alpha}{\cw \sqrt{4\cw -1}}\left( 6 x_{0}\cw + \sqrt{4\cw -1} - 2 x_{0}\right)  - \frac{\ii \omega}{\cw} \sqrt{4\cw -1} + {\cal O}\left(\omega^{-1}\right)\!, \nonumber \\
	1 + 2 \ii z_{2} &\sim&1 + \ii \sqrt{4 \cw -1} + {\cal O}\left(\omega^{-1}\right)\!,\nonumber\\
	-\omega g(z_{3}) &\sim& - \frac{\alpha}{4}\left(\sqrt{4\cw -1} + 2 x_{0}\right)^{\!2} \label{eq:gz3approx} \\
	& & \mbox{} - \frac{\ii\omega}{2}\sqrt{4\cw -1} + 2 \ii \cw \omega \tan^{-1}\!\left(\sqrt{4\cw -1}\right)  + {\cal O}\left(\omega^{-1}\right)\!,\nonumber \\
	\omega g''(z_{3}) &\sim & 8 \alpha - \frac{2\alpha}{\cw \sqrt{4\cw -1}}\left( -6 x_{0}\cw + \sqrt{4\cw -1} + 2 x_{0}\right)  + \frac{\ii \omega}{\cw} \sqrt{4\cw -1} + {\cal O}\left(\omega^{-1}\right)\!,\nonumber \\
	1 + 2 \ii z_{3} &\sim& 1 - \ii \sqrt{4 \cw -1} + {\cal O}\left(\omega^{-1}\right)\!.\nonumber
\end{Eqnarray}
It follows from the first term in \eqref{eq:gz3approx} that, as $x_{0}$ increases, $\ee^{-\omega g(z_{3})}$ gets smaller than $\ee^{-\omega g(z_{2})}$. 
Substituting the terms inside \eqref{eq:mt_bn_firstapprox},  the following asymptotic expansion holds for $ \cw = n/\omega> \frac{1}{4}$, $n \geq 0$, $\alpha,x_{0},\omega \geq 0$,
\begin{Eqnarray*}
b_n(\omega) &\sim& (-\ii)^n\!
\left\{\exp \!\left( - \frac{\alpha}{4 }\left(\sqrt{4\cw\! -\!1} -2 x_{0} \!\right)^{\!2} \!+\! \frac{\ii \omega}{2}\sqrt{4\cw\! -\!1} - 2 \ii \cw \omega \tan^{-1}\!\left(\sqrt{4\cw \!-\!1}\right)  + {\cal O}\left(\omega^{-1}\right) \!\right)\right. \\
&& \mbox{}\times \left[8 \alpha - \frac{2\alpha}{\cw \sqrt{4\cw -1}}\left( 6 x_{0}\cw + \sqrt{4\cw -1} - 2 x_{0}\right)  - \frac{\ii \omega}{\cw} \sqrt{4\cw -1} + {\cal O}\left(\omega^{-1}\right)  \right]^{-1/2} \\
& &  \mbox{} \times \left[1 + \ii \sqrt{4 \cw -1} + {\cal O}\left(\omega^{-1}\right) \right]^{-1} \left[1 + \mathcal{O}\left(\omega^{-1}\right) \right]  \\
&& \mbox{}+ \exp\! \left( - \frac{\alpha}{4}\left(\sqrt{4\cw -1} + 2 x_{0}\right)^{2} - \frac{\ii\omega}{2}\sqrt{4\cw -1} + 2 \ii \cw \omega \tan^{-1}\!\left(\sqrt{4\cw -1}\right)  + {\cal O}\left(\omega^{-1}\right)\!\right)  \\
& & \mbox{} \times \left[ 8 \alpha - \frac{2\alpha}{\cw \sqrt{4\cw -1}}\left( -6 x_{0}\cw + \sqrt{4\cw -1} + 2 x_{0}\right)  + \frac{\ii \omega}{\cw} \sqrt{4\cw -1} + {\cal O}\left(\omega^{-1}\right)\right]^{-1/2} \\
& &  \left.\mbox{} \times \left[1 - \ii \sqrt{4 \cw -1} + {\cal O}\left(\omega^{-1}\right) \right]^{-1} \left[1 + \mathcal{O}\left(\omega^{-1}\right) \right]  
\right\}, \\ 
&\sim & {\cal C}\ee^{-\alpha\left(\sqrt{\cw -1/4} - x_{0} \right)^{2} } \left[1 + \mathcal{O}\left(\omega^{-1}\right) \right]  
\end{Eqnarray*}

The final line of the above equation vividly demonstrates that our asymptotics exhibit exponential decay. This can also be seen in appendix \ref{sec:boundbnw} where we find a bound for $|b_{n}(\omega)|$ which leads to the result in \eqref{eq:mtestfull}. All that remains is to consider the $n<0$ and $\omega <0$ cases.

\subsection{Reduction to the case $n \geq 0 $ and $\omega > 0 $} \label{sec:reduction}

At the outset we have assumed $n,\omega >0$, in this subsection we demonstrate that this  can be extended to $n\in \R$, $\omega \in \R$. 

\subsubsection{A symmetry argument}

Noting that $\varphi_{-n-1} (x)= (-\ii)^{2 n + 1} \varphi_{n}(- x)$ for all $n \in \Z$, consider
\begin{displaymath}
  b_{-n-1}(-\omega) = \frac{1}{2} \int_{-\infty}^{\infty} \ee^{-\alpha(z- x_{0})^{2} - \ii \omega z } \overline{\varphi_{-(n+1)}(z)}\dd z \nonumber =\frac{\ii ^{2n+1}}{2} \int_{-\infty}^{\infty} \ee^{-\alpha(z- x_{0})^{2} - \ii \omega z } \overline{\varphi_{n}(-z)}\dd z.
\end{displaymath}
Let $z' = -z$, $\dd z' = - \dd z$. Relabelling, we obtain
\begin{equation}
   b_{-n-1}(-\omega)= -\ii^{2n+1} \frac{1}{2} \int^{\infty}_{-\infty} \ee^{-\alpha(z- x_{0})^{2} + \ii \omega z} \overline{\varphi_{n}(z)} \dd z = -\ii^{2n+1} b_{n}(\omega) = (-1)^{n-1} \ii b_{n}(\omega). \label{eq:mtsymm}
\end{equation}
Using this symmetry we  deduce the results for $n \leq -1, \omega \leq 0$ from $n \geq 0, \omega \geq 0$ for all $ \ n \in \Z$.

\subsubsection{Laguerre coefficients argument}

We can further reduce our range by considering the case where $n \leq -1, \omega > 0$ by using the relationship of MT functions to Laguerre polynomials, considered in \cite{iserles2019family}. Taking an inverse Fourier transform of our wave packet
\begin{displaymath}
	 \ee^{-\alpha ( z- x_{0})^{2} + \ii \omega z} =\frac{1}{\sqrt{2 \alpha}} \frac{1}{\sqrt{2\pi}} \int_{-\infty}^{\infty}  \ee^{- (\omega + \xi)^{2}/(4\alpha) + \ii x_{0}(\omega + \xi)} \ee^{-\ii \xi z} \dd z = \frac{1}{\sqrt{2\alpha}} \mathcal{F}^{-1}\!\left[ \ee^{- (\omega + \xi)^{2}/(4\alpha) + \ii x_{0}(\omega + \xi)}  \right]\!(z).
\end{displaymath}
It is proved in \cite{iserles2019family} that for $n\leq -1$ the MT functions are
\begin{displaymath}
  \varphi_n(x) = -\frac{(-\ii)^{n}}{\sqrt{2\pi}}\int_{0}^\infty \mathrm{L}_{-n-1}(\xi) \ee^{-\xi/2-\ii x \xi} \, \dd \xi = -(-\ii)^{n} \mathcal{F}^{-1}\!\left[ \chi_{[0,\infty)}(\xi) \mathrm{L}_{-(n+1)}(\xi) \ee^{-\xi/2} \right]\!(x)
\end{displaymath}
Therefore, substituting into the expression for the coefficients in terms of $\ee^{\ii\omega x -\alpha( z-x_{0})^2}$ and $\varphi_n$, we get
\begin{displaymath}
  b_n(\omega) = \frac{-\ii^{n}}{2^{\frac32}\sqrt{\alpha}}\int_{-\infty}^\infty \mathcal{F}^{-1}\!\left[ \ee^{-( \xi + \omega)^{2}/(4\alpha) + \ii x_{0}( \xi+\omega)} \right]\!(z)  \times \overline{\mathcal{F}^{-1}\!\left[ \chi_{[0,\infty)}(\xi) \mathrm{L}_{-(n+1)}(\xi) \ee^{-\xi/2} \right]\!(z)} \, \dd z
\end{displaymath}
By Parseval's identity, we deduce that
\begin{displaymath}
  b_n(\omega) = \frac{-\ii^{n}}{2^{\frac32}\sqrt{\alpha}}\int_{0}^\infty \ee^{-(\xi+\omega)^2/(4\alpha) + \ii x_{0}(\xi + \omega)} \, \mathrm{L}_{-n-1)}(\xi) \ee^{-\xi/2} \, \dd \xi.
\end{displaymath}

As both $\mathrm{L}_{n}(\xi) \ee^{-\xi/2}$ and $\ee^{\ii x_{0}(\xi - \omega)} $ have $\mathrm{L}_2(\mathbb{R})$ unit norm, by Cauchy--Schwartz,
\begin{Eqnarray*}
  |b_n(\omega)| &\leq& \frac{2^{-3/2}}{\sqrt{\alpha}} \left(\int_{0}^\infty \ee^{-(\xi+\omega)^2/(2\alpha)} \left|\ee^{\ii x_{0}(\xi + \omega)} \right|^{2} \, \dd \xi\right)^{\!1/2}\!\left(\int_{0}^\infty \left|\mathrm{L}_{-n-1}(\xi)\right|^2 \ee^{-\xi} \, \dd \xi\right)^{\!1/2} \\
  &=& \frac{2^{-3/2}}{\sqrt{\alpha}} \left(\int_{0}^\infty \ee^{-(\xi+\omega)^2/(2\alpha	)} \, \dd \xi\right)^{\!1/2} = \frac{2^{-3/2}}{\sqrt{\alpha}} \ee^{-\omega^2/(4\alpha)} \!\left(\int_{0}^\infty \ee^{-\xi^2/(2\alpha)} \ee^{-(\xi\omega)/\alpha} \, \dd \xi\right)^{\!1/2}.
\end{Eqnarray*}
Now, if $\omega \geq 0$ then $\ee^{-\xi\omega} < 1$. Finally, $\int_{0}^{\infty} \ee^{-\xi^{2}/(2\alpha)}\dd \xi = \sqrt{\pi \alpha/2}$ leads to
\begin{equation*}
  |b_n(\omega)| \leq \pi^{\frac14}2^{-\frac74} \alpha^{-1/4} \ee^{-\omega^2/(4\alpha)} \text{  for  } \omega \geq 0,\;\; n\leq -1.
\end{equation*}
Thus by the symmetry argument in  \eqref{eq:mtsymm}, this result holds for $\omega \leq 0, n \geq 0$. 

As we have shown in our crude bound above, the coefficients for $b_{n}(\omega)$ are small for $\omega \geq 0$ and $n \leq -1$. Moreover, we can use the above symmetry argument to evaluate $b_{n}(-\omega)$.

Finally, tying everything together, we formulate our estimate for $a_{n}(\omega)$ in equation \eqref{eq:anmtfunc} throughout the entire range of parameters.
 
\begin{figure}[bt]
  \begin{center}
  \includegraphics[width=.42\textwidth]{mta1w100x00.png}\hspace*{20pt}
  \includegraphics[width=.42\textwidth]{mta2w200x02.png}
  \includegraphics[width=.42\textwidth]{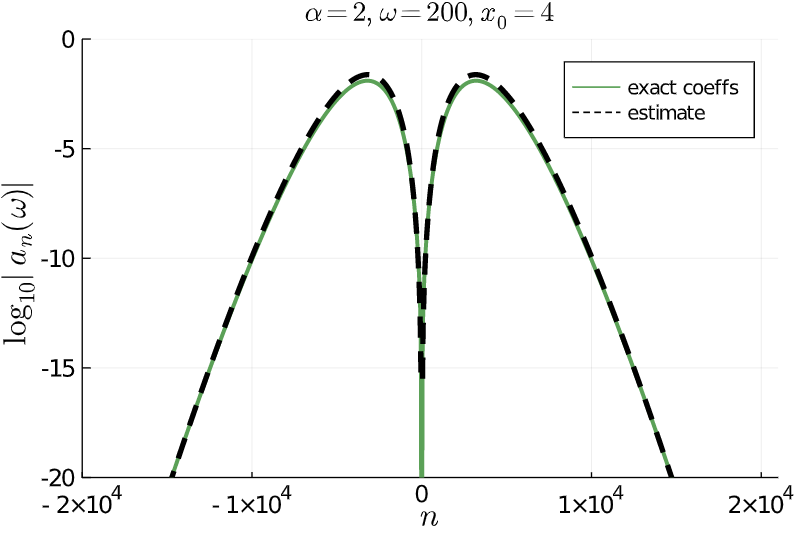}\hspace*{20pt}
  \includegraphics[width=.42\textwidth]{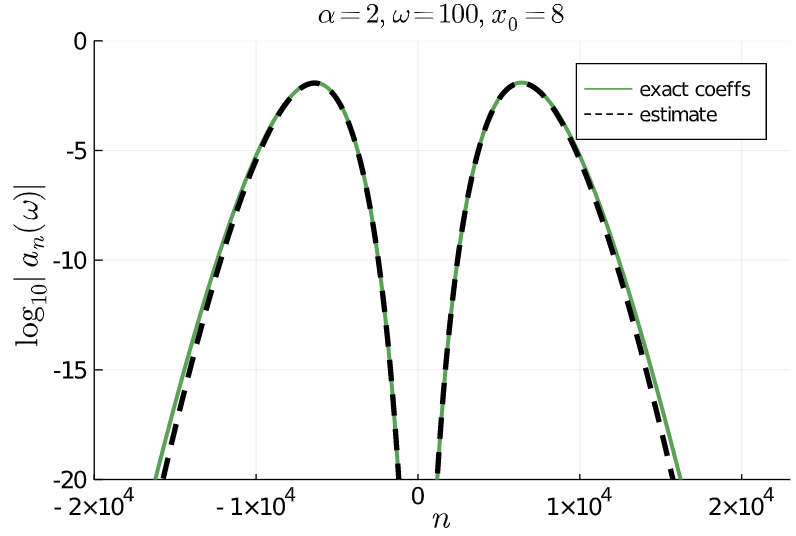}
    \end{center}
  \caption{\textit{MT functions:} Plot of  \eqref{eq:mtestfull} for varying $\alpha, \omega $ and $x_{0}$. $x_{0} = 0$, demonstrates linear spectral decay. As we increase $x_{0}$, the width grows. }\label{fig:4.6}
\end{figure}

\begin{figure}[bth]
  \begin{center}
  \includegraphics[width=.42\textwidth]{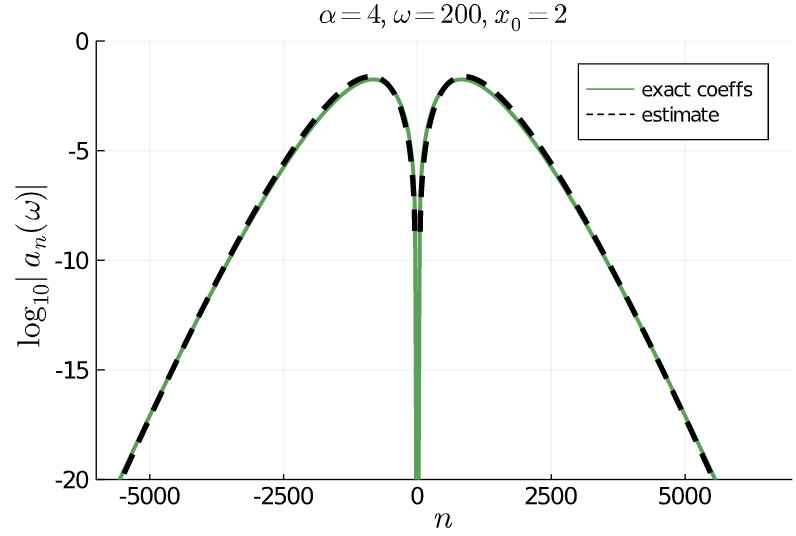}\hspace*{20pt}
  \includegraphics[width=.42\textwidth]{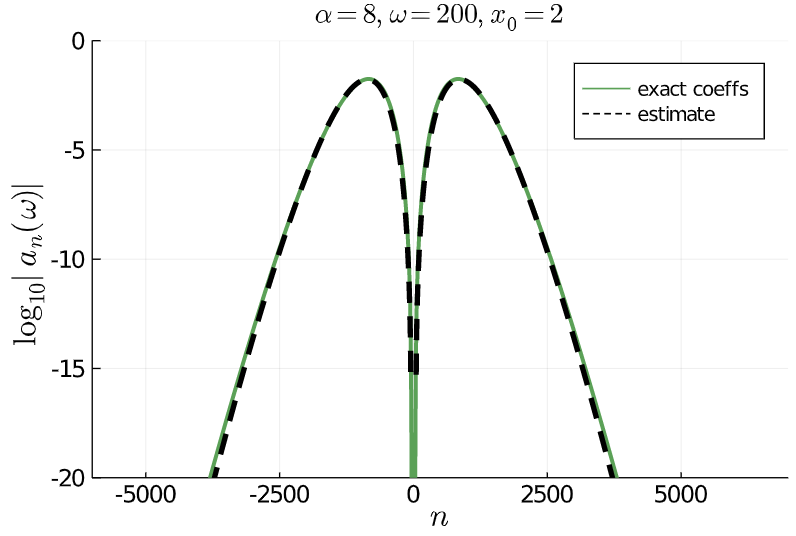}
  \includegraphics[width=.42\textwidth]{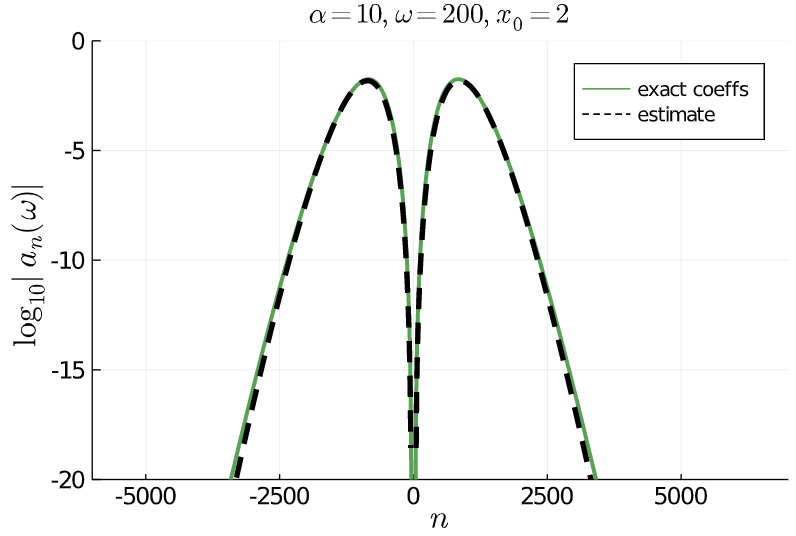}\hspace*{20pt}
  \includegraphics[width=.42\textwidth]{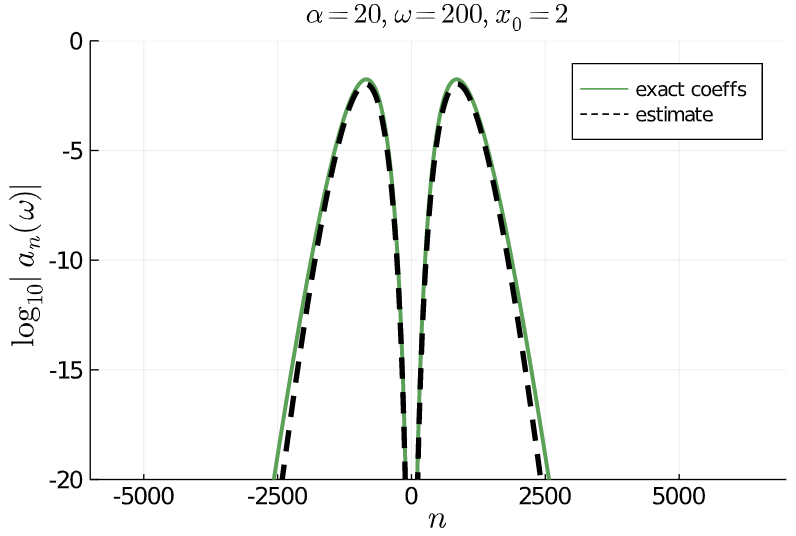}
    \end{center}
  \caption{\textit{MT functions:} Plot of \eqref{eq:mtestfull} for varying $\alpha$, and fixed $x_{0}$ and $\omega$. As we increase $\alpha$ we see that the spectral decay rate increases.  }\label{fig:4.7}
\end{figure}

\begin{theorem}[Malmquist--Takenaka function]\label{thm:mtfunc}
  Suppose that $f(z) = \ee^{-\alpha( z-x_{0})^{2}}\cos(\omega z)$, where $\alpha >0$, while $x_{0}$ and $\omega$ are real. The asymptotic behaviour of the coefficients, $a_{n}$, with Malmquist--Takenaka basis functions 
  \begin{Eqnarray*}
 	\varphi_{n}(z) = \ii^{n}\sqrt{\frac{2}{\pi}}\frac{(1 + 2\ii z)^{n}}{(1 - 2\ii z )^{n+1}}, \quad n \in \Z,
  \end{Eqnarray*}
  is subject to the estimate 
  \begin{Eqnarray*}
  	|a_{n}(\omega)| \leq |b_{n}(\omega)| + |b_{n}(-\omega)|,
  \end{Eqnarray*}
 where for $\omega \geq 0$ 
  \begin{Eqnarray}
 	n < 0: &&
  |b_n(-\omega)| \leq \pi^{\frac14}2^{-\frac74} \alpha^{-1/4}\cdot \ee^{-\omega^2/(4\alpha)}, \textnormal{ and } \nonumber \\
  n\geq 0: &&
  \label{eq:mtestfull}
		|b_{n}(\omega)| \leq \!\frac{\left[ 1 + \mathcal{O}\left(\omega^{-1}\right)\right]}{2\sqrt{2 \omega}\left( \left| n/\omega \right| -1/4 \right)^{1/4}}\!\left\{ 
			\!\exp\!\left( - \alpha\! \left(\sqrt{\left| \frac{n}{\omega}\right| -\frac{1}{4}} - x_{0} \right)^{2}  \right) \!+ \exp\!\left( - \alpha\! \left(\sqrt{\left| \frac{n}{\omega}\right| -\frac{1}{4}} + x_{0} \right)^{2}  \right)\! \right\}
			\!. \hspace*{25pt}
  \end{Eqnarray} 
  The bound in \eqref{eq:mtestfull} is valid uniformly for $n$ satisfying
$1/4 \leq \left| \frac{n}{\omega}\right| \leq c$ for any given constant $c$. The remaining case, $\omega <0$, can be found by applying the symmetry argument in \eqref{eq:mtsymm}.  These bounds hold as $|\omega|\rightarrow \infty$.
\end{theorem}
\begin{remark}
For practical purposes, the values of $n$ satisfying $1/4 \leq \left| \frac{n}{\omega}\right| \leq c$ are the only values of $n$ we care about because, with reasonable assumptions on the size of $c$, we can ensure that $|a_{n}| < \varepsilon$ for $n$ outside this range.
\end{remark}
Figs \ref{fig:4.6} and \ref{fig:4.7} depict the asymptotic bound in \eqref{eq:mtestfull}. We see that they capture exceedingly well the characteristics of the exact coefficients and can observe asymptotic spectral decay.

\section*{Acknowledgements}

KL wishes to acknowledge UK Engineering and Physical Sciences Research Council (EPSRC) for the grant EP/L016516/1 for the University of Cambridge Centre for Doctoral Training, the Cambridge Centre for Analysis, which supported the research and writing of this paper. KL is also grateful to the Department of Computer Science, KU Leuven for their hospitality during a visit where some of the research for this paper was done. MW thanks FWO Research Foundation Flanders for the postdoctoral research fellowship he enjoyed while some of the research for this paper was done. MW also thanks the Polish National Science Centre (SONATA-BIS-9), project no. 2019/34/E/ST1/00390, for the funding that supported the research. The authors would like to thank Bruno Salvy and Andr\'e Weideman for helpful correspondence on the steepest descent method, and are especially grateful to Nico Temme for the helpful discussions on Hermite functions and the method of steepest descent which helped simplify the section on the Malmquist--Takenaka functions.
\bibliographystyle{plain}
\bibliography{approx_wave_packets_real_line} 

\appendix
\section{Hermite}
\subsection{Different values of $\alpha$ } \label{sec:alphacases}
Starting from equation \eqref{eqn:Hermiteantilde}, it is helpful to define
\begin{equation}\label{eqn:atilden}
  \tilde{a}_{n}(\omega) = \ee^{-\alpha x_{0}^{2}} \mathcal{H}_{n,(1+2\alpha)^{-1},\omega-2\ii\alpha x_0},
  \end{equation}
which satisfies
\begin{equation}\label{eqn:Hermiteantilde1}
  a_{n}(\omega)=\frac{1}{\sqrt{2^{n}n! \pi^{1/2}}} \left(\frac{2\pi}{1+2\alpha}\right)^{\!1/2} \textnormal{Re}\,\tilde{a}_{n}(\omega),
  \end{equation}
and distinguish between the cases $\alpha = \tfrac12$ and $\alpha \neq \tfrac12$.

\vspace{6pt}
\noindent\textbf{Case I:} $\alpha=\frac12$. We need to use  \eqref{eq:8.20}, namely
\begin{Eqnarray*}
  \tilde{a}_n(\omega) &=& \exp\!\left(-\frac{\omega^2}{4}-\frac{x_0^2}{4} + \ii x_0\omega/2\right)\!(x_0 + \ii \omega)^n \\
                      &=& (x_0^2+\omega^2)^{n/2} \ee^{-(x_0^2+\omega^2)/4} \exp\left( \ii\frac{\omega x_{0}}{2} + \ii n \tan^{-1}\!\left(\frac{\omega}{x_0}\right)   \right).
  \end{Eqnarray*}
Substituting this into equation \eqref{eqn:Hermiteantilde1} yields
\begin{Eqnarray*}
  a_n(\omega) &=& \sqrt{\frac{\pi^{1/2}}{2^{n}n!}}(x_0^2 + \omega^2)^{n/2} \ee^{-(x_0^2+\omega^2)/4} \cos\left(\frac{\omega x_{0}}{2} + n \tan^{-1}\!\left(\frac{\omega}{x_0}\right) \right).
\end{Eqnarray*}

\vspace{6pt}
\noindent \textbf{Case II:} $\alpha>\frac12$. Using equation \eqref{eq:8.18} and \eqref{eqn:atilden}, we obtain,
\begin{Eqnarray*}
  \tilde{a}_{n}(\omega)&=& \left(\frac{2\alpha-1}{1+2\alpha}\right)^{\!n/2}\!\! \exp\!\left(-\frac{\alpha x_0^2}{1+2\alpha}-\frac{\omega^2}{2(1+2\alpha)} + \frac{2\ii\alpha \omega x_0}{1+2\alpha}\right) \HH_n\!\left(\frac{\omega-2\ii\alpha x_0}{(1-4\alpha^2)^{1/2}}\right)\\
   &=& \left|\frac{1-2\alpha}{1+2\alpha}\right|^{\!n/2}\!\! \exp\!\left(\frac{|1-2 \alpha
 |}{4\alpha}(X^2+2\alpha Y^2) - \ii|1-2 \alpha| XY \!\right)\HH_n(X+\ii Y)
\end{Eqnarray*}
where
\begin{Eqnarray*}
  X=\frac{2\alpha x_0}{|1-4\alpha^2|^{1/2}},\qquad Y=\frac{\omega}{|1- 4\alpha^2|^{1/2}}.
\end{Eqnarray*}
Note that $Y\gg X$ in the high frequency regime.

\vspace{6pt}
\noindent\textbf{Case III:} $\alpha<\frac12$. We obtain,
\begin{Eqnarray*}
  \tilde{a}_{n}(\omega)&=& \ii^n \left(\frac{1-2\alpha}{1+2\alpha}\right)^{\!n/2}\!\! \exp\!\left(-\frac{\alpha x_0^2}{1+2\alpha}-\frac{\omega^2}{2(1+2\alpha)} + \frac{2\ii\alpha \omega x_0}{1+2\alpha}\right) \HH_n\!\left(\frac{\omega-2\ii\alpha x_0}{(1-4\alpha^2)^{1/2}}\right)\\
  &=& \left|\frac{1-2\alpha}{1+2\alpha}\right|^{\!n/2}\!\! \exp\!\left(-\frac{|1-2 \alpha|}{4\alpha}(X^2+2\alpha Y^2) + \ii|1-2 \alpha| XY \!\right) \ii^n\HH_n(Y-\ii X).
\end{Eqnarray*}

\subsection{$\alpha > 1/2$: asymptotics} \label{sec:hermiteworking1}
The terms $w_{-}\phi''(w_{-})$ and $\phi(w_{-})$ in term of $X$ and $Y$ are
\begin{Eqnarray*}
  w_{-}\phi''(w_{-}) &=& - 2\ii \sqrt{1-\zeta^{2}} =-\frac{2\ii}{\nu}\sqrt{-X^{2}-2\ii XY + Y^{2} + \nu^{2}} \\
  \phi(w_{-}) &=& w_{-}^{2} - 2\zeta w_{-} + \frac{1}{2}\log(w_{-})\\
  &=&\mbox{}-\frac{1}{4} - \frac{Y - \ii X}{2\left[ Y - \ii X + \sqrt{-X^{2} - \ii 2 X Y + Y^{2} + \nu^{2}}\right]} \\
  &&\mbox{}-\frac{1}{2}\log (2\ii) - \frac{1}{2}\log\left| \frac{Y - \ii X+ \sqrt{-X^{2} - \ii 2 X Y + Y^{2} + \nu^{2}}}{\nu}\right|,
\end{Eqnarray*}
where
\begin{displaymath}
	w_{-} = \frac{1}{2}\left(\zeta - \ii \sqrt{1 - \zeta^{2}}\right), \qquad \nu = \sqrt{2 n +1}.
\end{displaymath}

Substituting these terms into  \eqref{eq:A1} and  \eqref{eq:B1},
\begin{subequations}
\begin{Eqnarray}
  A&=&\frac{1}{\ii \pi^{1/4}\sqrt{-2 \ii (2\alpha+1)}}\exp\!\left(-\frac{\alpha x_0^2}{1+2\alpha}+\frac{2\ii\alpha\omega x_0}{1+2\alpha}\right) \frac{1}{\sqrt{\nu}} \frac{1}{\left(-X^{2}-2\ii XY + Y^{2} + \nu^{2}\right)^{1/4}}, \label{eq:A2}\\
  B&=&\left(\frac{n!}{2^n}\right)^{\!1/2} \frac{1}{\nu^n} \exp\!\left(-\frac{\omega^{2}}{2(2\alpha +1)}+\frac{n}{2}\log \left|\frac{2\alpha-1}{2\alpha+1}\right| -\nu^2\phi(w_-)\right). \label{eq:B2}
\end{Eqnarray}
\end{subequations}

We need to be careful when choosing which large parameter we wish to expand around. Recalling that
\begin{Eqnarray*}
	\zeta = \frac{X + \ii Y}{\nu} = \frac{2\alpha x_{0} + \ii \omega}{\sqrt{(4\alpha^{2}-1)(2n + 1)}},
\end{Eqnarray*}
while $\alpha$ and $x_{0}$ are fixed variables, we want to preserve the relationship between our large parameters $n$ and $\omega$ so
\begin{displaymath}
	\sqrt{n} \propto \omega,\qquad n = c_{n,\omega} \omega^{2}.
\end{displaymath}
This agrees with our numerical experiment seen in figure \ref{fig:hermitea12asymp}.
 
We start by expanding  \eqref{eq:A2}, substituting $n =c_{n,\omega} \omega^{2} $:
\begin{Eqnarray*}
  \nu^{-1/2}=  \left(2 c_{n,\omega} \omega^{2} + 1\right)^{1/4} &=& \left( 2 c_{n,\omega} \omega^{2} \right)^{\!-1/4} \left(1 + \frac{1}{2 c_{n,\omega} \omega^{2} }\right)^{-1/4} \\
  &=&\frac{1}{(2 c_{n,\omega})^{1/4} \omega^{1/2}} - \frac{1}{8 \cdot 2^{1/4} c^{5/4} \omega^{5/2}}  \\
  & & \mbox{}+ \frac{5}{128 \cdot 2^{1/4} c^{9/4} \omega^{9/2}}+ \mathcal{O}(\omega^{-13/2}), \\
  \left(-X^{2}-\ii 2 X Y + Y^{2} +\nu^{2}\right)^{-1/4} & \sim & \frac{1}{\sqrt{\omega}} \left(\frac{4\alpha^{2}-1}{1+2c_{n,\omega}(4\alpha^{2}-1)}\right)^{\!1/4} \\
  & &  \mbox{}+ \frac{\ii x_{0}\alpha \left(4\alpha^{2}-1\right)^{1/4}}{\omega^{3/2}\left[1+2c_{n,\omega}\left(4\alpha^{2} -1\right)\right]^{5/4}} + \mathcal{O}\left(\omega^{-5/2}\right)\!, \\
  \nu^{-1/2} \left(-X^{2}-\ii 2 X Y + Y^{2} +\nu^{2}\right)^{-1/4} &=& \left[ \frac{4\alpha^{2}-1}{1+2c_{n,\omega}\left(4\alpha^{2}-1\right)}\right]^{\!1/4} \frac{1}{(2c_{n,\omega})^{1/4}\omega} \\
  & & \mbox{}\times \left[ 1 + \frac{\ii \alpha x_{0}}{\omega \left(1 + 2c_{n,\omega}\left(4\alpha	^{2}-1\right)\right)} + \mathcal{O}\left(\omega^{-2}\right) \right].
\end{Eqnarray*}

Substituting terms back into \eqref{eq:A2}, we obtain
\begin{Eqnarray*}
  A &\sim & \frac{1}{\ii \pi^{1/4} \sqrt{-2\ii (2\alpha + 1)}}\exp\left(\frac{{-\alpha x_{0}^{2} + 2\ii \omega\alpha x_{0}}}{2\alpha + 1}\right) \left[ \frac{4\alpha^{2}-1}{1+2c_{n,\omega}\left(4\alpha^{2}-1\right)}\right]^{1/4} \frac{1}{(2c_{n,\omega})^{1/4}\omega}  \\
  & &\mbox{} \times \left\{ 1 + \frac{\ii \alpha x_{0}}{\omega \left[1 + 2c_{n,\omega}\left(4\alpha	^{2}-1\right)\right]} + \mathcal{O}\left(\omega^{-2}\right) \right\}.
\end{Eqnarray*}

Evaluating $B$ in \eqref{eq:B2}, by the Sterling formula, 
\begin{displaymath}  
  \sqrt{\frac{n!}{2^{n}}}\frac{1}{\nu^{n}} \sim \frac{(2 \pi)^{1/4} n^{n/2 + 1/4} \ee^{n/2}}{2^{n/2}} \lim_{n \rightarrow \infty}(2n + 1)^{-n/2} \sim  (2 \pi n )^{1/4} \ee^{-\frac{1}{2}\left(n + \frac{1}{2}\right)}2^{-n}. 
\end{displaymath}
Thus, $B$ is estimated by
\begin{subequations}
\begin{Eqnarray}
  B&\sim& (2\pi n)^{1/4} 2^{1/2} \exp \!\left(\frac{\ii \pi}{2}\left(n + \frac{1}{2}\right) - \frac{\omega^{2}}{2(2\alpha+1)} + \frac{n}{2}\log \left| \frac{2\alpha + 1}{2\alpha - 1} \right| \right.\nonumber \\
  & & \left.\mbox{}+ \frac{(2 n +1)(Y-\ii X)}{2\left(Y - \ii X + \sqrt{-X^{2} - \ii 2 X Y + Y^{2} + \nu^{2}} \right)} \right. \label{eq:Brationalterm} \\
  & & \left. \mbox{}+ \frac{(2n +1)}{2}\log \!\left( \frac{Y - \ii X + \sqrt{-X^{2} - \ii 2 X Y + Y^{2} + \nu^{2}}}{\nu}\right)\!\right)\!. \label{eq:Blogterm}
\end{Eqnarray}
\end{subequations}

Setting $n = c_{n,\omega} \omega^{2}$ and expanding about $\omega$,
\begin{Eqnarray*}
  \text{Eqn }\eqref{eq:Brationalterm} &\sim& \frac{(2n + 1)}{2 \left[ 1 + \sqrt{1 + 2c_{n,\omega}\left(4\alpha^{2} -1 \right)}\right]} \\
  & & \mbox{} - \frac{2\ii c_{n,\omega} x_{0} \alpha (2n +1) \left(4\alpha^{2} -1\right)}{\omega\!\left[1+\sqrt{1+2c_{n,\omega}\left(4\alpha^{2}-1\right)}\right]\!\left[1+2c_{n,\omega}\left(4\alpha^{2}-1\right)+\sqrt{1+2c_{n,\omega}\left(4\alpha^{2}-1\right)}\right]} +\mathcal{O}\left(\omega^{-2}\right), \\
  \text{Eqn }\eqref{eq:Blogterm} &\sim& \frac{-1 +2 c_{n,\omega} + 8 c_{n,\omega} (x_{0}^{2}-1)\alpha^{2}}{4 \left[1 + 2c_{n,\omega}\left( 4\alpha^{2}-1\right)\right]^{3/2}} + \frac{1}{2} \log \left| \frac{1+\sqrt{1+2c_{n,\omega}\left(4\alpha^{2}-1\right)}}{\sqrt{2 c_{n,\omega}\left( 4\alpha^{2}-1\right)}}\right| \\
  & & \mbox{}+ \cw \omega^{2} \log \left| \frac{1 + \sqrt{1 + 2 \cw\left(4 \alpha^{2}-1\right)}}{\sqrt{2 \cw \left(4\alpha^{2} -1\right)}}\right| - \frac{\ii2c_{n,\omega}\alpha x_{0}\omega}{\sqrt{1 +2c_{n,\omega}\left(4\alpha^{2}-1\right)}} \\
  & & \mbox{}- \frac{\ii x_{0}\alpha}{3 \omega \left[1 + 2c_{n,\omega}\left(4\alpha^{2}-1\right)\right]^{5/2}}  \left\{ 3 + c_{n,\omega} \left( -9 + 6 c_{n,\omega} \right. \right.  \\
  & & \left.\left.\quad \mbox{}-4\alpha^{2}\left[-9+12c_{n,\omega} + 2x_{0}^{2}(1 + c_{n,\omega})\right] + 32  \alpha^{4} c_{n,\omega}(3 + x_{0}^{2})\right)\right\} + \mathcal{O}\left(\omega^{-2}\right). 
\end{Eqnarray*}
Assembling all this,
\begin{Eqnarray*}
  B&\sim&(2\pi n)^{1/4} 2^{1/2} \ee^{\ii \pi\left(n +1/2 \right)/2} \\
  & & \mbox{}\times \exp\! \left( \omega^{2}\! \left[-\frac{1}{2(2\alpha + 1)}+ c_{n,\omega} \!\left(\frac{1}{1+\sqrt{1+2c_{n,\omega}\left(4\alpha^{2}-1\right)}} + \log \left|\frac{1\!+\!\sqrt{1\!+\!2c_{n,\omega}\left(4\alpha^{2}\!-\!1\right)}}{\sqrt{2 c_{n,\omega}}\left(2\alpha\! + \!1\right)} \right| \right)\right] \right. \\
  & & \left. \quad \mbox{}- \frac{2\ii\alpha x_{0} \omega \left( -1+\sqrt{1+2c_{n,\omega}\left(4\alpha^{2}-1\right)}\right)}{4\alpha^{2}-1} + \frac{2 x_{0}^{2} \alpha^{2}}{ 4\alpha^{2} -1} \left(1 - \frac{1}{\sqrt{1 + 2c_{n,\omega} \left(4\alpha^{2} -1\right)}}\right) \right. \\
  & & \quad \left.\mbox{}+ \frac{1}{2} \log \left| \frac{1 + \sqrt{1 + 2 \cw \left(4\alpha^{2} -1\right)}}{\sqrt{2\cw \left(4\alpha^{2} -1\right)}} \right|+ \frac{\ii x_{0}\alpha \left[-3 + 6c_{n,\omega} + 8 c_{n,\omega}\left(-3 + x_{0}^{2}\right)\alpha^{2}\right]}{3\omega\left[1+2c_{n,\omega}\left(4\alpha^{2}-1\right)\right]^{3/2}} + \mathcal{O}\left(\omega^{-2}\right)\right)\!.
\end{Eqnarray*}
This results in
\begin{Eqnarray*}
  \tilde{a}_{n}(\omega) &\sim & \frac{1}{2 \sqrt{\omega(2\alpha +1)}} \left[ \frac{4\alpha^{2} -1}{1 + 2c_{n,\omega}(4\alpha^{2}-1)} \right]^{1/4} \left[ 1 + \frac{\ii \alpha	x_{0}}{\omega \left[ 1 + 2c_{n,\omega	}\left(4\alpha	^{2} -1\right) \right]} + \mathcal{O}(\omega^{-2}) \right] \\
  & & \mbox{}\times \exp \!\left( \omega^{2} \left[ - \frac{1}{2(2\alpha + 1)} \right. \right. \\
  & &  \left. \left. \quad \quad \quad\mbox{}+ c_{n,\omega}\left(\frac{\ii \pi}{2} +\frac{1}{1+\sqrt{1 +2c_{n,\omega}\left(4\alpha^{2}-1\right)}} + \log\left| \frac{1 + \sqrt{1 +2c_{n,\omega}\left(4\alpha^{2}-1\right)}}{\sqrt{2c_{n,\omega}}(2\alpha +1)}\right| \right)\right] \right. \\
  & & \quad \left. \mbox{}- \frac{2\ii \alpha x_{0} \omega \left(-2\alpha + \sqrt{1 + 2c_{n,\omega}\left(4\alpha^{2}-1\right)}\right)}{4\alpha^{2}-1} \right. \\
  & & \left. \quad \mbox{}+ \frac{1}{2}\log \left| \frac{1 +\sqrt{1+2c_{n,\omega}\left(4\alpha^{2}-1\right)}}{\sqrt{2\cw \left(4\alpha^{2}-1 \right)}} \right|+ \frac{\alpha x_{0}^{2}}{4\alpha^{2}-1}\left[1 -  \frac{2 \alpha}{\sqrt{1+2c_{n,\omega}\left(4\alpha^{2}-1\right)}}\right]   \right. \\
  & & \left. \quad \mbox{}+ \frac{\ii x_{0} \alpha \left[-3 + 6c_{n,\omega} + 8 c_{n,\omega}\alpha^{2}\left(x_{0}^{2} - 3\right)\right]}{3\omega\left[1 + 2 c_{n,\omega}\left(4\alpha^{2}-1\right)\right]^{3/2}} + \mathcal{O}\left(\omega^{-2}\right)\!\right)\!.
\end{Eqnarray*}
Restoring  $c_{n,\omega} = n/\omega^{2}$ and recalling that $a_{n}(\omega) \sim \textnormal{Re}\,\tilde{a}_{n}(\omega),$ we obtain
\begin{Eqnarray*}
  a_{n}(\omega) &\sim & \frac{1}{\sqrt{\omega(2\alpha +1)}} \left[ \frac{4\alpha^{2} -1}{1 + 2\frac{n}{\omega^{2}}(4\alpha^{2}-1)} \right]^{1/4} \\
  & &\mbox{} \times \left| \frac{1 + \sqrt{1 +2\frac{n}{\omega^{2}}\left(4\alpha^{2}-1\right)}}{\sqrt{2\frac{n}{\omega^{2}}}(2\alpha +1)}\right|^{n} \left| \frac{1 +\sqrt{1+2\frac{n}{\omega^{2}}\left(4\alpha^{2}-1\right)}}{\sqrt{2\frac{n}{\omega^{2}} \left( 4\alpha^{2}-1 \right) }} \right|^{1/2} \\
  & & \mbox{}\times \exp \!\left(- \frac{\omega^{2}}{2(2\alpha + 1)} + \frac{n}{1+\sqrt{1 +2\frac{n}{\omega^{2}}\left(4\alpha^{2}-1\right)}}+  \frac{\alpha x_{0}^{2}}{4\alpha^{2}-1} \left[ 1 - \frac{2 \alpha}{\sqrt{1+2\frac{n}{\omega^{2}}\left(4\alpha^{2}-1\right)}}\right]    \right.  \\
  & & \left. \quad \mbox{} -\frac{\left[-1 + 2\frac{n}{\omega^{2}} + 8\frac{n}{\omega^{2}}\alpha^{2}\left(x_{0}^{2}-1\right)\right]^{2}}{16 n \left[1+ 2\frac{n}{\omega^{2}}\left(4\alpha^{2}-1\right)\right]^{5/2}} \right)\\
  & & \mbox{}\times \left[ \cos \!\left( \frac{n \pi}{2} - \frac{2\alpha x_{0} \omega \left(-2\alpha + \sqrt{1 + 2\frac{n}{\omega^{2}}\left(4\alpha^{2}-1\right)}\right)}{4\alpha^{2}-1}+ \frac{x_{0} \alpha \left[-3 + 6\frac{n}{\omega^{2}} + 8 \frac{n}{\omega^{2}}\alpha^{2}\left(x_{0}^{2} - 3\right)\right]}{3\omega\left[1 + 2 \frac{n}{\omega^{2}}\left(4\alpha^{2}-1\right)\right]^{3/2}}\right) \right.\\
  & & \left. \quad \mbox{}- \frac{\alpha	x_{0}}{\omega \left[ 1 + 2\frac{n}{\omega^{2}}\left(4\alpha	^{2} -1\right) \right]} \sin\! \left( \frac{n \pi}{2} - \frac{2\alpha x_{0} \omega \left(-2\alpha + \sqrt{1 + 2\frac{n}{\omega^{2}}\left(4\alpha^{2}-1\right)}\right)}{4\alpha^{2}-1} \right. \right. \\
  & & \left. \left. \quad \quad \mbox{}+ \frac{x_{0} \alpha \left[-3 + 6\frac{n}{\omega^{2}} + 8 \frac{n}{\omega^{2}}\alpha^{2}\left(x_{0}^{2} - 3\right)\right]}{3\omega\left[1 + 2 \frac{n}{\omega^{2}}\left(4\alpha^{2}-1\right)\right]^{3/2}} \right)\!\right]\!. 
\end{Eqnarray*}

\subsection{$\alpha < 1/2$: asymptotics} \label{sec:hermiteworking2}

Take $\zeta = (Y - \ii X)/\nu$ and let our contour of integration pass through the root $w_{+}$,
\begin{Eqnarray*}
  w_{+}\phi''(w_{+}) &=& 2\ii \sqrt{1-\zeta^{2}} =\frac{2\ii}{\nu}\sqrt{X^{2}+2\ii XY - Y^{2} + \nu^{2}}, \\
  \phi(w_{+}) &=& w_{+}^{2} - 2\zeta w_{+} + \frac{1}{2}\log(w_{+}),\\
  &=&\mbox{}-\frac{1}{4} - \frac{X+\ii Y}{2\left[X + \ii Y + \sqrt{X^{2} + \ii 2 X Y - Y^{2} + \nu^{2}}\right]}, \\
  &&\mbox{}-\frac{1}{2}\log (-2\ii) - \frac{1}{2}\log\left| \frac{X + \ii Y+ \sqrt{X^{2} + \ii 2 X Y - Y^{2} + \nu^{2}}}{\nu}\right|,
\end{Eqnarray*}
where
\begin{displaymath}
	w_{+} = \frac{1}{2}\left(\zeta + \ii \sqrt{1 - \zeta^{2}}\right), \qquad \nu = \sqrt{2 n +1}.
\end{displaymath}

\begin{figure}[htb]
  \begin{center}
      \includegraphics[width=.45\textwidth]{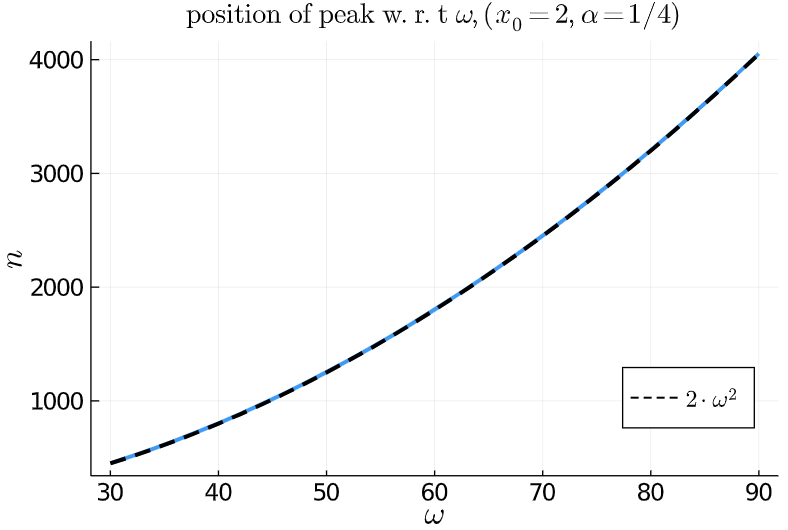}
      \includegraphics[width=.45\textwidth]{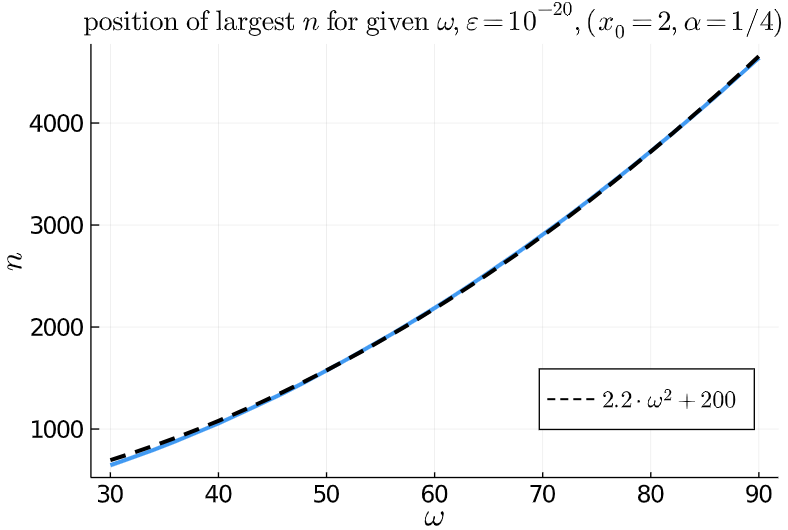}
    \end{center}
    \caption{Numerical experiments to determine the relationship between $n$ and $\omega$. The first figure is a plot of the position of the peak, $n$, as $\omega$ varies. The relationship follows $n = 2 \omega^{2}$. The second figure is the experiment where we found the largest of coefficients such that $|a_{n}| > 10^{-20}$. This relationship follows $n = 2.2 \omega^{2} +200$. The two experiments are consistent with a $\sqrt{n} \propto \omega$ relationship. \label{fig:numexpforalessthan12}}
  \end{figure}

Just like in $\alpha > 1/2$, numerical experiments show that the relationship between $n$ and $\omega$ follows
\begin{equation*}
	n \propto \omega^{2}.
\end{equation*}
See figure \ref{fig:numexpforalessthan12}.

Here we return to section \ref{sec:alessthan12} where we observe that for $\alpha < 1/2$, our asymptotics require us to consider $\cw \neq \cw^{*}$, where
\begin{Eqnarray}
	\cw^{*} = \frac{1}{2\left(1-4\alpha^{2}\right)}. \label{eq:cstar}
\end{Eqnarray}

 \subsubsection{$\cw \neq \cw^{*}$: full derivation} \label{sec:neqcstar}
 
Commencing from eqn \eqref{eq:AA1a}, substituting $n = \cw \omega^{2}$ and expanding about $\omega$, we get
\begin{Eqnarray*}
	\nu^{-1/2} (X^{2}+2\ii XY - Y^{2} + \nu^{2})^{-1/4} & = & \nu^{-1/2} \!\left[1 + \omega^{2} \left(2 \cw - \frac{1}{|1-4\alpha^{2}|} \right) + \frac{4\alpha^{2}x_{0}^{2}}{|1-4\alpha^{2}|} + \frac{4\ii \alpha x_{0} \omega}{|1-4\alpha^{2}|} \right]^{-1/4} \nonumber\\
		&\sim&\frac{|1-4\alpha^{2}|^{1/4}}{ \left(2\cw|1-4\alpha^{2}| -1\right)^{1/4}} \frac{1}{(2\cw)^{1/4} \omega} \\
		& & \times \left[1 - \frac{\ii \alpha x_{0}}{\omega \left(2 \cw |1 - 4\alpha^{2}| -1\right)} + \mathcal{O}\left(\omega^{-2}\right)\right]\!.
\end{Eqnarray*}
We substitute this back into \eqref{eq:AA1a}. 
Now, moving on to the asymptotics of $B$, expanding about $\omega$ having substituted $n=\cw\omega^{2}$,
\begin{Eqnarray*}
	\textnormal{Eqn } \eqref{eq:BB1a} &\sim& \frac{\cw \omega^{2} \left(1 + \ii \sqrt{2\cw |1- 4\alpha^{2}| -1}\right)}{2 \cw |1- 4\alpha^{2}|} + \frac{2\alpha x_{0} \omega\left(\cw |1- 4\alpha^{2}| - 1s\right)}{|1-4\alpha^{2}| \sqrt{2\cw |1-4\alpha^{2}|-1}|} - \frac{2 \ii \alpha x_{0} \omega }{|1-4\alpha^{2}|}\\
		& &\mbox{}- \frac{1}{4|1-4\alpha^{2}|\left(2\cw| 1 - 4\alpha^{2}| -1 \right)^{3/2}} \left[8\ii \alpha^{2}x_{0}^{2} + \ii |1 - 4\alpha^{2}| - 24 \ii \cw \alpha^{2}x_{0}^{2}|1-4\alpha^{2}| \right. \\
		& & \left. \mbox{}-2\ii \cw |1 - 4\alpha^{2}|^{2} + 8 x_{0}^{2} \alpha^{2} \left(2 \cw |1-4 \alpha^{2}| -1\right)^{3/2} \right] \\
		& & \mbox{}+ \frac{\alpha x_{0} \cw \left|1-4\alpha^{2}\right| \left(8 \alpha^{2}x_{0}^{2}\cw + 2\cw |4\alpha^{2} -1| -1\right)}{\omega \left(2\cw |1-4\alpha^{2}|-1\right)^{5/2}} + \mathcal{O}\left(\omega^{-2}\right), \\
		\textnormal{Eqn } \eqref{eq:BB1b} &\sim & \left(\frac{1}{2} + \cw \omega^{2}\right)\log \left| \frac{ \left(\ii + \sqrt{2\cw |1- 4\alpha^{2}|-1}\right)}{\sqrt{2 \cw \left|1-4\alpha^{2} \right|} }\right|+ \frac{2 \cw \alpha x_{0} \omega}{\sqrt{2\cw |4\alpha^{2}-1| -1}} \\
		& & \mbox{}- \frac{\ii \left(8 \alpha^{2}x_{0}^{2}\cw + 2\cw |1-4\alpha^{2}| -1\right)}{4 \left(2\cw |1-4\alpha^{2}| -1\right)^{3/2}} \\
		& & \mbox{}+ \frac{\alpha x_{0} \left[ 3 - 8 \alpha^{2}x_{0}^{2} \cw - \cw \left(9 + 8 \alpha^{2}x_{0}^{2}\cw\right)|1-4\alpha^{2}| + 6 \cw^{2} |1-4\alpha^{2}|^{2}\right] }{3\omega \left(2\cw |1-4\alpha^{2} | -1\right)^{5/2}} \ + \mathcal{O}\!\left(\omega^{-2}\right)\!.
\end{Eqnarray*} 
Combining the two results,
\begin{Eqnarray*}
	\textnormal{Eqn } \eqref{eq:BB1a}  + \textnormal{Eqn } \eqref{eq:BB1b} &\sim& \left(\cw \omega^{2} + \frac{1}{2}\right) \log \left(\frac{\ii + \sqrt{2\cw |1 - 4\alpha^{2}| -1}}{\sqrt{2 \cw |1 - 4\alpha^{2}|}} \right) \\
	& & \mbox{}+ \frac{ \cw \omega^{2}}{2\cw |1 - 4\alpha^{2}|} \left(1 + \ii \sqrt{2\cw |1 - 4\alpha^{2}| -1}\right) \\
	& &\mbox{}+ \frac{2 \alpha x_{0} \omega }{|1-4\alpha^{2}|}\left(-\ii + \sqrt{2 \cw |1-4\alpha^{2}| -1}\right) \\
	& & \mbox{}+ \frac{2 \alpha^{2} x_{0}^{2} \left(\ii -\sqrt{2 \cw |1 - 4\alpha^{2}| -1} \right)}{|1 - 4\alpha^{2}|\sqrt{2 \cw |1-4\alpha^{2}| -1}} \\
	& & \mbox{}+ \frac{\alpha x_{0} \left(8\cw \alpha^{2}x_{0}^{2} + 6 \cw |1 - 4\alpha^{2}| -3\right)}{3 \omega \left(2 \cw |1 - 4\alpha^{2}| -1\right)^{3/2}} + \mathcal{O}\!\left(\omega^{-2}\right)\!,
\end{Eqnarray*}
which we substitute back into $B$. Bringing all the terms together, 
\begin{Eqnarray*}
	a_{n}(\omega )= \textnormal{Re}(AB) &\sim& \textnormal{Re} \left(-\frac{1}{\sqrt{\omega(1+2\alpha)}}\left| \frac{2\alpha -1}{2 \alpha + 1}\right|^{n/2} \left(\frac{|1 - 4\alpha^{2}|}{2 \frac{n}{\omega^{2}} |1-4\alpha^{2}| -1}\right)^{\!1/4}  \right. \\
	& &\mbox{}\times  \left[1 - \frac{\ii \alpha x_{0}}{\omega \left(2\frac{n}{\omega^{2}} |1 - 4\alpha^{2}| -1\right)} + \mathcal{O}\! \left(\omega^{-2}\right)\right]  \left( \frac{1+ \sqrt{1-2\frac{n}{\omega^{2}} |1-4\alpha^{2} |}}{\sqrt{2 \frac{n}{\omega^{2}} |1-4\alpha^{2}|}}\right)^{\!n+1/2}\\
		& & \mbox{}\times \exp\left\{\frac{\alpha}{|1-4\alpha^{2}|}\left( \omega^{2}- x_{0}^{2} +2 x_{0} \omega \sqrt{2 \frac{n}{\omega^{2}}|1-4\alpha^{2}| -1} \right) \right. \\
		& & \mbox{}+ \frac{\alpha x_{0} \left(8 \alpha^{2} x_{0}^{2}\frac{n}{\omega^{2}} + 6 \frac{n}{\omega^{2}} |1-4\alpha^{2}| - 3\right)}{3\omega \left(2\frac{n}{\omega^{2}}|1-4\alpha^{2}|-1\right)^{3/2}} \\
		& & \mbox{}\times \ii \left[\left(n + \frac{1}{2} \right)\frac{\pi}{2} - \frac{ 4 \alpha^{2} x_{0} \omega}{|1-4\alpha^{2}|} + \frac{\omega^{2}}{2 |1-4\alpha^{2}|}\sqrt{2 \frac{n}{\omega^{2}} |1-4\alpha^{2}| -1}  \right. \\
		& & \left.\left. \left. \mbox{}+ \frac{2 \alpha^{2} x_{0}^{2}}{|1-4\alpha^{2}|\sqrt{2\frac{n}{\omega^{2}} |1-4\alpha^{2}|-1}}\right]\right\} \!\right)\!.
\end{Eqnarray*}

\section{Malmquist--Takenaka}

\subsection{Determining $a_{2k}$'s} \label{sec:a2ks}

The final step is to determine the coefficients $a^{(*)}_{2k}$. To do this, we first need to find the coefficients of the expansion near the saddle points,
\begin{equation}\label{eq:note14}
z=z^{*}+\sum_{k=0}^\infty b_k^{(*)} t^k.
\end{equation}
They can be obtained from the transformations in \eqref{eq:note10},  where we  rewrite the transformations as values of $z$ near the saddle points
\begin{equation*}
\sum_{k=2}^\infty \frac{g^{(k)}(z^{*})}{k!} (z-z^{*})^k=\tfrac12 t^2,
\end{equation*}
or,  taking the square root,
\begin{equation}\label{eq:note16}
t=(z-z^{*})\left[g^{(2)}(z^{*})+2\sum_{k=3}^\infty \frac{g^{(k)}(z^{*})}{k!} (z-z^{*})^{k-2}\right]^{\!1/2}\!.
\end{equation}
This equation can be inverted by using standard inversion methods. Substituting the expansions given in \eqref{eq:note14} and collecting equal powers of $t$, we can find the coefficients $b_k^{(*)}$ of \eqref{eq:note14}, the first  being
\begin{equation*}
\tfrac12 g^{(2)}(z^{*}){b_1^{(*)} t}^2=\tfrac12 t^2\quad \Longrightarrow \quad b_1^{(*)} =\frac{1}{\sqrt{g^{(2)}(z^{*})}}.
\end{equation*}
Here and in \eqref{eq:note16} we have taken the positive square root, this comes from the expansion in \eqref{eq:note14}, in fact $z-z^{*}\sim b_1^{(*)} t$, and we have assumed conditions as explained below \eqref{eq:note12}.

More coefficients can be obtained by using computer algebra. The next ones are
\begin{Eqnarray*}
b_2^{(*)}&=&\displaystyle{-\frac{g^{(3)}(z^{*}) {b_1^{(*)}}^2}{6g^{(2)}(z^{*})}}=\displaystyle{-\frac{g^{(3)}(z^{*})}{6\left[g^{(2)}(z^{*}) \right]^{2}},}  \\
b_3^{(*)}&=&\displaystyle{-\frac{g^{(4)}(z^{*}){b_1^{(*)}}^4 + 12g^{(3)}(z^{*})z_2^{*}{b_1^{(*)}}^2+ 12g^{(2)}(z^{*})b_{2}^{(*)}}{24g^{(2)}(z^{*}) b_1^{(*)}},}\\
	&=& \displaystyle{-\frac{3 g^{(4)}(z^{*}) g^{(2)}(z^{*}) - 5 \left[ g^{(3)}(z^{*}) \right]^{2}}{72 \left[ g^{(2)}(z^{*})\right]^{7/2}}}.
\end{Eqnarray*}
To obtain $a_{2k}$'s we use \eqref{eq:note12} and \eqref{eq:note14}
\begin{Eqnarray*}
	\frac{\dd z}{\dd t} \frac{1}{1+ 2 \ii z} &=&\left( \sum_{k=0}^{\infty} k b_{k}^{(*)} t^{k-1} \right) \left[ \sum_{k=0}^{\infty} \frac{1}{k!}\left( \frac{\dd^{k}}{\dd z^{k}} \frac{1}{1+2 \ii z} \right)_{z= z*} (z-z*)^{k} \right] = \sum_{k=0}^{\infty}a_{k}^{(*)}t^{k}
\end{Eqnarray*}
Equating powers of $t$ we obtain the first coefficients of the expansion in \eqref{eq:note13}, which are
\begin{Eqnarray*}
a_0^{(*)}&=&\frac{b_1^{(*)}}{1+2\ii z^{*}} = \frac{1}{(1 + 2 \ii z^{*})\sqrt{g''(z^{*})}} ,\\
a_2^{(*)}&=&\frac{-12{z^{*}}^2b_3^{(*)} + 12z^{*}b_1^{(*)}b_2^{(*)}- 4{b_1^{(*)}}^3 + 12\ii z^{*}b_3^{(*)} - 6\ii b_1^{(*)}b_2^{(*)}+ 3b_3^{(*)}}{(1 + 2\ii z^{*})^3},  \\
	&=& \frac{1}{2(\ii - 2 z^{*}) {g''}^{7/2}(z^{*})} \left\{ 8 {g''}^2(z^{*})+ g^{(2)}(z^{*})\!\left[ 2 g^{(3)}(z^{*})(2z^{*}- \ii) -\frac{g^{(4)}(z^{*})}{4}(\ii - 2 z^{*})^{2} \right] \right.\\
	& &\mbox{} \left. +\frac{5}{12} {g^{(3)}}^2\!(z^{*}) (1 - 2 \ii z^{*})^{2} \right\}\!.
\end{Eqnarray*}

\subsection{Bounding $|b_{n}(\omega)|$ for $\omega \geq 0$ and $n\geq0$} \label{sec:boundbnw}
The result obtained in \eqref{eq:mtestfull} is not obvious from the asymptotic expansions in section \ref{sec:mtasyexpan} as they require further manipulation.  In this section, we walk through the steps taken to bound \eqref{eq:mt_bn_firstapprox}.

Applying the triangle inequality to \eqref{eq:mt_bn_firstapprox} results in
\begin{Eqnarray} \label{eq:bnboundgen}
\left| b_n(\omega) \right| \lesssim
\left| \frac{\ee^{-\omega g(z_3)}}{(1 + 2 \ii z_{3})\sqrt{\omega g''(z_{3})}} \right|
+ \left|  \frac{\ee^{-\omega g(z_2)}}{(1 + 2 \ii z_{2})\sqrt{\omega g''(z_{2})}} \right|.
\end{Eqnarray}

Following the same expansion steps as described in section \ref{sec:mtasyexpan}, we obtain asymptotics that would help us bound $|b_{n}(\omega)|$,
\begin{Eqnarray*}
	\ee^{-\omega g(z_{2})} &\sim& \exp \left( -\frac{\alpha}{\omega} \left( \sqrt{\cw - 1} - 2x_{0}\right)^{2}\right)\left[\cos \left( \frac{\omega \sqrt{4\cw -1}}{2} - 2\cw \omega \tan^{-1} \left( \sqrt{4\cw -1}\right) + \mathcal{O}\left(\omega^{-1}\right)\right) \right. \\
	& & \quad \left. + \ii \sin \left( \frac{\omega \sqrt{4\cw -1}}{2} - 2\cw \omega \tan^{-1} \left( \sqrt{4\cw -1} \right) + \mathcal{O}\left(\omega^{-1} \right)\right) \right], \\
		\ee^{-\omega g(z_{3})} &\sim& \exp \left( -\frac{\alpha}{\omega} \left( \sqrt{\cw - 1} + 2x_{0}\right)^{2}\right)\left[\cos \left(- \frac{\omega \sqrt{4\cw -1}}{2} + 2\cw \omega \tan^{-1} \left( \sqrt{4\cw -1} \right)+ \mathcal{O}\left(\omega^{-1}\right) \right) \right. \\
	& & \quad \left. + \ii \sin \left(- \frac{\omega \sqrt{4\cw -1}}{2} +2\cw \omega \tan^{-1} \left( \sqrt{4\cw -1} \right)+ \mathcal{O}\left(\omega^{-1}\right)\right)  \right], \\
	\left(1 + 2\ii z_{2}\right)^{-1} &\sim& \frac{1}{4 \cw} \left\{ 1 - \ii \sqrt{4\cw -1} - \frac{\alpha}{\omega \sqrt{4\cw -1}} \left[(1-2\cw)( \sqrt{4\cw -1}-2x_{0})  \right. \right. \\
	& & \left. \left. \quad + \ii \left(2x_{0}\sqrt{4\cw -1} - 4\cw +1 \right) \right] + \mathcal{O}\left(\omega^{-1}\right)\right\}, \\
	\left(1 + 2\ii z_{3}\right)^{-1} &\sim& \frac{1}{4 \cw} \left\{ 1 + \ii \sqrt{4\cw -1} - \frac{\alpha}{\omega \sqrt{4\cw -1}} \left[(1-2\cw)(2x_{0} + \sqrt{4\cw -1}) \right. \right. \\
	& & \left. \left. \quad + \ii \left(2x_{0}\sqrt{4\cw -1} + 4\cw -1 \right) \right] + \mathcal{O}\left(\omega^{-1}\right)\right\}, \\
	\left( \omega g''(z_{2})\right)^{-1/2} &\sim&\frac{\left(1 + \ii\right) \sqrt{\cw}}{\sqrt{2\omega}\left(4\cw -1\right)^{1/4}} \left\{ 1 - \frac{\ii \alpha}{\omega \left(4 \cw -1\right)^{3/2}} \left[16 \left(\cw - \frac{1}{4}\right)^{2}-6x_{0}\left(\cw - \frac{1}{3} \right)\sqrt{4\cw -1}  \right] + \mathcal{O}\left(\omega^{-2} \right)\right\},\\
	\left( \omega g''(z_{3})\right)^{-1/2} &\sim&\frac{\left(-1 + \ii\right)\sqrt{\cw}}{\sqrt{2 \omega}\left(4\cw -1\right)^{1/4}}\left\{ 1 + \frac{\ii \alpha}{\omega \left(4 \cw -1\right)^{3/2}} \left[6x_{0}\left(\cw - \frac{1}{3} \right)\sqrt{4\cw -1} + 16 \left(\cw - \frac{1}{4}\right)^{2} \right] + \mathcal{O}\left(\omega^{-2} \right)\right\}.
\end{Eqnarray*}
Substituting the above terms into \eqref{eq:bnboundgen}, separating each term into real and imaginary parts,  and using the definition of the $\LL_{2}$ norm  $|x+ \ii y| = \sqrt{x^{2} + y^{2}}$ we obtain the following bounds 
\begin{Eqnarray*}
	\left| \frac{\ee^{-\omega g(z_2)}}{(1 + 2 \ii z_{2})\sqrt{\omega g''(z_{3})}} \right| &\leq& 
	\frac{\exp\left(-\frac{\alpha}{4} \left(\sqrt{4\cw -1} -2x_{0}\right)^{2}\right)}{2 \sqrt{\omega} \left(4\cw -1\right)^{1/4}} \left[ 1 - \frac{\alpha\left(\sqrt{4\cw -1} - 2x_{0}\right)}{\omega\sqrt{4\cw -1}} + \mathcal{O}\left(\omega^{-2}\right)\right]^{1/2},\\
	\left| \frac{\ee^{-\omega g(z_3)}}{(1 + 2 \ii z_{3})\sqrt{\omega g''(z_{3})}} \right| &\leq&
	\frac{\exp\left(-\frac{\alpha}{4} \left(\sqrt{4\cw -1} +2x_{0}\right)^{2}\right)}{2 \sqrt{\omega} \left(4\cw -1\right)^{1/4}} \left[ 1 - \frac{\alpha\left(\sqrt{4\cw -1} +2x_{0}\right)}{\omega\sqrt{4\cw -1}} + \mathcal{O}\left(\omega^{-2}\right)\right]^{1/2}.
\end{Eqnarray*}
From this point we can see that we obtain the bound in \eqref{eq:mtestfull}.
\end{document}